\documentclass[10pt,a4paper]{smfart}
\usepackage[latin1]{inputenc}
\usepackage{amsmath}
\usepackage{pstricks}
\usepackage{amsfonts}
\usepackage{amssymb}
\usepackage{amsthm}
\usepackage{lmodern}
\usepackage[all]{xy}
\usepackage{mathrsfs}
\usepackage[english,francais]{babel}
\newtheorem{defi}{D\'efinition}[section]
\newtheorem{thm}[defi]{Th\'eor\`eme}
\newtheorem{prop}[defi]{Proposition}
\newtheorem{lem}[defi]{Lemme}

\newtheorem{coro}[defi]{Corollaire}

\newtheorem{rem}[defi]{Remarque}
\newtheorem{ex}[defi]{Exemple}
\newtheorem{propr}[defi]{Propri\'et\'e}

\author{Jean-Christophe SAN SATURNINO}
\title[Uniformisation locale des sch\'emas quasi-excellents 
]{Uniformisation locale des sch\'emas quasi-excellents de caract\'eristique nulle}
\address{Universit\'e Toulouse III Paul Sabatier\\
Institut de Math\'ematiques de Toulouse\\ 
118 route de Narbonne\\
31062 Toulouse cedex 9 (France)}
\email{san@math.ups-tlse.fr}
\urladdr{http://www.math.univ-toulouse.fr/$\sim$san/}

\begin{document}
\begin{abstract}
Nous d\'emontrons un th\'eor\`eme d'uniformisation locale plong\'ee pour une valuation centr\'ee en un point d'un sch\'ema quasi-excellent de caract\'eristique nulle. La preuve se r\'eduit au cas des valuations de rang $1$ et consiste \`a d\'esingulariser l'id\'eal form\'e des \'el\'ements de valeur infinie et \`a monomialiser les polyn\^omes-cl\'es. On d\'emontre alors un th\'eor\`eme de monomialisation valable en toute caract\'eristique sous certaines conditions, notamment celle de ne pas avoir de polyn\^ome-cl\'e limite, fait qui se produit toujours en caract\'eristique nulle.
\end{abstract}
\begin{altabstract}
We prove an embedded local uniformization theroem for a valuation centered on a point of a quasi-excellent scheme of characteristic zero. The proof reduces to valuations of rank $1$ and consists in desingularizing the ideal formed by the  elements of infinite value and monomializing the key polynomials. We then prove a monomialization theorem valid in all characteristic under certain conditions, including that of non-existence of limit key polynomials, a condition that is always satified in characteristic zero.
\end{altabstract}
\maketitle
\setcounter{section}{-1}
\section{Introduction}
\indent Les premiers r\'esultats en r\'esolution des singularit\'es sont \`a attribuer \`a Newton au XVII$^{\textup{\`eme}}$ si\`ecle et \`a Puiseux au XIX$^{\textup{\`eme}}$ si\`ecle. Leurs r\'esultats permettent de r\'esoudre les singularit\'es des courbes d\'efinies sur $\mathbb C$. En 1939, Zariski (voir \cite{zariskiunifsurface}) propose une nouvelle m\'ethode pour r\'esoudre les singularit\'es d'une surface d\'efinie sur un corps de caract\'eristique nulle: on r\'esout le probl\`eme localement le long d'une valuation puis on recolle au niveau de la vari\'et\'e de Riemann-Zariski. Le recollement n'est actuellement possible que jusqu'\`a la dimension $3$ en toute caract\'eristique (voir par exemple \cite{piltantrecole}). La r\'esolution locale le long d'une valuation, ou uniformisation locale, n'est possible, en caract\'eristique positive ou mixte, \'egalement que jusqu'\`a la dimension $3$ (voir \cite{cosspilt1}, \cite{cosspilt2}, \cite{cossartpiltantmixtehal} et \cite{cossartpiltantmixtehal2}). En caract\'eristique nulle, la r\'esolution des singularit\'es a \'et\'e d\'emontr\'ee par Hironaka en 1964 (\cite{hiro}) pour des vari\'et\'es de dimension quelconque. Ce r\'esultat a par la suite \'et\'e red\'emontr\'e par Villamayor en 1989 (\cite{vill89}), Bierstone et Milman en 1990 (\cite{birsmilm}), Encinas et Villamayor en 2001 (\cite{encvilla}), Encinas et Hauser en 2002 (\cite{encinashauser}), W{\l}odarczyk en 2005 (\cite{wlod}) et Temkin en 2008 (\cite{tem}).
\\\indent Ces derni\`eres ann\'ees, une nouvelle approche a \'et\'e propos\'ee par Spivakovsky (\cite{spiva}) et Teissier (\cite{teissiertorique}) pour r\'esoudre le probl\`eme de la r\'esolution des singularit\'es en caract\'eristique positive et mixte via la m\'ethode de Zariski. La premi\`ere \'etape \'etant de d\'emontrer l'uniformisation locale d'une valuation, ils se sont int\'eress\'es \`a l'alg\`ebre gradu\'ee qui lui est naturellement associ\'ee ainsi qu'\`a ses g\'en\'erateurs, appel\'es \textit{polyn\^omes-cl\'es}.
\\\indent Dans cet article, nous adaptons l'approche de Spivakovsky (\cite{spiva}) pour l'uniformisation locale plong\'ee des anneaux quasi-excellents \'equicaract\'eristiques au cas o\`u le corps r\'esiduel est de caract\'eristique nulle. Ce r\'esultat est d\'ej\`a bien connu, l'int\'er\^et de notre d\'emonstration est qu'elle permet de d\'ecrire \`a l'avance toutes les coordonn\'ees de tous les points infiniment proches et ceci pour tous les \'eclatements interm\'ediaires. La m\'ethode est la suivante: par \cite{spivanova}, on sait qu'il suffit d'obtenir l'uniformisation locale pour des valuations de rang $1$ centr\'ees sur un anneau $S$ local noeth\'erien int\`egre. On consid\`ere ensuite l'id\'eal $\overline H$ compos\'e de tous les \'el\'ements de valuation infinie dans le compl\'et\'e de $S$, not\'e $\widehat{S}$. Cet id\'eal est premier d'apr\`es \cite{ideimp}, c'est l'\textit{id\'eal premier implicite}. Si $S$ est quasi-excellent, c'est-\`a-dire, si le morphisme de compl\'etion est r\'egulier, alors $\widehat{S}_{\overline H}$ est r\'egulier (Th\'eor\`eme \ref{RHestreg}). Pour conclure, il suffit de montrer qu'il existe une suite d'\'eclatements locaux qui rendent $\widehat{S}/\overline H$ r\'egulier et que tout \'el\'ement de cet anneau s'\'ecrit comme le produit d'un mon\^ome par une unit\'e (Lemme \ref{lemmetechniquespiva} et Th\'eor\`eme \ref{uniflocalerang1car0}). Par le th\'eor\`eme de structure de Cohen, on sait qu'il existe un anneau local r\'egulier complet $R$ et un morphisme surjectif $R\twoheadrightarrow \widehat{S}/\overline H$. Ce morphisme induit un isomorphisme entre $R/H$ et $\widehat{S}/\overline H$, o\`u $H$ est le noyau du morphisme. Pour obtenir notre r\'esultat, il ne reste plus qu'\`a l'obtenir sur $R/H$ (Th\'eor\`eme \ref{thmprelimcar0}). 
\\\indent La strat\'egie va consister \`a faire d\'ecro\^itre la dimension de plongement de $R/H$ si $H$ est non-nul, sinon \`a monomialiser les \'el\'ements de $R$. 
\\\indent Dans le cas \'equicaract\'eristique, on sait que $R$ est un anneau de s\'eries formelles, on peut alors utiliser la th\'eorie des polyn\^omes-cl\'es. Par r\'ecurrence sur la dimension de plongement, on montre que l'id\'eal $H$ est, \`a une suite d'\'eclatements locaux pr\`es, principal et engendr\'e par un polyn\^ome unitaire. La suite d'\'eclatements choisie n'est pas quelconque, les param\`etres r\'eguliers de l'anneau d'arriv\'ee sont des polyn\^omes-cl\'es et d\`es qu'un \'el\'ement s'est transform\'e en mon\^ome, la suite d'\'eclatements le conserve sous cette forme. 
\\\indent Comme tout polyn\^ome de $R$ peut s'\'ecrire de mani\`ere unique comme une somme finie de polyn\^omes-cl\'es, il suffit de monomialiser les polyn\^omes-cl\'es. S'il n'y a pas de premier polyn\^ome-cl\'e limite, une simple r\'ecurrence nous fournit le r\'esultat (Proposition \ref{eclatpolyclecar0}). Cette situation correspond au cas o\`u, apr\`es un nombre fini d'\'eclatements, on fait un nombre fini de translations pour obtenir un r\'esultat de monomialisation (en termes de polyn\^omes-cl\'es et de valuations, cette situation est celle o\`u, apr\`es un nombre fini $i_0$ d'\'etapes, la valuation de d\'epart est \'egale \`a la valuation monomiale correspondante au polyn\^ome-cl\'e de l'\'etape $i_0$). En terme de d\'efaut d'une extension, cette situation correspond au cas o\`u il n'y a pas de d\'efaut; dans un article en pr\'eparation (\cite{jcdefaut}) on montrera que le d\'efaut peut \^etre compris tr\`es pr\'ecis\'ement en terme de degr\'e du premier polyn\^ome-cl\'e limite.
\\\indent Or cette situation se produit toujours en \'equicaract\'eristique nulle: il n'existe pas de polyn\^ome-cl\'e limite pour des valuations de rang $1$ (Corollaire \ref{polycleencar0}). 
\\\indent Dans le cas o\`u il existe un premier polyn\^ome-cl\'e limite, c'est-\`a-dire lorsque l'on fait un nombre infini de translations, autrement dit, dans le cas o\`u l'extension consid\'er\'ee poss\`ede un d\'efaut, on ne sait pas conclure. Dans \cite{jcthese}, on a montr\'e qu'il suffit de monomialiser le premier polyn\^ome-cl\'e limite qui peut \^etre vu comme un polyn\^ome d'Artin-Schreier g\'en\'eralis\'e.
\\\indent En caract\'eristique mixte, la situation est semblable, $R$ est un anneau de s\'eries formelles \`a coefficients sur un anneau de valuation discr\`ete, ou un quotient d'un anneau de s\'eries formelles de ce type. On peut \'egalement conclure lorsqu'il n'y a pas de premier polyn\^ome-cl\'e limite en utilisant le cas \'equicaract\'eristique mais il nous faut supposer que la valuation de $p$, la caract\'eristique du corps r\'esiduel de $R$, n'est pas divisible par $p$ dans le groupe des valeurs.
\\ \\\indent Pour commencer, nous rappelons les d\'efinitions de centre d'une valuation ainsi que celles des diff\'erentes alg\`ebres gradu\'ees utilis\'ees.
\\\indent Dans la deuxi\`eme partie, nous redonnons la notion de quasi-excellence pour un anneau ainsi que pour un sch\'ema. Nous d\'efinissons \'egalement les diff\'erentes propri\'et\'es d'uniformisation locale et d'uniformisation locale plong\'ee que nous allons d\'emontrer dans le cas de caract\'eristique nulle.
\\\indent Par la suite, apr\`es avoir rappel\'e la notion de polyn\^ome-cl\'e, nous montrons qu'il n'y a pas de polyn\^omes-cl\'es limite en caract\'eristique nulle pour des valuations de rang $1$.
\\\indent Dans la quatri\`eme partie, nous d\'eveloppons les outils n\'ecessaires pour les diff\'erentes preuves des th\'eor\`emes de monomialisation et d'uniformisation locale. Nous donnons tout d'abord la d\'efinition d'\'eclatement local encadr\'e qui va imposer un syst\`eme de g\'en\'erateurs de l'id\'eal maximal; ce type d'\'eclatement va conserver la propri\'et\'e d'\^etre un produit d'un mon\^ome par une unit\'e. Par la suite, on en construit un explicitement ayant la propri\'et\'e de transformer les polyn\^omes-cl\'es en param\`etres r\'eguliers. Nous rappelons ensuite les r\'esultats essentiels sur l'id\'eal premier implicite. Nous terminons cette partie par la monomialisation d'\'el\'ements non-d\'eg\'en\'er\'es (c'est-\`a-dire dont leur valuation est \'egale \`a la valuation monomiale), qui est un cas particulier du jeu d'Hironaka (voir \cite{jeuhiro} et \cite{spivajeuhiro}), ainsi que par un th\'eor\`eme d'uniformisation locale pour des hypersurfaces quasi-homog\`enes satisfaisant certaines propri\'et\'es.
\\\indent Nous d\'emontrons, dans la cinqui\`eme partie, un th\'eor\`eme de monomialisation dans le cas \'equicaract\'eristique. \`A chaque \'etape de l'algorithme, un \'eclatement local est suivi d'une compl\'etion. On d\'esingularise l'id\'eal $H$ et on monomialise les polyn\^omes-cl\'es dans le cas o\`u il n'y a pas de premier polyn\^ome-cl\'e limite.
\\\indent La sixi\`eme partie est identique \`a la pr\'ec\'edente sauf que l'on se place dans le cadre de la caract\'eristique mixte, avec l'hypoth\`ese suppl\'ementaire que la valuation de $p$, o\`u $p$ est la caract\'eristique du corps r\'esiduel de $R$, n'est pas divisible par $p$ dans le groupe des valeurs.
\\\indent Dans la septi\`eme partie, nous d\'emontrons \`a nouveau le m\^eme r\'esultat de monomialisation que dans les deux parties pr\'ec\'edentes mais sans compl\'eter apr\`es chaque \'eclatement. 
\\\indent Nous terminons par la d\'emonstration du r\'esultat principal d'uniformisation locale plong\'ee d'une valuation centr\'ee en un point d'un sch\'ema quasi-excellent dont le corps r\'esiduel de l'anneau local en ce point est de caract\'eristique nulle.
\\ \\\indent Je tiens \`a remercier Mark Spivakovsky qui a d\'evelopp\'e la plupart des outils de cet article pendant de longues ann\'ees et qui m'a permis de les appliquer dans le cadre de la caract\'eristique nulle.
\\ \\ \noindent\textbf{Notations.} Soit $\nu$ une valuation sur un corps $K$. Notons $R_\nu=\lbrace f\in K\:\vert\:\nu(f)\geqslant 0\rbrace$, c'est un anneau local d'id\'eal maximal $m_\nu=\lbrace f\in K\:\vert\:\nu(f)> 0\rbrace$. On note alors $k_\nu=R_\nu/m_\nu$ le corps r\'esiduel de $R_\nu$ ansi que $\Gamma_\nu=\nu(K^*)$.
\\Si $R$ est un anneau, on notera $car(R)$ sa caract\'eristique. Si $(R,\mathfrak m)$ est un anneau local on notera $\widehat{R}$ le compl\'et\'e $\mathfrak m$-adique de $R$.
\\Pour tout $P\in Spec(R)$, on note $\kappa(P)=R_P/PR_P$ le corps r\'esiduel de $R_P$.
\\Pour $\alpha\in\mathbb Z^n$ et $u=(u_1,...,u_n)$ un ensemble d'\'el\'ements de $R$, on note:
\[u^\alpha=u_1^{\alpha_1}...u_n^{\alpha_n}.\]
Pour $P,Q\in R\left[ X\right] $ avec $P=\sum\limits_{i=0}^{n}a_{i}Q^{i}$ et $a_{i}\in R[X]$ tels que le degr\'e de $a_{i}$ est strictement inf\'erieur \`a celui de $Q$, on note:
\[d_{Q}^{\:\circ}(P)=n.\]
Si $Q=X$, on notera plus simplement $d^{\:\circ}(P)$ au lieu de $d_{X}^{\:\circ}(P)$.
\\Enfin, si $R$ est un anneau int\`egre, on notera $Frac(R)$ son corps des fractions.

\section{Centre d'une valuation, alg\`ebres gradu\'ees associ\'ees}
\begin{defi}\label{centrevaldef}
Soient $R$ un anneau et $P$ un id\'eal premier. Une valuation de $R$ \textbf{centr\'ee en} $\boldsymbol P$ est la donn\'ee d'un id\'eal premier minimal $P_\infty$ de $R$ contenu dans $P$ et d'une valuation du corps des fractions de $R/P_\infty$ centr\'ee en $P/P_\infty$. L'id\'eal $P_\infty$ est alors le support de la valuation.
\\Si $R$ est un anneau local d'id\'eal maximal $\mathfrak{m}$, on dira que $\nu$ est \textbf{centr\'ee en} $\boldsymbol R$ pour dire que $\nu$ est centr\'ee en $\mathfrak{m}$.
\\Soit $X$ un sch\'ema int\`egre de corps des fonctions $K(X)$. Une valuation $\nu$ de $K(X)$ est \textbf{centr\'ee en un point} $\boldsymbol{\xi}$ de $X$ si $\nu$ est centr\'ee en $\mathcal{O}_{X,\xi}$. On dira alors que $\xi$ est le \textbf{centre} de $\nu$.
\end{defi}
\begin{defi}\label{alggrad}
Soient $R$ un anneau et $\nu:R\rightarrow \Gamma\cup\lbrace\infty\rbrace$ une valuation centr\'ee en un id\'eal premier de $R$. Pour tout $\alpha\in\nu(R\setminus\lbrace 0\rbrace)$, on d\'efinit les id\'eaux:
\[P_\alpha=\lbrace f\in R\:\vert\:\nu(f)\geqslant\alpha\rbrace;\]
\[P_{\alpha,+}=\lbrace f\in R\:\vert\:\nu(f)>\alpha\rbrace.\]
L'id\'eal $P_\alpha$ est appel\'e le $\boldsymbol{\nu}$\textbf{-id\'eal} de $R$ de valuation $\alpha$.
\\On d\'efinit alors \textbf{l'alg\`ebre gradu\'ee de} $\boldsymbol{R}$ \textbf{associ\'ee \`a} $\boldsymbol{\nu}$ par:
\[gr_{\nu}(R)=\bigoplus\limits_{\alpha\in\nu(R\setminus\lbrace 0\rbrace)}P_{\alpha}/P_{\alpha,+}.\]
L'alg\`ebre $gr_\nu(R)$ est un anneau int\`egre.
\\Pour $f\in R\setminus\lbrace 0\rbrace$, on d\'efinit son image dans $gr_\nu(R)$, not\'ee $in_\nu(f)$, comme \'etant l'image naturelle de $f$ dans $P_{\nu(f)}/P_{\nu(f),+}\subset gr_\nu(R)$; c'est un \'el\'ement homog\`ene de degr\'e $\nu(f)$.
\\Enfin, on d\'efinit une valuation naturelle sur $gr_\nu(R)$ de groupe des valeurs $\nu(R\setminus\lbrace 0\rbrace)$, not\'ee $ord$, par:
\[ord(f)=\min \alpha,\]
o\`u $f\in gr_\nu(R)$ s'\'ecrit comme une somme finie $f=\sum\limits_{\alpha\in\nu(R\setminus\lbrace 0\rbrace)}f_{\alpha}$, $f_\alpha\in P_\alpha/P_{\alpha,+}$.
\end{defi}
Si $R$ est un anneau local int\`egre, on d\'efinit une autre alg\`ebre gradu\'ee comme suit:
\begin{defi}
Soient $R$ un anneau local int\`egre, $K=Frac(R)$ et $\nu:K^\times\twoheadrightarrow \Gamma\cup\lbrace\infty\rbrace$ une valuation de $K$ centr\'ee en $R$. Pour tout $\alpha\in\Gamma$, on d\'efinit les $R_\nu$-sous-modules de $K$ suivants:
\[\boldsymbol{P_\alpha}=\lbrace f\in K\:\vert\:\nu(f)\geqslant\alpha\rbrace;\]
\[\boldsymbol{P_{\alpha,+}}=\lbrace f\in K\:\vert\:\nu(f)>\alpha\rbrace.\]
On d\'efinit alors \textbf{l'alg\`ebre gradu\'ee associ\'ee \`a} $\boldsymbol{\nu}$ par:
\[G_{\nu}=\bigoplus\limits_{\alpha\in\Gamma}\boldsymbol{P_{\alpha}}/\boldsymbol{P_{\alpha,+}}.\]
Pour $f\in K^\times$, on d\'efinit son image dans $G_\nu$, not\'ee $in_\nu(f)$, comme dans la D\'efinition \ref{alggrad}. 
\\Enfin, on d\'efinit une valuation naturelle sur $G_\nu$ de groupe des valeurs $\Gamma$, not\'ee $ord$, comme dans la D\'efinition \ref{alggrad}.
\end{defi}
\begin{rem}
\textup{On a l'injection naturelle:
\[gr_\nu(R)\hookrightarrow G_\nu.\]}
\end{rem}
\begin{defi}
Soit $G$ une alg\`ebre gradu\'ee n'ayant pas de diviseurs de z\'ero. On appelle \textbf{satur\'e de} $\boldsymbol G$ l'ag\`ebre gradu\'ee $G^*$ d\'efinie par:
\[G^{*}=\left\lbrace \left.\dfrac{f}{g}\:\right|\:f,g\in G,\:g\textit{ homog\`ene},\:g\neq 0\right\rbrace.\]
On dit que $G$ est \textbf{satur\'ee} si $G=G^*$.
\end{defi}
\begin{rem}
\textup{Pour toute alg\`ebre gradu\'ee $G$, on a:
\[G^*=\left(G^*\right)^*.\]
Autrement dit, $G^*$ est toujours satur\'ee.}
\end{rem}
\begin{ex}
\textup{Soit $\nu$ une valuation centr\'ee en un anneau local $R$. Alors:
\[G_\nu=\left(gr_\nu(R)\right)^*.\]
En particulier, $G_\nu$ est satur\'ee.}
\end{ex}
\section{Quasi-excellence, \'eclatements locaux et uniformisation locale}
Pour plus de clart\'e nous rappelons la notion de quasi-excellence ainsi que diff\'erentes notions d'uniformisation locale, que ce soit pour des sch\'emas ou pour des anneaux. L'uniformisation locale est la version locale de la r\'esolution des singularit\'es. R\'esoudre les singularit\'es d'un sch\'ema $X$ noeth\'erien irr\'eductible et r\'eduit revient \`a trouver un morphisme propre et birationnel $X'\rightarrow X$ tel que $X'$ soit r\'egulier. Ainsi, l'uniformisation locale d'une valuation $\nu$ de $K$, le corps des fractions d'un anneau local int\`egre $R$ o\`u est centr\'ee la valuation, revient \`a trouver un anneau $R'$ r\'egulier qui domine birationnellement $R$ et tel que $R'\subset R_\nu\subset K$.
\begin{defi}
Un anneau noeth\'erien $R$ est \textbf{quasi-excellent} si les deux conditions suivantes sont v\'erifi\'ees:
\begin{enumerate}
\item Pour tout $P\in Spec(R)$, le morphisme de compl\'etion $R_P\rightarrow\widehat{R_P}$ est r\'egulier;
\item Le lieu r\'egulier de toute $R$-alg\`ebre de type fini est ouvert.
\end{enumerate}
Un sch\'ema localement noeth\'erien est dit \textbf{quasi-excellent} s'il existe un recouvrement form\'e d'ouverts affines $(U_{\alpha})$, $U_\alpha=Spec(R_\alpha)$, tel que, pour tout $\alpha$, $R_\alpha$ soit un anneau quasi-excellent.
\end{defi}
\begin{rem}
\textup{\begin{enumerate}
\item Un anneau local est quasi-excellent si et seulement si la condition 1. est v\'erifi\'ee.
\item La notion de quasi-excellence est conserv\'ee par localisation, passage au quotient et passage aux alg\`ebres de type fini.
\item Un corps est quasi-excellent.
\item les anneaux de s\'eries formelles sur un anneau de Cohen sont des anneaux quasi-excellents.
\end{enumerate}}
\end{rem}
\begin{defi}
Soient $X$ un sch\'ema noeth\'erien et $Y$ un sous-sch\'ema de $X$. Soit $\mathcal{I}_Y$ le faisceau d'id\'eaux d\'efinissant $Y$ dans $X$.
\\On dit que $X$ est \textbf{normalement plat} le long de $Y$ si, pour tout point $\xi\in Y$, $\bigoplus\limits_{n\geqslant 0}\mathcal{I}_{Y,\xi}^n/\mathcal{I}_{Y,\xi}^{n+1}$ est un $\mathcal{O}_{Y,\xi}$-module libre.
\end{defi}
\begin{propr}\textup{\textbf{(d'uniformisation locale des sch\'emas).}}
Soit $S$ un sch\'ema noeth\'erien (non n\'ecessairement int\`egre). Soient $X$ une composante irr\'eductible de $S_{red}$ et $\nu$ une valuation de $K(X)$
centr\'ee en un point $\xi\in X$. Il existe alors un \'eclatement $\pi: S'\rightarrow S$ le long d'un sous-sch\'ema de $S$, ne contenant aucune composante irr\'eductible de $S_{red}$ et ayant la propri\'et\'e suivante : \\Soient $X'$ le transform\'e strict de $X$ par $\pi$ et $\xi'$ le centre de $\nu$ sur $X'$, alors $\xi'$ est un point r\'egulier de $X'$ et $S'$ est normalement plat le long de $X'$ en $\xi'$.
\end{propr}
Le probl\`eme \'etant local, on peut l'exprimer en termes d'anneaux. Avant cela, nous allons rappeler la notion d'\'eclatement local par rapport \`a une valuation.
\begin{defi}\label{defeclatlocal}
Soit $(R,\mathfrak{m})$ un anneau local noeth\'erien int\`egre de corps des fractions $K$. Soit $\nu$ une valuation de $K$ centr\'ee en $R$. Soient $u_1,...u_r\in R$ et $v_1,...,v_r\in R$ tels que $\nu(v_i)\leqslant\nu(u_i)$ pour tout $i\in\lbrace 1,...,r\rbrace$. Notons $R'$ l'anneau:
\[R'=R\left[\dfrac{u_1}{v_1},...,\dfrac{u_r}{v_r}\right].\]
Alors l'anneau $R^{(1)}=R'_{m_\nu\cap R'}$ est un anneau local d'id\'eal maximal $\mathfrak{m}^{(1)}=(m_\nu\cap R')R'_{m_\nu\cap R'}$.
\\Un \textbf{\'eclatement local de} $\boldsymbol R$ \textbf{par rapport \`a} $\boldsymbol{\nu}$ est un morphisme local d'anneaux locaux de la forme:
\[\pi:(R,\mathfrak m)\rightarrow (R^{(1)},\mathfrak{m}^{(1)}).\]
Soient $I$ un id\'eal de $R$ et $u_0\in I$ tel que $\nu(u_0)\leqslant\nu(f)$, pour tout $f\in I$. Compl\'etons $u_0$ en un ensemble $\lbrace u_0,u_1,...,u_s\rbrace$ de g\'en\'erateurs de $I$. Le morphisme pr\'ec\'edent est appel\'e un \textbf{\'eclatement local de} $\boldsymbol R$ \textbf{par rapport \`a} $\boldsymbol{\nu}$ \textbf{le long de} $\boldsymbol I$ si $r=s$ et $v_i=u_0$ pour tout $i\in\lbrace 1,...,s\rbrace$; conditions auxquelles on peut toujours se ramener sans perte de g\'en\'eralit\'e en posant $u_0=v_1...v_r$ et $u_i=\frac{u_i}{v_i}u_0$, $i\in\lbrace 1,...,r\rbrace$.
\end{defi}
\begin{rem}
\textup{\`A isomorphisme pr\`es, la d\'efinition pr\'ec\'edente est ind\'ependante du choix de l'ensemble de g\'en\'erateurs de $I$, c'est-\`a-dire qu'un autre choix de g\'en\'erateurs donne un anneau isomorphe.}
\end{rem}
\begin{propr}\textup{\textbf{(d'uniformisation locale des anneaux locaux).}}\label{uniflocanneaunonint}
Soient $(S,\mathfrak m)$ un anneau local noeth\'erien (non n\'ecessairement int\`egre),  $P$ un id\'eal premier minimal de $S$ et $\nu$ une valuation du corps des fractions de $S/P$ centr\'ee en $S/P$. Alors, il existe un \'eclatement local $\pi:S\rightarrow S'$ par rapport \`a $\nu$ tel que $S'_{red}$ soit r\'egulier et $Spec(S')$ soit normalement plat le long de $Spec(S'_{red})$.
\end{propr}
Nous finissons avec la notion de croisements normaux et d'uniformisation locale plong\'ee.
\begin{defi}
Soient $(R,\mathfrak m)$ un anneau local noeth\'erien et $\nu$ une valuation centr\'ee en $R$, au sens de la D\'efinition \ref{centrevaldef}, de groupe des valeurs $\Gamma$. Soit $u=\lbrace u_1,...,u_n\rbrace\subset\mathfrak m$ tel que $(u)+\sqrt{(0)}=\mathfrak m+\sqrt{(0)}$. Enfin, pour $f\in R$, on note $\overline f\in R/\sqrt{(0)}=R_{red}$ l'image de $f$ dans $R_{red}$ par le morphisme de passage au quotient.
\begin{enumerate}
\item Un mon\^ome $u^\alpha=u_1^{\alpha_1}...u_n^{\alpha_n}$ est dit \textbf{minimal par rapport \`a} $\boldsymbol \nu$ si la famille $\lbrace\nu(u_j)\:\vert\:\alpha_j\neq 0\rbrace_{1\leqslant j\leqslant n}$ est $\mathbb Z$-libre dans $\Gamma$.
\item Soit $I$ un id\'eal de $R$. On dit que le triplet $(R,I,u)$ est \`a \textbf{croisements normaux} si:
\begin{enumerate}
\item $R_{red}$ est un anneau local r\'egulier et $(\overline{u_1},...,\overline{u_n})$ est un syst\`eme r\'egulier de param\`etres de $R_{red}$;
\item $Spec(R)$ est normalement plat le long de $Spec(R_{red})$;
\item $I/\left(I+\sqrt{(0)}\right)$ est un id\'eal principal engendr\'e par un mon\^ome en $\overline{u_1},...,\overline{u_n}$ (avec la possibilit\'e que $I=(1)$ et donc $I/\left(I+\sqrt{(0)}\right)=(1)$).
\end{enumerate}
\item Soit $I$ un id\'eal de $R$, le triplet $(R,I,u)$ est \`a \textbf{croisements normaux standards} par rapport \`a $\nu$ si $(R,I,u)$ est \`a croisements normaux et $I/\left(I+\sqrt{(0)}\right)$ est engendr\'e par un mon\^ome minimal par rapport \`a $\nu$.
\item Soit $I$ un id\'eal de $R$, on dit que $(R,I)$ est \`a \textbf{croisements normaux} (resp. \`a \textbf{croisements normaux standards}) s'il existe $u$ tel que $(R,I,u)$ soit \`a croisements normaux (resp. \`a croisements normaux standards).
\item On dit que $R$ est \textbf{d\'esingularis\'e} si $(R,R)$ est \`a croisements normaux.
\end{enumerate}
\end{defi}
\begin{defi}\label{defiuniflocpourpaire}
Soient $(R,\mathfrak m)$ un anneau local noeth\'erien et $\nu$ une valuation centr\'ee en $R$ au sens de la D\'efinition \ref{centrevaldef}. Soit $I$ un id\'eal de $R$, on dit que la paire $(R,I)$ admet une \textbf{uniformisation locale plong\'ee} (resp. une \textbf{uniformisation locale plong\'ee standard}) s'il existe une suite:
\[ \xymatrix{ R \ar[r]^-{\pi_{0}} & R^{(1)} \ar[r]^-{\pi_{1}} & \ldots \ar[r]^-{\pi_{l-2}} & R^{(l-1)}  \ar[r]^-{\pi_{l-1}} & R^{(l)}} \]
o\`u, pour $1\leqslant i \leqslant l$, $\pi_i$ est un \'eclatement local par rapport \`a $\nu$ le long d'un id\'eal $J^{(i)}$ ayant les propri\'et\'es suivantes:
\begin{enumerate}
\item Pour $1\leqslant i \leqslant l$, $J^{(i)}\not\subset P_{\infty}^{(i)}$, $P_{\infty}^{(i)}$ \'etant le support de $\nu$ dans $R^{(i)}$.
\item $\left(R^{(i)},IR^{(i)}\right)$ est \`a croisements normaux (resp. \`a croisements normaux standards).
\end{enumerate}
Enfin, on dit que $R$ admet une \textbf{uniformisation locale plong\'ee} (resp. une \textbf{uniformisation locale plong\'ee standard}) si, pour tout id\'eal $I$ de $R$,  $(R,I)$ admet une uniformisation locale plong\'ee (resp. une uniformisation locale plong\'ee standard).
\end{defi}
\begin{propr}
\textup{\textbf{(d'uniformisation locale plong\'ee des sch\'emas).}}
Soit $S$ un sch\'ema noeth\'erien (non n\'ecessairement int\`egre). Soient $X$ une composante irr\'eductible de $S_{red}$ et $\nu$ une valuation de $K(X)$
centr\'ee en un point $\xi\in X$. Il existe alors un \'eclatement $\pi: S'\rightarrow S$ le long d'un sous-sch\'ema de $S$, ne contenant aucune composante irr\'eductible de $S_{red}$ et ayant la propri\'et\'e suivante : \\Soient $X'$ le transform\'e strict de $X$ par $\pi$, $\xi'$ le centre de $\nu$ sur $X'$ et $D$ le diviseur exceptionnel de $\pi$, alors $(\mathcal{O}_{X',\xi'},\mathcal{I}_{D,\xi'})$ admet une uniformisation locale plong\'ee.
\end{propr}
Dans le cas des anneaux locaux noeth\'eriens int\`egres, on peut \'enoncer la propri\'et\'e de mani\`ere un peu plus simple:
\begin{propr}\textup{\textbf{(d'uniformisation locale plong\'ee des anneaux locaux int\`egres).}}\label{uniflocplongint}
Soient $(R,\mathfrak{m})$ un anneau local noeth\'erien int\`egre et $\nu$ une valuation de $K$, le corps des fractions de $R$, centr\'ee en $R$. On dit que $\nu$ admet une \textbf{uniformisation locale plong\'ee} si, pour un nombre fini d'\'el\'ements de $R$, $f_1,...,f_q\in R$ tels que $\nu(f_1)\leqslant ...\leqslant\nu(f_q)$, il existe une suite d'\'eclatements locaux par rapport \`a $\nu$:
\[ \xymatrix{ R \ar[r]^-{\pi_{0}} & R^{(1)} \ar[r]^-{\pi_{1}} & \ldots \ar[r]^-{\pi_{l-2}} & R^{(l-1)}  \ar[r]^-{\pi_{l-1}} & R^{(l)}} \]
telle que $R^{(l)}$ soit r\'egulier et telle qu'il existe un syst\`eme r\'egulier de param\`etres $u^{(l)}=\left(u_1^{(l)},...,u_d^{(l)}\right)$ de $R^{(l)}$ tel que les $f_i$, $1\leqslant i\leqslant q$, soient des mon\^omes en $u^{(l)}$ multipli\'es par une unit\'e de $R^{(l)}$ et $f_1/.../f_q$ dans $R^{(l)}$.
\end{propr}
\section{Polyn\^omes-cl\'es en caract\'eristique nulle}
Consid\'erons $K\hookrightarrow K(x)$ une extension de corps simple et transcendante. Soit $\mu'$ une valuation de $K(x)$, notons $\mu:=\mu'_{\vert\:K}$. On note $G$ le groupe des valeurs de $\mu'$ et $G_{1}$ celui de $\mu$. On suppose de plus que $\mu$ est de rang $1$, $\mu'(x)>0$ et $car(k_\mu)=0$. Enfin, pour $\beta\in G$, on pose:
\[P'_{\beta}=\lbrace f\in K(x)\:\vert\:\mu'(f)\geqslant\beta\rbrace\cup\lbrace 0\rbrace;\]
\[P'_{\beta,+}=\lbrace f\in K(x)\:\vert\:\mu'(f)>\beta\rbrace\cup\lbrace 0\rbrace;\]
\[G_{\mu'}=\bigoplus\limits_{\beta\in G}P'_{\beta}/P'_{\beta,+};\]
et $in_{\mu'}(f)$ l'image de $f\in K(x)$ dans $G_{\mu'}$.
\begin{defi}
Un \textbf{ensemble complet de polyn\^omes-cl\'es} pour $\mu'$ est une collection bien ordonn\'ee:
\[\textbf{Q}=\lbrace Q_{i}\rbrace_{i\in\Lambda}\subset K[x]\]
telle que, pour tout $\beta\in G$, le groupe additif $P'_{\beta}\cap K[x]$ soit engendr\'e par des produits de la forme $a\prod\limits_{j=1}^{s}Q_{i_{j}}^{\gamma_{j}}$, $a\in K$, tels que $\sum\limits_{j=1}^{s}\gamma_{j}\mu'\left(Q_{i_{j}}\right)+\mu(a)\geqslant\beta$.
\\L'ensemble est dit \textbf{1-complet} si la condition a lieu pour tout $\beta\in G_1$.
\end{defi}
\begin{thm}(\cite{spivaherrera}, Th\'eor\`eme 62)\label{existpolycle}
Il existe une collection $\textbf{Q}=\lbrace Q_{i}\rbrace_{i\in\Lambda}$ qui soit un ensemble $1$-complet de polyn\^omes-cl\'es.
\end{thm}
Par le Th\'eor\`eme \ref{existpolycle}, on sait qu'il existe un ensemble $1$-complet de polyn\^omes-cl\'es $\textbf{Q}=\lbrace Q_{i}\rbrace_{i\in\Lambda}$ et que le type d'ordre de $\Lambda$ est au plus $\omega\times\omega$. Si $K$ est sans d\'efaut, on va voir que le type d'ordre de $\Lambda$ est au plus $\omega$ et que, par cons\'equent, il n'y a pas de polyn\^omes-cl\'es limites, ce sera en particulier le cas si $car(k_\mu)=0$. Pour tout $i\in\Lambda$, notons $\beta_i=\mu'(Q_i)$.

Soit $l\in\Lambda$, on note:
\[\alpha_{i}=d_{Q_{i-1}}^{\:\circ}(Q_{i}),\:\forall\:i\leqslant l;\]
\[\boldsymbol{\alpha_{l+1}}=\lbrace \alpha_{i}\rbrace_{i\leqslant \:l};\]
\[\textbf{Q}_{l+1}=\lbrace Q_{i}\rbrace_{i\leqslant \:l}.\]
On utilise \'egalement la notation $\overline{\gamma}_{l+1}=\lbrace \gamma_{i}\rbrace_{i\leqslant\:l}$ o\`u les $\gamma_{i}$ sont tous nuls sauf pour un nombre fini d'entres eux, $\textbf{Q}_{l+1}^{\overline{\gamma}_{l+1}}=\prod\limits_{i\leqslant\: l}Q_{i}^{\gamma_{i}}$.

\begin{defi}
Un multi-indice $\overline{\gamma}_{l+1}$ est dit \textbf{standard par rapport \`a} $\boldsymbol{\alpha_{l+1}}$ si $0\leqslant \gamma_{i}<\alpha_{i+1}$, pour $i\leqslant l$.
\\ Un \textbf{mon\^ome l-standard en} $\boldsymbol{Q_{l+1}}$ est un produit de la forme $c_{\overline{\gamma}_{l+1}}\textbf{Q}_{l+1}^{\overline{\gamma}_{l+1}}$, o\`u $c_{\overline{\gamma}_{l+1}}\in K$ et $\overline{\gamma}_{l+1}$ est standard par rapport \`a $\boldsymbol{\alpha_{l+1}}$.
\\ Un \textbf{d\'eveloppement l-standard n'impliquant pas} $\boldsymbol{Q_{l}}$ est une somme finie $\sum\limits_{\beta}S_{\beta}$ de mon\^omes $l$-standards n'impliquant pas $Q_{l}$, o\`u $\beta$ appartient \`a un sous-ensemble fini de $G_{+}$ et $S_{\beta}=\sum\limits_{j} d_{\beta,j}$ est une somme de mon\^omes standards de valuation $\beta$ v\'erifiant $\sum\limits_{j} in_{\mu'}(d_{\beta,j})\neq 0$.
\end{defi}
\begin{defi}
Soient $f\in K[x]$ et $i\leqslant l$, un \textbf{d\'eveloppement i-standard de f} est une expression de la forme:
\[f=\sum\limits_{j=0}^{s_{i}}c_{j,i}Q_{i}^{j},\]
o\`u $c_{j,i}$ est un d\'eveloppement $i$-standard n'impliquant pas $Q_{i}$.
\end{defi}
\begin{rem}\textup{
Un tel d\'eveloppement existe, par division Euclidienne et est unique dans le sens o\`u les $c_{j,i}\in K[x]$ sont uniques. Plus pr\'ecis\'ement, si $i\in\mathbb N$, on montre par r\'ecurrence que le d\'eveloppement $i$-standard est unique.}
\end{rem}
\begin{defi}
Soient $f\in K[x]$, $i\leqslant l$ et $f=\sum\limits_{j=0}^{s_{i}}c_{j,i}Q_{i}^{j}$ un d\'eveloppement $i$-standard de $f$. On d\'efinit la \textbf{i-troncature de} $\boldsymbol{\mu'}$, not\'ee $\mu_{i}'$, comme \'etant la pseudo-valuation:
\[\mu_{i}'(f)=\min_{0\leqslant j\leqslant s_{i}}\lbrace j\mu'(Q_{i})+\mu'(c_{j,i})\rbrace.\]
\end{defi}
\begin{rem}\label{inegalitetronque}\textup{
On peut montrer que c'est en fait une valuation. On a de plus:
\[\forall \:f\in K[x],\: i\in\Lambda, \:\mu_{i}'(f)\leqslant \mu'(f).\]}
\end{rem}

\indent La construction des polyn\^omes-cl\'es se fait par r\'ecurrence (voir \cite{mahboub}, \cite{mahboubthese}, \cite{spiva} \textsection 9 et \cite{spivaherrera}). Pour $l\in\mathbb N^*$, on construit un ensemble de polyn\^omes-cl\'es $\textbf{Q}_{l+1}=\lbrace Q_{i}\rbrace_{1\leqslant i\leqslant \:l}$; deux cas se pr\'esentent:
\begin{enumerate}
\item[(1)] $\exists\:l_0\in\mathbb N,\:\beta_{l_0}\notin G_1$;
\item[(2)] $\forall\:l\in\mathbb N,\:\beta_l\in G_1$.
\end{enumerate}
Dans le cas (1), on stoppe la construction; l'ensemble $\textbf{Q}_{l_0}=\lbrace Q_{i}\rbrace_{1\leqslant i\leqslant \:l_0-1}$ est par d\'efinition un ensemble $1$-complet de polyn\^omes-cl\'es et $\Lambda=\lbrace 1,...,l_0-1\rbrace$. Remarquons de plus que l'ensemble $\textbf{Q}_{l_0+1}$ est quant \`a lui un ensemble complet de polyn\^omes-cl\'es.
\\Dans le cas (2), l'ensemble $\textbf{Q}_{\omega}=\lbrace Q_{i}\rbrace_{i\geqslant 1}$ est infini et $\Lambda=\mathbb N^*$. Les propositions qui suivent nous assurent que dans ce cas, l'ensemble des polyn\^omes-cl\'es obtenu est \'egalement $1$-complet.
\\\indent Le lemme qui suit, tr\`es utile dans la suite, nous permet de remarquer qu'il n'existe pas de suite croissante born\'ee dans le cadre des valuations de rang $1$.
\begin{lem}\label{pasdesuitecroissantebornee}
Soit $\nu$ une valuation de rang $1$ centr\'ee en un anneau local noeth\'erien $R$. Notons $P_\infty$ le support de $\nu$. Alors, $\nu\left( R\setminus P_\infty\right)$ ne contient aucune suite infinie croissante et born\'ee.
\end{lem}
\begin{rem}
\textup{La donn\'ee d'une valuation $\nu$ centr\'ee en un anneau local $(R,\mathfrak{m})$ est la donn\'ee d'un id\'eal premier minimal $P_\infty$ de $R$ (le support de la valuation) et d'une valuation $\nu'$ du corps des fractions de $R/P_\infty$ telle que $R/P_\infty\subset R_{\nu'}$ et $\mathfrak m/P_\infty=(R/P_\infty)\cap m_{\nu'}$.}
\end{rem}
\noindent\textit{Preuve}: Soit $\left(\beta_i\right)_{i\geqslant 1}$ une suite croissante infinie de $\nu\left( R\setminus P_\infty\right)$ born\'ee par $\beta$. Cette suite correspond \`a une suite infinie d\'ecroissante d'id\'eaux de $R/P_\beta$. Il nous suffit donc de montrer que $R/P_\beta$ est de longueur finie. Notons $\mathfrak m$ l'id\'eal maximal de $R$, $\nu(\mathfrak m)=\min\left\lbrace\nu\left( R\setminus P_\infty\right)\setminus\lbrace 0\rbrace\right\rbrace$ et $\Gamma$ le groupe des valeurs de $\nu$. Remarquons que le groupe $\nu\left( R\setminus P_\infty\right)$ est archim\'edien. En effet, par l'absurde, si $\nu\left( R\setminus P_\infty\right)$ n'est pas archim\'edien, il existe $\alpha,\:\beta\in\nu\left( R\setminus P_\infty\right)$, $\beta\neq 0$ tels que, pour tout $n\geqslant 1$, $n\beta\leqslant\alpha$. En particulier, l'ensemble:
\[\lbrace \gamma\in\Gamma\:\vert\:\exists\:n\in\mathbb{N}\setminus\lbrace 0\rbrace,\:-n\beta<\gamma<n\beta\rbrace\]
est un sous-groupe isol\'e non trivial de $\Gamma$.
\\On en d\'eduit qu'il existe $n\in\mathbb N$ tel que:
\[\beta\leqslant n\nu(\mathfrak m).\]
Ainsi, $\mathfrak m^n\subset P_\beta$ et donc , il existe une application surjective:
\[R/\mathfrak m^n\twoheadrightarrow R/P_\beta.\]
On en d\'eduit que $R/P_\beta$ est de longueur finie, ce qui est absurde. On en conclut que $\nu\left( R\setminus P_\infty\right)$ ne contient aucune suite infinie croissante born\'ee.\\\qed

\begin{prop}\label{sibornealorscomplet}(\cite{spiva}, Proposition 9.30) Supposons que l'on ait construit un ensemble infini de polyn\^omes-cl\'es $\textbf{Q}_{\omega}=\lbrace Q_{i}\rbrace_{i\geqslant 1}$ tel que, pour tout $i\in\mathbb N^*$, $\beta_i\in G_1$. Supposons de plus que la suite $\lbrace\beta_i\rbrace_{i\geqslant 1}$ n'est pas born\'ee dans $G_1$. Alors, l'ensemble de polyn\^omes-cl\'es $\textbf{Q}_{\omega}$ est $1$-complet.
\end{prop}
\noindent\textit{Preuve}: Il suffit de montrer que, pour tout $\beta\in G_1$ et pour tout $h\in K\left[x\right]$ tels que $\mu'(h)=\beta$, $h$ est dans le $R_\mu$-sous-module de $K\left[x\right]$ engendr\'e par tous les mon\^omes de la forme $a\prod\limits_{j=1}^{s}Q_{i_{j}}^{\gamma_{j}}$, $a\in K$, tels que $\mu'\left(a\prod\limits_{j=1}^{s}Q_{i_{j}}^{\gamma_{j}}\right)\geqslant\beta$.
\\Consid\'erons donc $h\in K\left[x\right]$ tel que $\mu'(h)\in G_1$. En notant $h=\sum\limits_{j=0}^dh_jx^j$, on peut supposer, sans perte de g\'en\'eralit\'e, que:
\[\forall\:j\in\lbrace 0,...,d\rbrace,\:\mu(h_j)\geqslant 0.\]
En effet, dans le cas contraire, il suffit de multiplier $h$ par un \'el\'ement de $K$ choisi convenablement.
\\Comme la suite $\lbrace\beta_i\rbrace_{i\geqslant 1}$ n'est pas born\'ee dans $G_1$, il existe $i_0\in\mathbb N^*$ tel que:
\[\mu'(h)<\beta_{i_0}.\]
Notons alors $h=\sum\limits_{j=0}^{s_{i_0}}c_{j,i_0}Q_{i_0}^j$, le d\'eveloppement $i_0$-standard de $h$. Or, comme ce d\'eveloppement est obtenu par division euclidienne, vu le choix fait sur les coefficients de $h$ et, comme la suite $\left\lbrace\frac{\beta_i}{d^{\:\circ}\left(Q_i\right)}\right\rbrace_{i\geqslant 1}$ est strictement croissante (il suffit de regarder le d\'eveloppement $(i-1)$-standard de $Q_i$), on montre facilement que:
\[\forall\:j\in\lbrace 0,...,s_{i_0}\rbrace,\:\mu\left(c_{j,i_0}\right)\geqslant 0.\]
Rappelons que, par construction des polyn\^omes-cl\'es, pour $j\in\lbrace 0,...,s_{i_0}\rbrace$, $\mu'_{i_0}\left(c_{j,i_0}\right)=\mu'\left(c_{j,i_0}\right)$. On en d\'eduit alors que:
\[\forall\:j\in\lbrace 1,...,s_{i_0}\rbrace,\:\mu'\left(c_{j,i_0}Q_{i_0}^j\right)=\mu'_{i_0}\left(c_{j,i_0}Q_{i_0}^j\right)>\mu'(h).\]
Ainsi, $\mu'(h)=\mu'\left(c_{0,i_0}\right)$ et donc, $h$ est une somme de mon\^omes en $\textbf{Q}_{i_0+1}$ de valuation au moins $\mu'(h)$ (et en particulier, $\mu'_{i_0}(h)=\mu'(h)$).\\\qed
\\ \\\indent 
On consid\`ere alors deux cas:
\begin{enumerate}
\item[(1)] $\sharp\lbrace i\geqslant 1\:\vert\:\alpha_i>1\rbrace=+\infty$;
\item[(2)] $\sharp\lbrace i\geqslant 1\:\vert\:\alpha_i>1\rbrace<+\infty$.
\end{enumerate}
Dans le cas (1), \`a l'aide de la Proposition \ref{sialphainfini}, on montre que l'ensemble infini de polyn\^omes-cl\'es est toujours $1$-complet, ind\'ependamment de la caract\'eristique de $k_\mu$. Dans le cas (2), si la caract\'eristique de $k_\mu$ est nulle et si l'ensemble de polyn\^omes-cl\'es $\textbf{Q}_{\omega}=\lbrace Q_{i}\rbrace_{i\geqslant 1}$ n'est pas complet, on montre dans la Proposition \ref{sialphafini} que la suite $\lbrace\beta_i\rbrace_{i\geqslant 1}$ n'est jamais born\'ee. Dans ce cas-l\`a, gr\^ace \`a la Proposition \ref{sibornealorscomplet}, on en d\'eduit que l'ensemble de polyn\^omes-cl\'es $\textbf{Q}_{\omega}=\lbrace Q_{i}\rbrace_{i\geqslant 1}$ est \'egalement $1$-complet.
\begin{prop}\label{sialphainfini}(\cite{spiva}, Corollaire 11.8)
Supposons que l'on ait construit un ensemble infini de polyn\^omes-cl\'es $\textbf{Q}_{\omega}=\lbrace Q_{i}\rbrace_{i\geqslant 1}$ tel que, pour tout $i\in\mathbb N^*$, $\beta_i\in G_1$. Supposons de plus que l'ensemble $\lbrace i\geqslant 1\:\vert\:\alpha_i>1\rbrace$ est infini. Alors, $\textbf{Q}_{\omega}$ est un ensemble $1$-complet de polyn\^omes-cl\'es.
\end{prop}
\noindent\textit{Preuve}: Soit $h\in K[x]$, comme dans la preuve de la Proposition \ref{sibornealorscomplet}, il suffit de montrer que $\mu'_{i}(h)=\mu'(h)$ pour un certain $i\geqslant 1$. Or, si on note:
\[\delta_i(h)=\max S_i(h,\beta_i),\]
o\`u:
\[S_i(h,\beta_i)=\lbrace j\in\lbrace 0,...,s_{i}\rbrace\:\vert\:j\beta_i+\mu'\left(c_{j,i}\right)=\mu'_i(h)\rbrace,\]
\[h=\sum\limits_{j=0}^{s_{i}}c_{j,i}Q_{i}^j,\]
par le (1) de la Proposition 37 de \cite{spivaherrera} (Proposition 11.2 de \cite{spiva}), on a:
\[\alpha_{i+1}\delta_{i+1}(h)\leqslant\delta_i(h),\:\forall\:i\geqslant 1.\]
On en d\'eduit qu'\`a chaque fois que $\delta_i(h)>0$ et $\alpha_{i+1}>1$:
\[\delta_{i+1}(h)<\delta_i(h),\forall\:i\geqslant 1.\]
Comme l'ensemble $\lbrace i\geqslant 1\:\vert\:\alpha_i>1\rbrace$ est infini et que l'in\'egalit\'e pr\'ec\'edente ne peut se produire une infinit\'e de fois, on en conclut qu'il existe $i_0\geqslant 1$ tel que $\delta_{i_0}(h)=0$ et donc que $\mu'_{i_0}(h)=\mu'(h)$.\\\qed
\\ \\\indent \`A partir de maintenant, on suppose que l'on a construit un ensemble infini de polyn\^omes-cl\'es $\textbf{Q}_{\omega}=\lbrace Q_{i}\rbrace_{i\geqslant 1}$ tel que $\alpha_i=1$, pour tout $i$ suffisamment grand. Ainsi pour ces $i$, on a:
\[Q_{i+1}=Q_i+z_i,\]
o\`u $z_i$ est un d\'eveloppement $i$-standard homog\`ene, de valuation $\beta_i$, n'impliquant pas $Q_i$. 
\begin{prop}\label{sialphafini}(\cite{spiva}, Proposition 12.8)
Supposons que $car(k_\mu)=0$ et que l'on ait construit un ensemble infini de polyn\^omes-cl\'es $\textbf{Q}_{\omega}=\lbrace Q_{i}\rbrace_{i\geqslant 1}$ tel que, pour tout $i\in\mathbb N^*$, $\beta_i\in G_1$. Supposons de plus qu'il existe $h\in K[x]$ tel que, pour tout $i\geqslant 1$: $$\mu'_i(h)<\mu'(h).$$ 
Alors, la suite $\lbrace\beta_i\rbrace_{i\geqslant 1}$ n'est pas born\'ee dans $G_1$.
\end{prop}
\noindent\textit{Preuve}: Par la Proposition 37 de \cite{spivaherrera} (Proposition 11.2 de \cite{spiva}), la suite $\lbrace\delta_i(h)\rbrace_{i\geqslant 1}$ est d\'ecroissante, il existe donc $i_0\geqslant 1$ tel que $\delta_{i_0+t}(h)=\delta_{i_0}(h)$, pour tout $t\in\mathbb N$. Notons $\delta$ cette valeur commune. Si on note $h=\sum\limits_{j=0}^{s_{i}}c_{j,i}Q_{i}^j$ le d\'eveloppement $i$-standard de $h$ pour $i\geqslant i_0$, alors, par la Proposition 37 de \cite{spivaherrera} (Proposition 11.2 de \cite{spiva}), $\mu'_i(h)=\delta\beta_i+\mu'\left(c_{\delta,i}\right)$ et $\mu'\left(c_{\delta,i}\right)$ sont ind\'ependants de $i$. Il suffit donc de montrer que la suite $\lbrace\mu'_i(h)\rbrace_{i\geqslant 1}$ n'est pas born\'ee.
\\Notons:
\[\mu_i^+(h)=\min\left\lbrace\mu'\left.\left(c_{j,i} Q_i^j\right)\:\right\vert\:\delta<j\leqslant s_i\right\rbrace,\]
\[\varepsilon_i(h)=\min\left\lbrace j\in\lbrace\delta+1,...,s_i\rbrace\:\left\vert\:\mu'\left(c_{j,i} Q_i^j\right)=\mu_i^+(h)\right.\right\rbrace.\]
Toujours par la Proposition 37 de \cite{spivaherrera} (Proposition 11.2 de \cite{spiva}), la suite $\lbrace\varepsilon_i(h)\rbrace_{i\geqslant i_0}$ est d\'ecroissante, il existe donc $i_1\geqslant i_0$ tel que cette suite soit constante \`a partir de $i_1$. Notons alors $c_{\delta,i_1}^*\in K[x]$ l'unique polyn\^ome de degr\'e strictement inf\'erieur \`a $d^{\:\circ}\left(Q_{i_0}\right)=d^{\:\circ}\left(Q_{i_1}\right)$ tel que $c_{\delta,i_1}^*c_{\delta,i_1}-1$ soit divisible par $Q_{i_1}$ dans $K[x]$. On peut montrer que $\mu'_i\left(c_{\delta,i_1}^*\right)=\mu'\left(c_{\delta,i_1}^*\right)$, pour tout $i\geqslant i_1$. Multiplier $h$ par $c_{\delta,i_1}^*$ n'affecte pas $\delta$, donc multiplier $h$ par $c_{\delta,i_1}^*$ ne change rien au probl\`eme. On peut donc supposer que $in_{\mu'}\left(c_{\delta,i}\right)=in_{\mu'_i}\left(c_{\delta,i}\right)=1$ pour tout $i\geqslant i_1$.
\\Supposons que $h=Q_{i_1+1}$, rappelons que nous sommes dans la situation o\`u $Q_{i+1}=Q_i+z_i$, pour $i\geqslant i_1\geqslant i_0$. Les $z_i$ n'\'etant pas uniques, un choix possible de $z_i$, pour $i=i_1$, est:
\[z_{i_1}=\dfrac{c_{\delta-1,i_1}}{\delta}.\]
Par d\'efinition de $z_{i_1}$, $\mu'\left(z_{i_1}\right)=\beta_{i_1}$ et $\beta_{i_1}<\beta_{i_1+1}$. Par r\'ecurrence sur $t\in\mathbb N$, on construit $Q_{i_1+t}$.
\\Il faut tout de m\^eme montrer que la propri\'et\'e \og$\lbrace\mu'(Q_i+z_i+...+z_l)\rbrace_l$ n'est pas born\'ee\fg\: ne d\'epend pas du choix des $z_i,...,z_l$, $i\leqslant l$. En effet, supposons que l'on ait construit une autre suite de la forme $\lbrace\mu'(Q_i+z'_i+...+z'_{l'})\rbrace_{l'}$. Si pour tout $l$, il existe $l'$ tel que $\mu'(Q_i+z_i+...+z_{l})<\mu'(Q_i+z'_i+...+z'_{l'})$ alors la suite $\lbrace\mu'(Q_i+z'_i+...+z'_{l'})\rbrace_{l'}$ ne peut pas \^etre born\'ee car sinon la suite $\lbrace\mu'(Q_i+z_i+...+z_l)\rbrace_l$ le serait ce qui contredit l'hypoth\`ese de d\'epart. Supposons donc qu'il existe $l$ tel que, pour tout $l'$, $\mu'(Q_i+z'_i+...+z'_{l'})<\mu'(Q_i+z_i+...+z_l)$. Par la Proposition 9.29 de \cite{spiva}, il existe un d\'eveloppement $Q_i+z'_i+...+z'_{l'}+z''_{l'+1}+...+z''_{l''}$ tel que $Q_i+z''_i+...+z''_{l''}=Q_i+z_i+...+z_l$. Ainsi, on peut construire une troisi\`eme suite qui n'est pas born\'ee. 
\\Comme $car(k_\mu)=0$, le sous-corps premier de $K$ est $\mathbb Q$, consid\'erons alors $A$ la $\mathbb Q$-sous-alg\`ebre de $K$ engendr\'ee par tous les coefficients de $Q_{l_1}$, on a donc que, pour tout $t\in\mathbb N$, $Q_{i_1+t}\in A\left[x\right]$. L'anneau $A$ \'etant noeth\'erien, l'anneau $A\left[x\right]$ l'est aussi. La valuation $\mu'_{\vert\:A[x]}$ est alors centr\'ee en $A[x]$ et $\left\lbrace \mu'_{\vert\:A[x]}\left(Q_{i_1+t}\right)\right\rbrace_{t\in\mathbb N}\subset G_1$, $G_1$ \'etant de rang $1$. En appliquant le Lemme \ref{pasdesuitecroissantebornee}, on en d\'eduit que la suite $\lbrace\beta_i\rbrace_{i\geqslant 1}$ ne peut \^etre born\'ee dans $G_1$.\\\qed
\begin{coro}\label{polycleencar0}
Si $car(k_\mu)=0$, il existe un ensemble $1$-complet de polyn\^omes-cl\'es $\lbrace Q_i\rbrace_{i\in\Lambda}$ tel que $\Lambda$ est, soit un ensemble fini, soit $\mathbb N^*$. En particulier, il n'y a pas de polyn\^omes-cl\'es limites  pour des valuations de rang $1$ dont le corps r\'esiduel est de caract\'eristique nulle.
\end{coro}
\noindent\textit{Preuve}: On applique le processus de construction de \cite{spiva}, \textsection 9 et \cite{spivaherrera}. S'il existe $i_0\in\mathbb N$, tel que $\beta_{i_0}\notin G_1$, on pose $\Lambda=\lbrace 1,...,i_0-1\rbrace$ et, par d\'efinition, $\lbrace Q_i\rbrace_{i\in\Lambda}$ est $1$-complet. Sinon, pour tout $i\in\mathbb N$, $\beta_i\in G_1$ et on pose $\Lambda=\mathbb N^*$. Si $\sharp\lbrace i\geqslant 1\:\vert\:\alpha_i>1\rbrace=+\infty$, par la Proposition \ref{sialphainfini}, l'ensemble $\lbrace Q_i\rbrace_{i\in\Lambda}$ est $1$-complet. Si $\sharp\lbrace i\geqslant 1\:\vert\:\alpha_i>1\rbrace<+\infty$, par la Proposition \ref{sialphafini}, la suite $\lbrace\beta_i\rbrace_{i\geqslant 1}$ n'est pas born\'ee dans $G_1$ et donc, par la Proposition \ref{sibornealorscomplet}, l'ensemble $\lbrace Q_i\rbrace_{i\in\Lambda}$ est un ensemble $1$-complet de polyn\^omes-cl\'es.\\\qed

\section{Pr\'eliminaires}
\indent Les \'eclatements locaux sont un outil essentiel pour obtenir un r\'esultat d'uniformisation locale. Ces \'eclatements sont d\'ependants du choix des diff\'erents param\`etres r\'eguliers possibles pour l'anneau d'arriv\'ee. Les \'eclatements locaux encadr\'es vont imposer un syst\`eme de g\'en\'erateurs de l'id\'eal maximal d'arriv\'ee pour permettre de faire d\'ecro\^itre des invariants (dimension de plongement et rang rationnel d'un sous-groupe de l'enveloppe divisible du groupe des valeurs de la valuation).
\\Pour les preuves de tous les r\'esultats, on pourra consulter \cite{spiva}: \textsection 5, \textsection 6, \textsection 7, \textsection 8 ainsi que le Chapitre I de \cite{jcthese}.
\subsection{Suites locales encadr\'ees}\label{sectioneclatencad}
~\smallskip ~\\ \indent Soit $(R,\mathfrak{m},k)$ un anneau local noeth\'erien. Notons:
\[u=(u_1,...,u_n)\]
un ensemble de g\'en\'erateurs de $\mathfrak m$. Pour un sous-ensemble $I\subset\lbrace 1,...,n\rbrace$, notons:
\[u_I=\lbrace u_i\:\vert\:i\in I\rbrace.\]
Fixons un sous-ensemble $J\subset\lbrace 1,...,n\rbrace$ et un \'el\'ement $j\in J$. Notons:
\[J^c=\lbrace 1,...,n\rbrace\setminus J.\]
Pour tout $i\in \lbrace 1,...,n\rbrace$, consid\'erons les changements de variables suivants:
\[ u'_i=\left \{ \begin{array}{ccl}  u_i & \textup{si} & i\in J^c\cup\lbrace j\rbrace \\  \frac{u_i}{u_j} & \textup{si} & i\in J\setminus\lbrace j\rbrace \end{array} \right.\]
On note alors:
\[u'=(u'_1,...,u'_n).\]
Rappelons que pour $f\in R$, l'\textbf{annulateur de} $\boldsymbol f$, not\'e $Ann_R(f)$, est l'id\'eal de $R$ d\'efini par:
\[Ann_R(f)=\lbrace g\in R\:\vert\:gf=0\rbrace.\]
Pour tout $i\in \lbrace 1,...,n\rbrace$, notons:
\[Ann_R\left(u_i^\infty\right)=\bigcup\limits_{l\geqslant 1}Ann_R\left(u_i^l\right),\]
\[R_i=R/Ann_R\left(u_i^\infty\right)\textup{ et }R'=R_j\left[u'_{J\setminus\lbrace j\rbrace}\right].\]
Notons $(R^{(1)},\mathfrak{m}^{(1)},k^{(1)})$ le localis\'e de l'anneau $R'$ en un id\'eal premier de $R'$. 
\begin{rem}
\textup{Le sch\'ema $Spec\left(R'\right)$ est un sous-sch\'ema affine de l'\'eclat\'e de $Spec(R)$ le long de l'id\'eal $\left(u_J\right)$.}
\end{rem}
\indent Enfin, nous r\'ealisons une partition de $\lbrace 1,...,n\rbrace$ comme suit:
\[J^\times=\lbrace i\in J\setminus\lbrace j\rbrace\:\vert\:u'_i\in R^{(1)\times}\rbrace,\]
\[J^{\times c}=\lbrace i\in J\setminus\lbrace j\rbrace\:\vert\:u'_i\not\in R^{(1)\times}\rbrace.\]
On a donc:
\[\lbrace 1,...,n\rbrace=J^c\amalg J^\times \amalg J^{\times c}\amalg\lbrace j\rbrace,\]
\[u'=u'_{J^c}\cup u'_{J^\times}\cup u'_{J^{\times c}}\cup \lbrace u'_j\rbrace,\]
o\`u les r\'eunions sont disjointes dans la derni\`ere \'egalit\'e si $R$ est un anneau r\'egulier avec $u$ pour syst\`eme r\'egulier de param\`etres.
\\Notons $u^{(1)}=\left(u^{(1)}_1,...,u_{n_1}^{(1)}\right)$ un syst\`eme de g\'en\'erateurs de $\mathfrak m^{(1)}$ et \[\pi:(R,u)\rightarrow \left(R^{(1)},u ^{(1)}\right)\] le morphisme naturel entre ces deux anneaux locaux.
\begin{defi}\label{defeclatloc}
On dit que $\pi:(R,u)\rightarrow \left(R^{(1)},u ^{(1)}\right)$ est un \textbf{\'eclatement encadr\'e} de $(R,u)$ si $n_1\leqslant n$ et s'il existe un sous-ensemble $D_1\subset\lbrace 1,...,n_1\rbrace$ tel que:
\[u'_{\lbrace 1,...,n\rbrace\setminus J^\times}=u'_{J^c\cup J^{\times c}\cup\lbrace j\rbrace}=u_{D_1}^{(1)}.\]
Si de plus, $R$ est r\'egulier, $u$ est un syst\`eme r\'egulier de param\`etres de $R$ et $J^\times=\emptyset$ (c'est-\`a-dire si $n=n_1$ et $u'=u_{D_1}^{(1)}$), on dit que $\pi$ est un \textbf{\'eclatement monomial}.
\\Enfin, une \textbf{suite locale encadr\'ee} est une suite de la forme:
\[ \xymatrix{\left( R,u\right)=\left( R^{(0)},u^{(0)}\right) \ar[r]^-{\pi_{0}} & \left( R^{(1)},u^{(1)}\right) \ar[r]^-{\pi_{1}} & \ldots   \ar[r]^-{\pi_{l-1}} & \left( R^{(l)},u^{(l)}\right)}, \]
o\`u chaque $\pi_i:\left( R^{(i)},u^{(i)}\right)\rightarrow\left( R^{(i+1)},u^{(i+1)}\right)$, $0\leqslant i\leqslant l-1$, est un \'eclatement encadr\'e. Si de plus, pour tout $i$, les $\pi_i$ sont des \'eclatements monomiaux, on dit que la suite est \textbf{monomiale}.
\end{defi}
\begin{defi}
Soient $\pi:(R,u)\rightarrow \left(R^{(1)},u ^{(1)}\right)$ un \'eclatement encadr\'e et $T\subset\lbrace 1,...,n\rbrace$. Supposons que $R$ est r\'egulier et que $u$ est un syst\`eme r\'egulier de param\`etres de $R$.
\\On dit que $\pi$ est \textbf{ind\'ependant de} $\boldsymbol{u_T}$ si $T\cap J=\emptyset$ (c'est-\`a-dire, $T\subset J^c$).
\end{defi}
\begin{rem}
\textup{Si un \'eclatement encadr\'e est ind\'ependant de $u_T$, alors: \[u_T\subset\left\lbrace u^{(1)}_1,...,u_{n_1}^{(1)}\right\rbrace.\] }
\end{rem}
On d\'efinit par r\'ecurrence l'ind\'ependance pour une suite locale encadr\'ee en supposant qu'elle est d\'ej\`a d\'efinie pour des suites de longueur $l-1$.
\begin{defi}
Une suite locale encadr\'ee de la forme:
\[ \xymatrix{\left( R,u\right)=\left( R^{(0)},u^{(0)}\right) \ar[r]^-{\pi_{0}} & \left( R^{(1)},u^{(1)}\right) \ar[r]^-{\pi_{1}} & \ldots   \ar[r]^-{\pi_{l-1}} & \left( R^{(l)},u^{(l)}\right)} \]
est dite \textbf{ind\'ependante de} $\boldsymbol{u_T}$ si elle v\'erifie les deux conditions suivantes:
\begin{enumerate}
\item la suite $\pi_{l-2}\circ...\circ\pi_{0}$ est ind\'ependante de $u_T$;
\item si $u_T\subset\left\lbrace u^{(i)}_1,...,u_{n_i}^{(i)}\right\rbrace$, $0\leqslant i\leqslant l-1$, alors $\pi_{l-1}$ est ind\'ependante de $u_T$.
\end{enumerate}
\end{defi}
\begin{rem}\label{matriceeclat}
\textup{Soit $q\in\lbrace 1,...,n\rbrace$, on peut \'ecrire $u'_q$ sous la forme:
\[u'_q=u_1^{m_{1,q}}...u_{n}^{m_{n,q}},\]
o\`u $m_{p,q}\in\mathbb Z$, $p\in\lbrace 1,...,n\rbrace$. Le changement de variables $u\rightarrow u'$ est alors donn\'e par la matrice $M=(m_{p,q})_{p,q}\in SL_n(\mathbb Z)$ avec, par d\'efinition:
\[ m_{p,q}=\left \{ \begin{array}{rcl}  1 & \textup{si} & p=q \\  -1 & \textup{si} & p=j\textup{ et }q\in J
\\ 0 & \textup{si} & p\neq q \textup{ et, ou bien }q\neq j\textup{, ou bien }q\not\in J \end{array} \right.\]
En particulier, si $q\in J^c$, alors $u'_{q}\in u_{J^c}$ et si $q\in J$ alors $u'_q$ est un mon\^ome en $u_J$.
\\De m\^eme, le changement de variables $u'\rightarrow u$ est donn\'e par la matrice $N=(n_{p,q})_{p,q}=M^{-1}\in SL_n(\mathbb Z)$ avec:
\[ n_{p,q}=\left \{ \begin{array}{clc}  1 & \textup{si} & \textup{ou bien }p=q \textup{, ou bien }p=j\textup{ et }q\in J\\  0 & \textup{sinon} &  \end{array} \right.\]
En particulier, si $q\in J^c$, alors $u_q\in u'_{J^c}$ et si $q\in J$, alors $u_q$ est un mon\^ome en $u'_J$.
\\Si on note $e=\#(J^c\cup J^{\times c}\cup\lbrace j\rbrace)$, on en d\'eduit qu'il existe $\beta_q\in\mathbb N^q$ et $z_q\in R'^\times$ tels que:
\[u_q=\left(u_{J^c\cup J^{\times c}\cup\lbrace j\rbrace}\right)^{\beta_q}z_q.\]
De plus, si $q\in J$, alors $\left(u_{J^c\cup J^{\times c}\cup\lbrace j\rbrace}\right)^{\beta_q}$ est un mon\^ome en $u'_{J^{\times c}\cup\lbrace j\rbrace}$ uniquement. On a \'egalement:
\[\mathfrak mR'=\left(u_{J^c\cup\lbrace j\rbrace}\right)R'.\]
Enfin, si $J^\times=\emptyset$ alors, $z_q=1$.
\\Pour terminer cette remarque, on va \'etudier le cas o\`u l'\'eclatement encadr\'e est ind\'ependant d'un sous-ensemble. Soit $T\subset J^c$, notons:
\[t=\#(T)\textup{ et } r=n-t.\]
Soient $v=\lbrace v_1,...,v_t\rbrace=u_T$, $w=\lbrace w_1,...,w_r\rbrace=u_{\lbrace 1,...,n\rbrace\setminus T}$ et $u'=(v,w')$ o\`u $w'=\lbrace w'_1,...,w'_r\rbrace$. Pour $1\leqslant q\leqslant r$, on \'ecrit:
\[w'_q=w^{\gamma_q},\]
o\`u $\gamma_q\in\mathbb Z^r$. Alors les $r$ vecteurs $\gamma_1,...,\gamma_r$ forment une matrice de $SL_r(\mathbb Z)$ not\'ee $M_r$. Quitte \`a renum\'eroter les lignes de la matrice $M$, on peut \'ecrire $M$ sous la forme d'une matrice diagonale par blocs o\`u un bloc est $M_r$ et l'autre est $I_t$ la matrice identit\'e de taille $t$.
\\Ainsi, pour tout $\delta\in\mathbb Z^r$, on a:
\[w'^\delta=w^\gamma,\:\gamma=\delta F_r.\]
De m\^eme pour le changement de variables inverse, pour tout $\gamma\in\mathbb Z^r$, on a:
\[w^\gamma=w'^\delta,\:\delta=\gamma F_r^{-1}.\]}
\end{rem}
Nous allons g\'en\'eraliser cette remarque dans le cadre des suites locales encadr\'ees.
\begin{prop}\label{propgenerem}
Consid\'erons une suite locale encadr\'ee de la forme:
\[ \xymatrix{\left( R,u\right)=\left( R^{(0)},u^{(0)}\right) \ar[r]^-{\pi_{0}} & \left( R^{(1)},u^{(1)}\right) \ar[r]^-{\pi_{1}} & \ldots   \ar[r]^-{\pi_{l-1}} & \left( R^{(l)},u^{(l)}\right)}. \]
Pour $0\leqslant i\leqslant l-1$, notons $n_{i+1}$ l'entier correspondant \`a l'entier $n_1$ de la D\'efinition \ref{defeclatloc}, $D_{i+1}$ l'ensemble correspondant \`a $D_1$ et $e_{i+1}=\#(D_{i+1})$.
\\Soient $0\leqslant i<i'\leqslant l$, $q\in\lbrace 1,...,n_{i}\rbrace$, $q'\in\lbrace 1,...,n_{i'}\rbrace$. Alors:
\begin{enumerate}
\item $\exists\:\delta_{q}^{(i',i)}\in\mathbb N^{e_{i}},\:z_{q}^{(i',i)}\in R^{(i')\times}$ tels que $u_q^{(i)}=\left(u_{D_{i'}}^{(i')}\right)^{\delta_{q}^{(i',i)}}z_{q}^{(i',i)}$.
\item Supposons de plus que la suite soit ind\'ependante de $u_T$ avec $T\subset\lbrace 1,...,n\rbrace$ et $u_q^{(i)}\not\in u_T$. Alors $\left(u_{D_{i'}}^{(i')}\right)^{\delta_{q}^{(i',i)}}$ est un mon\^ome uniquement en $u_{D_{i'}}^{(i')}\setminus u_T$.
\item Supposons que pour tout $i''>0$ tel que $i\leqslant i''<i'$, $D_{i''}=\lbrace 1,...,n_{i''}\rbrace$ et $q'\in D_{i'}$. Il existe alors $\gamma_{q'}^{(i,i')}\in\mathbb Z^{n_i}$ tel que $u_{q'}^{(i')}=\left(u^{(i)}\right)^{\gamma_{q'}^{(i,i')}}$.
\item Supposons de plus que la suite soit ind\'ependante de $u_T$ avec $T\subset\lbrace 1,...,n\rbrace$ et $u_{q'}^{(i')}\not\in u_T$. Alors $u_{q'}^{(i')}$ est un mon\^ome uniquement en $u_{\lbrace 1,...,n_i\rbrace}^{(i)}\setminus u_T$.
\end{enumerate}
\end{prop}
\noindent\textit{Preuve}: Il suffit de montrer le cas o\`u $i'=i+1$, le cas g\'en\'eral se faisant par r\'ecurrence. Or ce cas n'est qu'une application des d\'efinitions et de la Remarque \ref{matriceeclat}.\\\qed
\begin{prop}\label{propmatriceeclat}
Consid\'erons les m\^emes hypoth\`eses que dans la Proposition \ref{propgenerem} et supposons de plus que la suite locale encadr\'ee est monomiale et ind\'ependante de $u_T$, $T\subset\lbrace 1,...,n\rbrace$. Notons $t=\#(T)$ et $r=n-t$. On pose:
\[v=\lbrace v_1,...,v_t\rbrace=u_T,\]
\[w=\lbrace w_1,...,w_r\rbrace=u_{\lbrace 1,...,n\rbrace\setminus T}.\]
Alors:
\begin{enumerate}
\item $\forall\:i\in [\![0,l]\!],\:n_i=n$.
\item $\forall\:i\in \:]\!]0,l]\!],\:D_i=\lbrace 1,...,n\rbrace$.
\item Pour $0\leqslant i<i'\leqslant l$, notons $u^{(i)}=\left(v,w^{(i)}\right)$ o\`u $w^{(i)}=\left(w_1^{(i)},...,w_r^{(i)}\right)$ et $u^{(i')}=\left(v,w^{(i')}\right)$ o\`u $w^{(i')}=\left(w_1^{(i')},...,w_r^{(i')}\right)$. Alors, pour tout $1\leqslant q\leqslant r$, $w_q^{(i)}$ est un mon\^ome en $w^{(i')}$ ayant des exposants positifs.
\item Pour $1\leqslant q\leqslant r$, notons $w_q^{(i')}=\left(w^{(i)}\right)^{\gamma_q}$, $\gamma_q\in\mathbb Z^r$. Alors, les $r$ vecteurs colonnes $\gamma_1,...,\gamma_r$ forment une matrice $F_r^{(i',i)}\in SL_r(\mathbb Z)$. R\'eciproquement, notons $w_q^{(i)}=\left(w^{(i')}\right)^{\delta_q}$, $\delta_q\in\mathbb N^r$. Alors, les $r$ vecteurs colonnes $\delta_1,...,\delta_r$ forment la matrice $\left(F_r^{(i',i)}\right)^{-1}\in SL_r(\mathbb Z)$.
\end{enumerate}
\end{prop}
\noindent\textit{Preuve}: Comme dans la preuve de la Proposition \ref{propgenerem}, il suffit de montrer le cas o\`u $i'=i+1$, le cas g\'en\'eral se faisant par r\'ecurrence (et en remarquant que $SL_r(\mathbb Z)$ est un groupe). Or ce cas n'est qu'une application des d\'efinitions et de la Remarque \ref{matriceeclat}.\\\qed
\subsection{Construction d'un \'eclatement local encadr\'e}
~\smallskip ~\\ \indent Gardons les m\^emes notations que dans la sous-section \ref{sectioneclatencad}. Nous d\'efinissons un \'eclatement encadr\'e $\pi:(R,u,k)\rightarrow\left(R^{(1)},u^{(1)},k^{(1)}\right)$ tr\`es utile par la suite. Nous allons d\'ecrire, en termes de g\'en\'erateurs et relations, l'extension de corps $k\hookrightarrow k^{(1)}$ induite par $\pi$.
\\\indent Rappelons que $R'$ est l'anneau:
\[R'=R_j\left[u'_{J\setminus\lbrace j\rbrace}\right].\]
Notons:
\begin{align*}
h&=\#(J),\\h^c&=\#\left(J^c\right)=n-h,\\h^{\times c}&=\#\left(J^{\times c}\right)+1,\\h^\times&=\#\left(J^\times\right)=h-h^{\times c}.
\end{align*}
Quitte \`a renum\'eroter les variables, on peut supposer que:
\begin{align*}
J&=\lbrace 1,...,h\rbrace,\\J^c&=\lbrace h+1,...,n\rbrace,\\j&=1,\\J^{\times c}&=\lbrace 2,...,h^{\times c}\rbrace,\\J^\times&=\lbrace h^{\times c}+1,...,h\rbrace.
\end{align*}
Les changements de variables deviennent alors:
\[ u'_i=\left \{ \begin{array}{ccl}  u_i & \textup{si} & i\in \lbrace 1\rbrace\cup\lbrace h+1,...,n\rbrace \\  \frac{u_i}{u_j} & \textup{si} & i\in\lbrace 2,...,h\rbrace \end{array} \right.\]
Comme on a vu pr\'ec\'edemment, prenons $\mathfrak m'\in Spec(R')$ tel que $u'_{J^c\cup J^{\times c}\cup\lbrace j\rbrace}\subset\mathfrak m'$, ainsi $R^{(1)}=R'_{\mathfrak m'}$ et $\mathfrak m_1=\mathfrak m'R^{(1)}$. De plus, $\mathfrak m=\mathfrak m_1\cap R=\mathfrak m'\cap R$.
\\Pour $1\leqslant i\leqslant n$, notons $z_i\in k^{(1)}$ l'image de $u'_1\in R'$ dans $k^{(1)}$. On remarque alors que:
\[z_i=0,\:\forall\:i\in J^c\cup J^{\times c}\cup\lbrace j\rbrace.\]
\begin{rem}\label{engparjc}
\textup{Notons $\overline R=R'/\mathfrak m R'$ et $\overline{u_i}\in\overline R$ l'image de $u'_i\in R'$ dans $\overline R$, $i\in J\setminus\lbrace j\rbrace$, alors $\overline R=k\left[\overline{u_{J^{\times c}}},\overline{u_{J^\times}}^{\:\pm\:1}\right].$ Les \'el\'ements $\overline{u_{J^{\times c}}}$ et $\overline{u_{J^\times}}^{\:\pm\:1}$ sont alg\'ebriquement ind\'ependants sur $k$, lorsque $R$ est r\'egulier avec $u$ comme syst\`eme r\'egulier de param\`etres. Or, on a les morphismes:
\[R\rightarrow R'\rightarrow R'_{\mathfrak m'}\rightarrow k^{(1)}.\]
En passant modulo $\mathfrak m$, on obtient:
\[k\rightarrow \overline R\rightarrow \overline R_{\overline{\mathfrak m}}\rightarrow k^{(1)},\]
o\`u $\overline{\mathfrak m}=\mathfrak m'/\mathfrak m R'$. On en d\'eduit que $k^{(1)}$ est engendr\'e sur $k$, en tant que corps, par $z_{J^\times}$.}
\end{rem}
Notons $t=deg.tr\left(k^{(1)}\vert k\right)+h^{\times c}$, par la Remarque \ref{engparjc}: \[deg.tr\left(k^{(1)}\vert k\right)\leqslant h^\times.\]
On en d\'eduit les in\'egalit\'es:
\[h^{\times c}\leqslant t\leqslant h^\times+h^{\times c}=h\leqslant n.\]
De plus, on peut supposer que $z_{h^{\times c}+1},...,z_t$ sont alg\'ebriquement ind\'ependants sur $k$ dans $k^{(1)}$, tant que $z_{t+1},...,z_h$ sont alg\'ebriques sur $k\left(z_{h^{\times c}+1},...,z_t\right)$.
\\Pour $t< i\leqslant h$, notons $P_i(X_i)$ le polyn\^ome minimal de $z_i$ sur $k\left(z_{h^{\times c}+1},...,z_{i-1}\right)$. On a l'isomorphisme:
\[k^{(1)}\simeq \dfrac{k\left(z_{h^{\times c}+1},...,z_t\right)\left[X_{t+1},...,X_h\right]}{\left(P_{t+1}(X_{t+1}),...,P_h(X_h)\right)}.\]
Quitte \`a r\'eduire au m\^eme d\'enominateur, pour $t< i\leqslant h$, on peut choisir $P_i\in k\left[z_{h^{\times c}+1},...,z_{i-1}\right]\left[X_i\right]$, mais alors les $P_i$ ne seront plus des polyn\^omes unitaires. Notons:
\[P_i(X_i)=\sum\limits_m p_{i,m}X_i^m,\]
o\`u $p_{i,m}\in k\left[z_{h^{\times c}+1},...,z_{i-1}\right]$, $t< i\leqslant h$. Notons alors $q_{i,m}$ l'\'el\'ement de $R\left[u'_{h^{\times c}+1},...,u'_{i-1}\right]$ obtenu \`a partir de $p_{i,m}$ en rempla\c{c}ant chaque $z_{i'}$ par $u_{i'}$, $t<i'<i$ et en rempla\c{c}ant chaque coefficient de $p_{i,m}$ par un repr\'esentant dans $R$ (on voit $p_{i,m}$ comme un polyn\^ome en $z_{i'}$ \`a coefficients dans $k=R/\mathfrak m$).
\\En particulier, on remarque que $p_{i,m}\equiv q_{i,m}\mod\mathfrak m^{(1)}$. Enfin, notons:
\[Q_i(X)=\sum\limits_m q_{i,m}X^m.\]
Pour $t< i\leqslant h$, comme $P_i(z_i)=0$ dans $k^{(1)}$, on en d\'eduit que:
\[Q_i(u'_i)\in\mathfrak m^{(1)}.\]
\begin{prop}\label{propconstruceclatencad}
Notons $n_1=n-t+h^{\times c}$ et posons le changement de variables suivant:
\[ u_i^{(1)}=\left \{ \begin{array}{lcl}  Q_{i+n-n_1}(u'_{i+n-n_1}) & \textup{si} & h^{\times c}< i\leqslant h-(n-n_1) \\  u'_i & \textup{si} & 1\leqslant i\leqslant j^{\times c}\\  u'_{i+n-n_1} & \textup{si} & h-(n-n_1)< i \leqslant n_1 \end{array} \right.\]
Alors:
\begin{enumerate}
\item $u^{(1)}=\left(u_1^{(1)},...,u_{n_1}^{(1)}\right)$ est un syst\`eme de g\'en\'erateurs de $\mathfrak m^{(1)}$.
\item $\pi:(R,u)\rightarrow\left(R^{(1)},u^{(1)}\right)$ est un \'eclatement local encadr\'e.
\item Si $R$ est r\'egulier avec $u$ pour syst\`eme r\'egulier de param\`etres, alors $u^{(1)}$ est un syst\`eme r\'egulier de param\`etres de $R^{(1)}$.
\end{enumerate}
\end{prop}
\noindent\textit{Preuve}: Nous allons donner une id\'ee de preuve. Pour (1), il suffit de remarquer que, par construction:
\[u_i^{(1)}\in\mathfrak m^{(1)},\:1\leqslant i\leqslant n_1.\]
R\'eciproquement, par la Remarque \ref{matriceeclat}:
\[\mathfrak m R^{(1)}=\left(u_1^{(1)},u_{h+1-(n-n_1)}^{(1)},...,u_{n_1}^{(1)}\right)R^{(1)}\subset \left(u_1^{(1)},...,u_{n_1}^{(1)}\right)R^{(1)}.\]
Rappelons que $\overline{u_2},...,\overline{u_{h^{\times c}}},Q_{t+1}\left(\overline{u_{t+1}}\right),...,Q_{h}\left(\overline{u_h}\right)$ sont les images de $u_2^{(1)},...,u_{n_1}^{(1)}$ dans $k\left(z_{h^{\times c}+1},...,z_t\right)\left[\overline{u_2},...,\overline{u_{h^{\times c}}},\overline{u_{t+1}},...,\overline{u_h}\right]$; en particulier, ce sont des \'el\'ements de $\overline{\mathfrak m}\overline{R}_{\overline{\mathfrak m}}$. Enfin, $\overline{u_2},...,\overline{u_{h^{\times c}}},Q_{t+1}\left(\overline{u_{t+1}}\right),...,Q_{h}\left(\overline{u_h}\right)$ engendrent un id\'eal maximal de $k\left(z_{h^{\times c}+1},...,z_t\right)\left[\overline{u_2},...,\overline{u_{h^{\times c}}},\overline{u_{t+1}},...,\overline{u_h}\right]$ vu que:
\[\dfrac{k\left(z_{h^{\times c}+1},...,z_t\right)\left[\overline{u_2},...,\overline{u_{h^{\times c}}},\overline{u_{t+1}},...,\overline{u_h}\right]}{\left(\overline{u_2},...,\overline{u_{h^{\times c}}},Q_{t+1}\left(\overline{u_{t+1}}\right),...,Q_{h}\left(\overline{u_h}\right)\right)}\simeq\dfrac{k\left(z_{h^{\times c}+1},...,z_{h^{\times c}}\right)\left[\overline{u_{t+1}},...,\overline{u_h}\right]}
{\left(Q_{t+1}\left(\overline{u_{t+1}}\right),...,Q_{h}\left(\overline{u_h}\right)\right)}\simeq k^{(1)}.\]
Comme:
\begin{align*}
\overline{R}_{\overline{\mathfrak m}}&\simeq k\left[\overline{u_2},...,\overline{u_h}\right]_{\overline{\mathfrak{m}}}
\\&\simeq k\left(z_{h^{\times c}+1},...,z_t\right)\left[\overline{u_2},...,\overline{u_{h^{\times c}}},\overline{u_{t+1}},...,\overline{u_h}\right]_{\overline{\mathfrak m}k\left(z_{h^{\times c}+1},...,z_t\right)\left[\overline{u_2},...,\overline{u_{h^{\times c}}},\overline{u_{t+1}},...,\overline{u_h}\right]}.
\end{align*}
Tout ceci montre que les images de $u_{2}^{(1)},...,u_{n_1}^{(1)}$ engendrent l'id\'eal maximal $\overline{\mathfrak m}\overline{R}_{\overline{\mathfrak m}}$ de $\overline{R}_{\overline{\mathfrak m}}$. Or par d\'efinition de $\overline{R}$ et de ${\overline{\mathfrak m}}$, on en d\'eduit que $u_1^{(1)},...,u_{n_1}^{(1)}$ engendrent l'id\'eal $\mathfrak m'R_{\mathfrak m'}\equiv\mathfrak m^{(1)}$ dans $R'_{\mathfrak m'}\equiv R^{(1)}$.
\\Par d\'efinition, (2) est \'evidente, l'ensemble $D_1$ \'etant:
\[D_1=\lbrace 1,...,h^{\times c}\rbrace\cup\lbrace h-(n-n_1)+1,...,n_1\rbrace.\]
Pour montrer (3), on remarque que, $R$ \'etant r\'egulier avec $u$ comme syst\`eme r\'egulier de param\`etres, alors, $\overline R$ est r\'egulier et $\overline{u_2},...,\overline{u_{h^{\times c}}},Q_{t+1}\left(\overline{u_{t+1}}\right),...,Q_{h}\left(\overline{u_h}\right)$ forment un syst\`eme r\'egulier de param\`etres de l'anneau local r\'egulier $\overline{\mathfrak m}\overline{R}_{\overline{\mathfrak m}}$ qui est de dimension $h-(n-n_1)-1$. Enfin, on montre par r\'ecurrence sur $n-h$ que:
\[(0)\subsetneq \left(u_1^{(1)}\right)\subsetneq\left(u_1^{(1)},u_{h-(n-n_1)+1}^{(1)}\right)\subsetneq ... \subsetneq \left(u_1^{(1)},u_{h-(n-n_1)+1}^{(1)},...,u_{n_1}^{(1)}\right)\]
forme une cha\^ine de $n-h+1$ id\'eaux premiers de $R^{(1)}$ distincts.\\\qed
\\ \\\indent Pour terminer, nous allons interpr\'eter les r\'esultats pr\'ec\'edents en termes d'\'eclatements encadr\'es par rapport \`a une valuation donn\'ee.
\\\indent Soient $(R,\mathfrak m,k)$ un anneau local noeth\'erien, $u$ un ensemble de g\'en\'erateurs de $\mathfrak m$ et $\nu$ une valuation centr\'ee en $R$. Pour $1\leqslant i\leqslant n$, notons:
\[\beta_i=\nu(u_i),\]
\[x_i=in_\nu(u_i).\]
Soient $T\subset\lbrace 1,...,n\rbrace$, $E=\lbrace 1,...,n\rbrace\setminus T$ et $k\left[x_E\right]$ la sous-alg\`ebre gradu\'ee de $G_\nu$. Notons:
\[G=k\left[x_E\right]^{*}=\left\lbrace \left.\dfrac{f}{g}\:\right|\:f,g\in k\left[x_E\right],\:g\textup{ homog\`ene},\:g\neq 0\right\rbrace.\]
Consid\'erons $J\subset E$ et choisissons $j\in J$ tel que:
\[\beta_j=\min_{i\in J}\lbrace \beta_i\rbrace.\]
Soit $\pi:(R,\mathfrak m)\rightarrow\left(R^{(1)},\mathfrak m^{(1)}\right)$ un \'eclatement local par rapport \`a $\nu$ (voir D\'efinition \ref{defeclatlocal}) et consid\'erons $R^{(1)}$ comme le localis\'e de $R'$ en le centre de $\nu$. On a donc:
\[J^{\times c}=\lbrace i\in J\:\vert\:\beta_i>\beta_j\rbrace,\]
\[J^\times=\lbrace i\in J\setminus\lbrace j\rbrace\:\vert\:\beta_i=\beta_j\rbrace.\]
Notons alors:
\[ x'_i=\left \{ \begin{array}{ccl}  x_i & \textup{si} & i\in J^c\cup\lbrace j\rbrace \\  \frac{x_i}{x_j} & \textup{si} & i\in J\setminus\lbrace j\rbrace \end{array} \right.\]
\[ \overline{z}_i=\left \{ \begin{array}{ccl}  u'_i & \textup{si} & i\in J^\times \\  1 & \textup{si} & i\in \lbrace 1,...,n\rbrace\setminus J^\times \end{array} \right.\]
\[ z_i=\left \{ \begin{array}{ccl}  x'_i & \textup{si} & i\in J^\times \\  1 & \textup{si} & i\in \lbrace 1,...,n\rbrace\setminus J^\times \end{array} \right.\]
\[\beta'_i=ord(x'_i), 1\leqslant i\leqslant n,\]
\[E'=E\setminus J^\times.\]
Pour $1\leqslant i\leqslant n$, $x'_i$ est homog\`ene et $ord(x'_i)\geqslant 0$. On a:
\[ord(x'_i)>0\Leftrightarrow\beta_i>\beta_j\Leftrightarrow i\in E',\]
\[ord(z_i)=0,\:\forall\:i\in J^\times.\]
Remarquons que $k^{(1)}=k\left(z_{J^\times}\right)$. Consid\'erons le morphisme $\rho:R'\rightarrow k^{(1)}$, extension de $R\rightarrow k$, d\'efini en envoyant $u'_i$ sur $z_i$ si $z\in J^\times$ et sur $0$ si $i\in J^{\times c}$. L'id\'eal $\mathfrak m'=\ker \rho$ est le centre de $\nu$ dans $R'$ et $R^{(1)}=R'_{\mathfrak m'}$.
\begin{defi}\label{eclatlocparrapnu}
Consid\'erons $u^{(1)}$ comme dans la Proposition \ref{propconstruceclatencad}, l'\'eclatement local encadr\'e $\pi:(R,\mathfrak m)\rightarrow\left(R^{(1)},\mathfrak m^{(1)}\right)$ qui en r\'esulte est appel\'e l'\textbf{\'eclatement local encadr\'e le long de} $\boldsymbol{\left(u_J\right)}$ \textbf{par rapport \`a} $\boldsymbol{\nu}$.
\end{defi}
Soit $\varphi:D_1\rightarrow J^c\cup J^{\times c}\cup\lbrace j\rbrace$ la bijection r\'esultante de l'\'eclatement encadr\'e. Notons alors:
\[E^{(1)}=\varphi^{-1}(E')\subset D_1,\]
\[x_i^{(1)}=x'_{\varphi(i)},\]
\[\beta_i^{(1)}=ord\left(x_i^{(1)}\right)=\beta'_{\varphi(i)}.\]
\begin{rem}\label{egalitealggrad}
\textup{Pour tout $i,i'\in E$, il existe $\delta_i\in\mathbb N^{\#(E^{(1)})}$, $\gamma_{i'}\in\mathbb Z^{\#(E)}$ tels que:
\[u_{i'}^{(1)}=u_E^{\gamma_{i'}},\]
\[u_i=\left(u_{E^{(1)}}^{(1)}\right)^{\delta_i}\overline{z_i}.\]
On a les m\^emes transformations au niveau des alg\`ebres gradu\'ees:
\[x_{i'}^{(1)}=x_E^{\gamma_{i'}},\]
\[x_i=\left(x_{E^{(1)}}^{(1)}\right)^{\delta_i}z_i.\]
On a \'egalement:
\[\beta_i=\left\langle\delta_i,\beta_{E^{(1)}}^{(1)}\right\rangle,\]
o\`u $\langle .,.\rangle$ repr\'esente le produit scalaire de vecteurs de taille $\#(D_1)$.
On en d\'eduit les \'egalit\'es d'alg\`ebres gradu\'ees suivantes:
\[k[x]^*=k\left[z_{J^\times},x_{D_1}^{(1)}\right]^*=k_1\left[x_{D_1}^{(1)}\right]^*,\]
\[k[x_E]^*=k\left[z_{J^\times},x_{E^{(1)}}^{(1)}\right]^*=k_1\left[x_{E^{(1)}}^{(1)}\right]^*.\]}
\end{rem}
\subsection{L'id\'eal premier implicite}
~\smallskip ~\\ \indent Pour une valuation donn\'ee, l'id\'eal premier implicite est un des objets central de l'uniformisation locale. En effet, cet id\'eal va \^etre l'id\'eal \`a d\'esingulariser. C'est un id\'eal du compl\'et\'e qui d\'ecrit les \'el\'ements de valuation infinie. Via le Lemme \ref{lemmetechniquespiva}, pour rendre $R$ r\'egulier, il nous suffit de rendre r\'egulier $\widehat{R}_H$ et $\widehat{R}/H$, o\`u $H$ est l'id\'eal premier implicite associ\'e \`a une valuation centr\'ee en $R$. L'int\'er\^et de l'id\'eal premier implicite est que $\widehat{R}_H$ est automatiquement r\'egulier sous l'hypoth\`ese de quasi-excellence. Ainsi il suffira de d\'emontrer l'uniformisation locale pour des valuations centr\'ees en $\widehat{R}/H$.
\\On pourra consulter \cite{ideimp}, \cite{jcthese} ou \cite{spiva} pour plus de d\'etails.
\begin{defi}(\cite{ideimp}, Definition 2.1)\label{defiidim}
Soient $(R,\mathfrak{m},k)$ un anneau local noeth\'erien, $P_\infty$ un id\'eal premier minimal de $R$ et $\nu$ une valuation de rang $1$ de $R_{P_\infty}$ centr\'ee en $R$. On appelle \textbf{id\'eal premier implicite} de $R$ associ\'e \`a $\nu$, not\'e $H(R,\nu)$ ou plus simplement $H$ s'il n'y a pas d'ambigu\"it\'e, l'id\'eal de $\widehat{R}$ d\'efinit par:
\[H=\bigcap\limits_{\beta\in\nu(R\setminus P_\infty)}P_\beta\widehat{R},\]
o\`u $P_\beta=\lbrace f\in R\:\vert\:\nu(f)\geqslant\beta\rbrace$.
\end{defi}
\begin{rem}\label{remcauchy}
\textup{\begin{enumerate}
\item Si l'on suppose de plus que $R$ est int\`egre, alors $P_\infty=(0)$.
\item Comme la valuation est de rang $1$ le groupe $\nu(R\setminus P_\infty)$ est archim\'edien (voir preuve du Lemme \ref{pasdesuitecroissantebornee}) et donc, pour tout $\beta\in\nu(R\setminus P_\infty)$, il existe $n\in\mathbb N$ tel que $\mathfrak m^n\subset P_\beta$. Il y a donc \'equivalence entre:
\begin{enumerate}
\item $f\in H$;
\item Il existe une suite de Cauchy $(f_n)_n\subset R$ telle que, si $f_n\underset{n\rightarrow +\infty}{\longrightarrow}f$, alors $\nu(f_n)\underset{n\rightarrow +\infty}{\longrightarrow}\infty$;
\item Pour toute suite de Cauchy $(f_n)_n\subset R$ telle que, si $f_n\underset{n\rightarrow +\infty}{\longrightarrow}f$, alors $\nu(f_n)\underset{n\rightarrow +\infty}{\longrightarrow}\infty$.
\end{enumerate}
\end{enumerate}}
\end{rem}
\begin{lem}\label{intersecimp}
Sous les hypoth\`eses de la D\'efinition \ref{defiidim}, si $H$ est l'id\'eal premier implicite de $R$ associ\'e \`a $\nu$, alors:
\[H\cap R=P_\infty\]
et il existe une inclusion naturelle:
\[R/P_\infty\hookrightarrow\widehat{R}/H.\]
\end{lem}
\begin{thm}(\cite{ideimp}, Theorem 2.1)\label{valetenuniqidealimpl}
Reprenons les hypoth\`eses de la D\'efinition \ref{defiidim}. Soit $H$ l'id\'eal premier implicite de $R$ associ\'e \`a $\nu$, alors:
\begin{enumerate}
\item $H$ est un id\'eal premier de $\widehat{R}$;
\item $\nu$ s'\'etend de mani\`ere unique en une valuation $\widehat{\nu}$ centr\'ee en $\widehat{R}/H$.
\end{enumerate}
\end{thm}
\begin{coro}(\cite{ideimp}, Proposition 2.8)\label{RHestreg}
Soient $R$ un anneau local quasi-excellent r\'eduit, $P_\infty$ un id\'eal premier minimal de $R$ et $\nu$ une valuation de rang $1$ de $R_{P_\infty}$ centr\'ee en $R$. Alors, $\widehat R_H$ est un anneau local r\'egulier.
\end{coro}
\begin{lem}(\cite{ideimp}, Lemma 2.4)\label{implistable}
Soit $(R,\mathfrak{m})\rightarrow (R',\mathfrak m')$ un morphisme local d'anneaux locaux noeth\'eriens. Soient $P_\infty$ un id\'eal premier minimal de $R$ et $\nu$ une valuation de rang $1$ de $R_{P_\infty}$ centr\'ee en $R$. Supposons qu'il existe un id\'eal premier minimal $P'_\infty$ de $R'$ tel que $P_\infty=P'_\infty\cap R$ et que $\nu$ s'\'etende en une valuation de rang $1$ $\nu'$ telle que son groupe des valeurs contienne celui de $\nu$.
\\Pour $\beta\in\nu'\left(R'\setminus\lbrace 0\rbrace\right)$, notons $P'_\beta=\lbrace f\in R'\:\vert\:\nu'(f)\geqslant\beta\rbrace$. Enfin, notons $H'=H(R',\nu')$ l'id\'eal premier implicite de $R'$ associ\'e \`a $\nu'$.
\\Alors, pour tout $\beta\in\nu\left(R\setminus\lbrace 0\rbrace\right)$,
\[\left(P'_\beta\widehat{R'}\right)\cap\widehat R=P_\beta\widehat R.\]
\end{lem}
\begin{coro}(\cite{ideimp}, Corollary 2.5)\label{impeclatloc}
Avec les hypoth\`eses du Lemme \ref{implistable}, on a:
\[H'\cap\widehat R=H.\]
\end{coro}
\begin{coro}(\cite{ideimp}, Corollary 2.7)\label{coroht}
Reprenons les hypoth\`eses du Lemme \ref{implistable} et supposons de plus que $(R,\mathfrak{m})\rightarrow (R',\mathfrak m')$ est un \'eclatement local par rapport \`a $\nu$ le long d'un id\'eal $J$ non nul de $R$ et que $\nu$ reste de rang $1$ sur $R'$. On a:
\[ht(H')\geqslant ht(H)\textit{ et}\]
\[\dim\left(\widehat{R'}/H'\right)\leqslant\dim\left(\widehat{R}/H\right).\]
\end{coro}
\subsection{Monomialisation d'\'el\'ements non d\'eg\'en\'er\'es}
~\smallskip ~\\ \indent Nous allons voir l'effet des \'eclatements encadr\'es sur les mon\^omes. Une cons\'equence sera qu'un \'el\'ement non d\'eg\'en\'er\'e, c'est-\`a-dire qu'en cet \'el\'ement, la valuation est \'egale \`a la valuation monomiale, peut \^etre transform\'e en un mon\^ome via une suite locale encadr\'ee.
\\Toute cette partie est un cas particulier du jeu d'Hironaka (voir \cite{jeuhiro} et \cite{spivajeuhiro}).
\\\indent Pour un \'el\'ement $\alpha=(\alpha_1,...,\alpha_n)\in \mathbb N^n$, on note:
\[\vert\alpha\vert=\alpha_1+...+\alpha_n.\]
Soient $\alpha=(\alpha_1,...,\alpha_n)$, $\gamma=(\gamma_1,...,\gamma_n)\in \mathbb N^n$. Pour $1\leqslant i\leqslant n$, notons:
\[\delta_i=\min\lbrace\alpha_i,\gamma_i\rbrace.\]
Posons alors $\delta=(\delta_1,...,\delta_n)\in \mathbb N^n$, $\tilde\alpha=\alpha-\delta$, $\tilde\gamma=\gamma-\delta$. Quitte \`a \'echanger $\alpha$ et $\gamma$, on peut supposer que $\vert\tilde\alpha\vert\leqslant\vert\tilde\gamma\vert$. On d\'efinit $\tau(\alpha,\gamma)$ par:
\[\tau(\alpha,\gamma)=\left(\vert\tilde\alpha\vert ,\vert\tilde\gamma\vert\right).\]
\begin{rem}\label{remtau}
\textup{
\begin{enumerate}
\item Si $\tilde\alpha=(0,...,0)$, alors $u^\alpha$ divise $u^\gamma$ dans $R$.
\item Quitte \`a renum\'eroter les variables de $\tilde\alpha$ et $\tilde\gamma$, on peut supposer qu'il existe $a\in\mathbb N$, $1\leqslant a<n$, tel que:
\[\tilde\alpha=(\tilde\alpha_1,...,\tilde\alpha_a,\underbrace{0,...,0}_{n-a}),\]
\[\tilde\gamma=(\underbrace{0,...,0}_{a},\tilde\gamma_{a+1},...,\tilde\gamma_n).\]
On peut \'egalement supposer que, pour $1\leqslant i\leqslant a$, $\tilde\alpha_i>0$.
\end{enumerate}}
\end{rem}
\indent Soit $(R,\mathfrak m)$ un anneau local noeth\'erien tel que $\mathfrak m$ soit non nilpotent et $u=(u_1,...,u_n)$ un ensemble de g\'en\'erateurs de $\mathfrak m$. Soit $\nu$ une valuation centr\'ee en $R$ de groupe des valeurs $\Gamma$.
\\Consid\'erons $J\subset\lbrace 1,...,n\rbrace$ le sous-ensemble le plus petit possible, au sens de l'inclusion, tel que:
\[\lbrace 1,...,a\rbrace\subset J\textup{ et }\sum\limits_{i\in J}\tilde\gamma_i\geqslant\vert\tilde\alpha\vert.\]
En reprenant les notations de la D\'efinition \ref{defeclatloc}, consid\'erons $\pi:(R,u)\rightarrow\left(R^{(1)},u^{(1)}\right)$ un \'eclatement encadr\'e le long de $\left(u_J\right)$, selon la D\'efinition \ref{eclatlocparrapnu}. Notons :
\[ \tilde\alpha'_i=\left \{ \begin{array}{ccc}  \tilde\alpha_i & \textup{si} & i\neq j \\  0 & \textup{si} & i=j \end{array} \right.\]
\[ \tilde\gamma'_i=\left \{ \begin{array}{ccc}  \tilde\gamma_i & \textup{si} & i\neq j \\  \sum\limits_{i\in J}\tilde\gamma_i-\vert\tilde\alpha\vert & \textup{si} & i=j \end{array} \right.\]
\[\tilde\alpha'=(\tilde\alpha'_1,...,\tilde\alpha'_n),\]
\[\tilde\gamma'=(\tilde\gamma'_1,...,\tilde\gamma'_n),\]
\[\delta'=(\delta_1,...,\delta_{j-1},\delta_j+\vert\tilde\alpha\vert,\delta_{j+1},...,\delta_n).\]
Avec ces notations on obtient:
\[u^\alpha=\left(u'\right)^{\delta'+\tilde\alpha'},\]
\[u^\gamma=\left(u'\right)^{\delta'+\tilde\gamma'}.\]
Posons $\alpha'=\delta'+\tilde\alpha'$ et $\gamma'=\delta'+\tilde\gamma'$.
\begin{prop}\label{propdesctau}
Avec les notations pr\'ec\'edentes, on a, pour l'ordre lexicographique:
\[\tau(\alpha',\gamma')<_{lex}\tau(\alpha,\gamma).\]
\end{prop}
\noindent\textit{Preuve}: Nous ne donnerons qu'une id\'ee de preuve. Il suffit de montrer que:
\[\left(\vert\tilde\alpha'\vert,\vert\tilde\gamma'\vert\right)<_{lex}\left(\vert\tilde\alpha\vert,\vert\tilde\gamma\vert\right).\]
Si $j\in\lbrace 1,...,a\rbrace$, alors par d\'efinition et par la Remarque \ref{remtau}, on a:
\[\vert\tilde\alpha'\vert=\vert\tilde\alpha\vert-\tilde\alpha_j<\vert\tilde\alpha\vert.\]
Si $j\in\lbrace a+1,...,n\rbrace$, alors $\vert\tilde\alpha'\vert=\vert\tilde\alpha\vert$. Par minimalit\'e de $J$, il vient que:
\[\sum\limits_{i\in J\setminus\lbrace j\rbrace}\tilde\gamma_i<\vert\tilde\alpha\vert.\]
On en conclut que $\vert\tilde\gamma'\vert<\vert\tilde\gamma\vert$.\\\qed
\begin{coro}\label{procesusstop}
Soit $s=\#\lbrace i\in\lbrace 1,...,n\rbrace\:\vert\:u'_i\not\in R^{(1)}{}^\times\rbrace$. Quitte \`a renum\'eroter les variables, on peut supposer que $u'_i$ n'est pas inversible dans $R^{(1)}$, pour $1\leqslant i\leqslant s$ et, inversible pour $s< i\leqslant n$.
\\Comme $\pi$ est un \'eclatement encadr\'e, $\lbrace u'_1,...,u'_s\rbrace\subset u^{(1)}$. Quitte \`a renum\'eroter les variables, on peut supposer que $u'_i=u_i^{(1)}$, $1\leqslant i\leqslant s$. Notons les vecteurs de taille $n_1$ par:
\[\alpha^{(1)}=(\tilde\alpha'_1,...,\tilde\alpha'_s,\underbrace{0,...,0}_{n_1-s}),\]
\[\tilde\gamma=(\gamma'_{1},...,\gamma'_s,\underbrace{0,...,0}_{n_1-s}).\]
Alors, on a:
\[\tau\left(\alpha^{(1)},\gamma^{(1)}\right)<_{lex}\tau(\alpha,\gamma).\]
\end{coro}
\noindent\textit{Preuve}: Par la Proposition \ref{propdesctau} et par d\'efinition:
\[\tau\left(\alpha^{(1)},\gamma^{(1)}\right)\leqslant_{lex}\tau(\alpha',\gamma')<_{lex}\tau(\alpha,\gamma).\]\qed
\begin{rem}\label{remindeptau}
\textup{Soit $T\subset\lbrace 1,..,n \rbrace$ tel que $\tilde\alpha_i=\tilde\gamma_i=0$, pour tout $i\in T$. Alors, tout \'eclatement encadr\'e le long de $(u_J)$, avec $J$ d\'efini comme pr\'ec\'edemment, est ind\'ependant de $u_T$.}
\end{rem}
\begin{coro}\label{coroexistsuitlocale}
Soient $(R,\mathfrak m)$ un anneau local noeth\'erien tel que $\mathfrak m$ soit non nilpotent et $u=(u_1,...,u_n)$ un ensemble de g\'en\'erateurs de $\mathfrak m$. 
Soit $r\in\mathbb N$ tel que $1\leqslant r\leqslant n$. Notons $u=(w,v)$ avec:
\[w=(w_1,...,w_r)=(u_1,...,u_r),\]
\[v=(v_1,...,v_{n-r}).\]
Soit $\nu$ une valuation centr\'ee en $R$, prenons $j$ dans $J$ v\'erifiant:
\[\nu(u_j)=\min_{i\in J}\lbrace\nu(u_i)\rbrace.\]
Soient $\alpha,\gamma\in\mathbb N^{r}$, il existe alors  une suite locale encadr\'ee par rapport \`a $\nu$ et ind\'ependante de $v$:
\[(R,u)\rightarrow \left(R^{(l)},u^{(l)}\right)\]
telle que $w^\alpha$ divise $w^\gamma$ ou bien $w^\gamma$ divise $w^\alpha$ dans $R^{(l)}$.
\end{coro}
\noindent\textit{Preuve}: On it\`ere le processus de construction de la Proposition \ref{propdesctau} en choisissant des \'eclatements locaux encadr\'es par rapport \`a $\nu$, qui sont, par construction et par la Remarque \ref{remindeptau}, ind\'ependants de $v$. Par le Corollaire \ref{procesusstop}, cette construction s'arr\^ete apr\`es un nombre fini d'it\'erations. On conclut alors gr\^ace au (1) de la Remarque \ref{remtau}.\\\qed
\begin{prop}\label{quidivisequi}
Gardons les notations du Corollaire \ref{coroexistsuitlocale}. Alors:
\[w^\alpha\textit{ divise }w^\gamma \textit{ dans }R^{(l)}\Leftrightarrow\nu\left(w^\alpha\right)\leqslant\nu\left(w^\gamma\right).\]
\end{prop}
\noindent\textit{Preuve}: Notons $u^{(l)}=\left(w_1^{(l)},...,w_{r_l}^{(l)},v\right)$. Par le (1) de la Proposition \ref{propgenerem}, il existe $\alpha^{(l)},\gamma^{(l)}\in\mathbb N^{r_l}$ et $y,z\in R^{(l)}{}^\times$ tels que:
\[w^\alpha=y\left(w^{(l)}\right)^{\alpha^{(l)}},\]
\[w^\gamma=z\left(w^{(l)}\right)^{\gamma^{(l)}}.\]
Comme $\nu\left(w_1^{(l)}\right),...,\nu\left(w_{r_l}^{(l)}\right)\geqslant 0$ et comme, par construction, l'un des $\alpha^{(l)}$ ou $\gamma^{(l)}$ est plus grand que l'autre, composante par composante, on a:
\[\left(w^{(l)}\right)^{\alpha^{(l)}}\textit{ divise }\left(w^{(l)}\right)^{\gamma^{(l)}} \textit{ dans }R^{(l)}\Leftrightarrow\nu\left(\left(w^{(l)}\right)^{\alpha^{(l)}}\right)\leqslant\nu\left(\left(w^{(l)}\right)^{\gamma^{(l)}}\right).\]\qed
\begin{coro}\label{coroideal}
Gardons les notations du Corollaire \ref{coroexistsuitlocale}. Soit $I$ un id\'eal de $R$ engendr\'e par des mon\^omes en $w$. Consid\'erons $\varepsilon_0,...,\varepsilon_b\in\mathbb N^{r}$ une collection minimale d'\'el\'ements de $\mathbb N^{r}$ telle que $\left(w^{\varepsilon_0},...,w^{\varepsilon_b}\right)=I$.
\\Enfin, supposons que $\nu\left(w^{\epsilon_0}\right)\leqslant\nu\left(w^{\epsilon_i}\right)$, $1\leqslant i\leqslant b$. Il existe alors une suite locale encadr\'ee par rapport \`a $\nu$ et ind\'ependante de $v$:
\[(R,u)\rightarrow \left(R^{(l)},u^{(l)}\right)\]
telle que:
\[I R^{(l)}=\left(w^{\varepsilon_0}\right)R^{(l)}.\]
\end{coro}
\noindent\textit{Preuve}: On d\'efinit l'entier suivant:
\[\tau(I,w)=\left(b,\min_{0\leqslant i<i'\leqslant b}\lbrace \tau\left(w^{\varepsilon_i},w^{\varepsilon_{i'}}\right)\rbrace\right).\]
On suppose que $\tau(I,w)=(0,1)$ si $b=0$. Si $b\geqslant 1$, on applique la Proposition \ref{propdesctau} \`a la paire $\left\lbrace w^{\varepsilon_i},w^{\varepsilon_{i'}}\right\rbrace$, pour laquelle le minimum est atteint dans $\lbrace \tau\left(w^{\varepsilon_i},w^{\varepsilon_{i'}}\right)\rbrace$. On obtient alors un sous-ensemble $J$ de $\lbrace 1,...,n\rbrace$ tel que tout \'eclatement encadr\'e le long de $(u_J)$ fasse d\'ecro\^itre $\tau(I,w)$ pour l'ordre lexicographique. On conclut en utilisant la Proposition \ref{quidivisequi}.\\\qed
\begin{lem}\label{lemvalmonomiale}
Soit $(R,\mathfrak{m},k)$ un anneau local r\'egulier. On suppose que $\mathfrak{m}=(u_1,...,u_n)=u$, o\`u $n$ est le nombre de g\'en\'erateurs de $\mathfrak m$. Soient $\Phi$ un semi-groupe ordonn\'e archim\'edien et $\beta_1,...,\beta_n\in\Phi$ tels que $\beta_i> 0$, $1\leqslant i\leqslant n$.
\\Notons $\Phi_*\subset\Phi$ le semi-groupe ordonn\'e suivant:
\[\Phi_*=\left\lbrace \left.\sum\limits_{i=1}^n\alpha_i\beta_i\:\right|\:\alpha_i\in\mathbb N\right\rbrace.\]
Pour $\gamma\in\Phi_*$, consid\'erons l'id\'eal de $R$:
\[I_\gamma=\left\langle\left\lbrace u_1^{\alpha_1}...u_n^{\alpha_n}\:\left|\:\sum\limits_{i=1}^n\alpha_i\beta_i\geqslant\gamma\right.\right\rbrace\right\rangle.\]
Alors, pour $f\in R\setminus\lbrace 0\rbrace$, l'ensemble:
\[\Phi_f=\lbrace \gamma\in\Phi_*\:\vert\:f\in I_\gamma\rbrace\]
est fini.
\end{lem}
\noindent\textit{Preuve}: Soit $f\in R\setminus\lbrace 0\rbrace$. Comme $\Phi$ est archim\'edien alors $\Phi_*$ l'est aussi. Remarquons que $\Phi_*$ est un ensemble d\'enombrable et que $\Phi_f$ est un ensemble bien ordonn\'e car \`a une suite d\'ecroissante de $\Phi_f$ correspond une suite croissante d'id\'eaux de $R$ de la forme $I_\gamma$ qui est forc\'ement finie vu que $R$ est noeth\'erien. Notons $\gamma_0$ le plus petit \'el\'ement non nul de $\Phi_f$ (en fait $0$ est le min de $\Phi_*$ et de $\Phi_f$, si ce dernier ensemble est r\'eduit \`a $0$, la preuve est termin\'ee). 
\\Comme $f$ est non nul, il existe $i\geqslant 0$ tel que $f\notin\mathfrak m^i$ et donc $\Phi_f\neq\Phi_*$. Il existe donc $\gamma_1=\sup \Phi_f\in\Phi_*$. Or $\Phi_*$ est archim\'edien, ainsi, il existe $N\in\mathbb N$ tel que $\gamma_1 < N\gamma_0$. Alors, pour tout \'el\'ement $\gamma=\sum\limits_{i=1}^n\alpha_i\beta_i\in\Phi_f$, $\alpha_i\in\mathbb N$, comme $\beta_i\in\Phi_f$, pour $1\leqslant i\leqslant n$, on en d\'eduit:
\[\left(\sum\limits_{i=1}^n\alpha_i\right)\gamma_0\leqslant\sum\limits_{i=1}^n\alpha_i\beta_i <\gamma_1 < N\gamma_0.\] 
N\'ecessairement, on a $\left(N-\sum\limits_{i=1}^n\alpha_i\right)\gamma_0>0$ et comme $\gamma_0>0$ on en d\'eduit que $\sum\limits_{i=1}^n\alpha_i\leqslant N$, c'est-\`a-dire qu'il n'y a qu'un choix fini de $n$-uplets $(\alpha_1,...,\alpha_n)$ et donc de $\gamma\in\Phi_f$. \\\qed
\begin{coro}\label{defivaluationmono}
Sous les hypoth\`eses du Lemme \ref{lemvalmonomiale}, il existe une unique valuation, not\'ee $\nu_{0,u}$, centr\'ee en un id\'eal premier de $R$, telle que:
\[\nu_{0,u}(u_j)=\beta_j,\:1\leqslant j\leqslant n;\]
\[\nu_{0,u}(f)=\max \lbrace\gamma\in\Phi_f\rbrace,\:\forall\:f\in R\setminus\lbrace 0\rbrace.\]
Cette valuation est appel\'ee la \textbf{valuation monomiale} de $R$ associ\'ee \`a $u$ et \`a $\beta_1,...,\beta_n$.
\\Soit $\nu$ une valuation de groupe des valeurs $\Gamma$ et centr\'ee en un id\'eal premier de $R$. On dit que $\nu$ est \textbf{monomiale} par rapport \`a $u$ s'il existe $\beta_1,...,\beta_n\in\Gamma_+$ tels que:
\[\forall\:f\in R\setminus\lbrace 0\rbrace,\:\nu(f)=\nu_{0,u}(f).\]
\end{coro}
\begin{ex}
\textup{Pour $\nu$ une valuation centr\'ee en $R=k\left[\left[u_1,...,u_n\right]\right]$, si $f=\sum c_{\alpha}u^\alpha$, alors $\nu_{0,u}(f)=\min\left\lbrace\left.\sum\limits_{i=1}^n\alpha_i\nu(u_i)\:\right|\:c_\alpha\neq 0\right\rbrace$ est la valuation monomiale associ\'ee \`a $u$ et \`a $\nu(u_1),...,\nu(u_n)$.}
\end{ex}
\begin{rem}\label{valmonpluspetite}
\textup{Si $\nu$ est une valuation centr\'ee en $R$ dont le groupe des valeurs est archim\'edien et si $\nu_{0,u}$ est la valuation monomiale associ\'ee \`a $u$ et \`a $\nu(u_1),...,\nu(u_n)$, alors, pour tout $\gamma\in\Phi_*$, $\nu(I_\gamma)=\min\lbrace \nu(f)\:\vert\:f\in I_\gamma\rbrace\geqslant\gamma$. Ainsi, pour tout $f\in R\setminus\lbrace 0\rbrace$:
\[\nu_{0,u}(f)\leqslant\nu(f).\]
De plus, la valuation $\nu$ est monomiale si et seulement si:
\[gr_\nu(R)=k\left[in_\nu(u_1),...,in_\nu(u_n)\right].\]}
\end{rem}
\begin{defi}\label{defnondeg}
Soient $R$ un anneau local r\'egulier et $u=(u_1,...,u_n)$ un syst\`eme r\'egulier de param\`etres de $R$. Soit $\nu$ une valuation centr\'ee en $R$. On dit que $f\in R$ est \textbf{non d\'eg\'en\'er\'e} par rapport \`a $\nu$ et $u$ si:
\[\nu_{0,u}(f)=\nu(f),\]
o\`u $\nu_{0,u}$ est la valuation monomiale de $R$ par rapport \`a $u$ (Corollaire \ref{defivaluationmono}).
\end{defi}
\begin{rem}
\textup{\begin{enumerate}
\item $f\in R$ est non d\'eg\'en\'er\'e par rapport \`a $u$ si et seulement s'il existe un id\'eal $I$ de $R$, monomial par rapport \`a $u$, tel que $\nu(f)=\min\limits_{g\in I}\lbrace\nu(g)\rbrace$.
\item Consid\'erons une suite locale encadr\'ee $(R,u)\rightarrow\left(R^{(1)},u^{(1)}\right)$ et $f\neq 0$. Par le (1) de la Proposition \ref{propgenerem}, chaque $u_j$ est un mon\^ome en $u^{(1)}$ multipli\'e par une unit\'e de $R^{(1)}$. Ainsi, si $f$ est non d\'eg\'en\'er\'e par rapport \`a $u$ alors, $f$ est aussi non d\'eg\'en\'er\'e par rapport \`a $u^{(1)}$.
\end{enumerate}}
\end{rem}
Le Th\'eor\`eme \ref{uniflocnondeg} suivant peut \^etre vu comme un th\'eor\`eme \og d'uniformisation locale plong\'ee \fg{} de $f$, $f$ \'etant un \'el\'ement non d\'eg\'en\'er\'e par rapport \`a $\nu$. 
\begin{thm}\label{uniflocnondeg}
Consid\'erons les m\^emes hypoth\`eses que celle de la D\'efinition \ref{defnondeg}. Soit $f$ un \'el\'ement non d\'eg\'en\'er\'e par rapport \`a $u$. Il existe alors une suite locale encadr\'ee $(R,u)\rightarrow\left(R^{(l)},u^{(l)}\right)$ telle que $f$ soit un mon\^ome en $u^{(l)}$ multipli\'e par une unit\'e de $R^{(l)}$.
\\De plus, soit $I$ un id\'eal de $R$ tel que $\nu(f)=\min\limits_{g\in I}\lbrace\nu(g)\rbrace$. Notons $u=(w,v)$ et supposons que $I$ est engendr\'e uniquement par des mon\^omes en $w$. Alors, la suite locale encadr\'ee pr\'ec\'edente peut \^etre choisie ind\'ependante de $v$.
\end{thm}
\noindent\textit{Preuve}: La suite locale encadr\'ee par rapport \`a $\nu$ provient du Corollaire \ref{coroideal}. Ainsi, comme $f\in I$, il existe $z\in R^{(l)}$ tel que $f=zw^{\varepsilon_0}$ (selon les notations du Corollaire \ref{coroideal}). Comme $I$ est engendr\'e par $w^{\varepsilon_0}$ (Corollaire \ref{coroideal}) et par hypoth\`eses, on en conclut que:
\[\nu(z)=\nu(f)-\nu(w^{\varepsilon_0})=\nu(f)-\min\limits_{g\in I}\lbrace\nu(g)\rbrace=0.\]
Or $\nu$ est centr\'ee en $R^{(l)}$, donc, $z\in R^{(l)}{}^\times$.\\\qed
\subsection{Suite \'el\'ementaire uniformisante}\label{sectionsuiteelemunif}
~\smallskip ~\\ \indent Nous allons construire une uniformisation locale, par rapport \`a une valuation $\nu$, d'une hypersurface quasi-homog\`ene satisfaisant certaines conditions vis-\`a-vis de l'alg\`ebre gradu\'ee $G_\nu=(gr_\nu(R))^*$, o\`u, pour une alg\`ebre gradu\'ee n'ayant pas de diviseurs de z\'ero $G$, l'ag\`ebre gradu\'ee $G^*$ est d\'efinie par:
\[G^{*}=\left\lbrace \left.\dfrac{f}{g}\:\right|\:f,g\in G,\:g\textit{ homog\`ene},\:g\neq 0\right\rbrace.\] 
\\ \\\indent  Soient ($R,\mathfrak m,k)$ un anneau local r\'egulier, $u=(u_1,...,u_n)$ un syst\`eme r\'egulier de param\`etres de $R$ et $\nu$ une valuation centr\'ee en $R$ de groupe des valeurs $\Gamma$. Notons:
\[\beta_i=\nu(u_i),\:1\leqslant i\leqslant n.\]
Pour $r\in\lbrace 1,...,n-1\rbrace$ posons $t=n-r-1$.
\\Supposons que $r=\dim_\mathbb Q\left(\sum\limits_{i=1}^n\mathbb Q\beta_i\right)$, c'est-\`a-dire, quitte \`a renum\'eroter les variables, que $\beta_1,...,\beta_r$ sont $\mathbb Q$-lin\'eairement ind\'ependants dans $\Gamma\otimes_\mathbb Z\mathbb Q$ et qu'en particulier $\beta_n$ est $\mathbb Q$-combinaison lin\'eaire de $\beta_1,...,\beta_r$.
\\Notons $u=(w,v)$ avec:
\[v=(v_1,...,v_t)=(u_{r+1},...,u_{n-1}),\]
\[w=(w_1,...,w_r,w_n)=(u_1,...,u_r,u_n).\]
Soit $\Delta=\langle \beta_1,...,\beta_r\rangle$ le sous-groupe de $\Gamma$ engendr\'e par $\beta_1,...,\beta_r$. Notons:
\[\overline\alpha=\min\lbrace m\in\mathbb N^*\:\vert\:m\beta_n\in\Delta\rbrace.\]
Par hypoth\`eses, $\overline\alpha<+\infty$. Pour $i\in\lbrace 1,...,r,n\rbrace$, on note $x_i=in_\nu(u_i)$, on a donc $ord(x_i)=\beta_i$.
\\Par le Corollaire 4.6 de \cite{spiva}, on peut montrer que les $x_1,...,x_r$ sont alg\'ebriquement ind\'ependants sur $k$ dans $G_\nu$. Si $x_n$ est alg\'ebrique sur $k\left[x_1,...,x_r\right]$, notons $P$ le polyn\^ome minimal de $x_n$ sur $k\left[x_1,...,x_r\right]^*$, choisi unitaire et de plus bas degr\'e possible; sinon posons $P=0$. Si $P\neq 0$, notons $\alpha=d^{\:\circ}(P)$. Soient $\alpha_1,...,\alpha_r\in\mathbb Z$ tels que:
\[\overline\alpha\beta_n-\sum\limits_{i=1}^r\alpha_i\beta_i=0.\]
On peut montrer (Lemme 4.5 de \cite{spiva}) que $d=\frac{\alpha}{\overline\alpha}\in\mathbb N$. On note alors:
\[\overline y=x_1^{\alpha_1}...x_r^{\alpha_r},\]
\[y=w_1^{\alpha_1}...w_r^{\alpha_r},\]
\[\overline z=\dfrac{x_n^{\overline \alpha}}{\overline y},\]
\[z=\dfrac{w_n^{\overline \alpha}}{y}.\]
Si $P\neq 0$, alors $P$ est de la forme:
\[P(X)=\sum\limits_{i=0}^dc_i\overline y^{d-i}X^{i\overline\alpha},\]
o\`u $c_i\in k$, pour $0\leqslant i\leqslant d$, $c_d=1$ et $\sum\limits_{i=0}^dc_iX^i$ est le polyn\^ome minimal de $\overline z$ sur $k$ dans $G_\nu$.
\\Enfin, pour $0\leqslant i\leqslant d$, fixons un \'el\'ement $b_i\in R$ tel que $c_i\equiv b_i\mod \mathfrak m$. On pose alors:
\[Q=\sum\limits_{i=0}^db_iy^{d-i}w_n^{i\overline\alpha}.\]
\begin{prop}\label{propcaspartpolycle}
Avec les hypoth\`eses et les notations pr\'ec\'edentes, il existe une suite locale encadr\'ee par rapport \`a $\nu$ et ind\'ependante de $v$:
\[ \xymatrix{\left( R,u\right)=\left( R^{(0)},u^{(0)}\right) \ar[r]^-{\pi_{0}} & \left( R^{(1)},u^{(1)}\right) \ar[r]^-{\pi_{1}} & \ldots   \ar[r]^-{\pi_{l-1}} & \left( R^{(l)},u^{(l)}\right)}, \]
telle que, pour $0\leqslant i\leqslant l$, si on note:
\[u^{(i)}=\left(u_1^{(i)},...,u_{n_i}^{(i)}\right)\]
et $k^{(i)}$ le corps r\'esiduel de $R^{(i)}$, alors:
\begin{enumerate}
\item Les \'eclatements encadr\'es $\pi_0,...,\pi_{l-2}$ sont monomiaux. En particulier, $n_i=n$, $k^{(i)}=k$, pour $1\leqslant i<l$. De plus, $z\in R^{(l)}{}^\times$.
\item $n_l=\left \{ \begin{array}{ccc}  n & \textup{si} & P\neq 0 \\  n-1 & \textup{si} & P=0 \end{array} \right.$
\item Notons:
\[u^{(l)}=\left \{ \begin{array}{lcc}  \left(w_1^{(l)},...,w_r^{(l)},v,w_n^{(l)}\right) & \textup{si} & P\neq 0 \\  \left(w_1^{(l)},...,w_r^{(l)},v\right) & \textup{si} & P=0 \end{array} \right.\]
Alors, pour $j\in\lbrace 1,...,r,n\rbrace$, $w_j$ est un mon\^ome en $w_1^{(l)},...,w_r^{(l)}$ multipli\'e par une unit\'e de $R^{(l)}$.
\item Pour $j\in\lbrace 1,...,r\rbrace$, $w_j^{(l)}$ est un mon\^ome en $w$ dont les exposants peuvent \^etre n\'egatifs.
\item Si $P\neq 0$, alors:
\[w_n^{(l)}=\sum\limits_{i=0}^d b_iz^i=\dfrac{Q}{y^d}.\]
\item $k^{(l)}\simeq k\left(\overline z\right)\simeq\left \{ \begin{array}{lcc}  k(X) & \textup{si} & P= 0 \\  k[X]/\left(\sum\limits_{i=0}^d c_iX^i\right) & \textup{si} & P\neq 0 \end{array} \right.$
\end{enumerate}
\end{prop}
\noindent\textit{Preuve}: Nous ne donnerons qu'une id\'ee de preuve. Sans perte de g\'en\'eralit\'e, on peut supposer qu'il existe $\delta=(\delta_1,...,\delta_r,\delta_n)$ et $\gamma=(\gamma_1,...,\gamma_r,\gamma_n)$ tels que $z=\frac{w^\delta}{w^\gamma}$ et $\nu(w^\gamma)=\nu(w^\delta)$. En appliquant le Corollaire \ref{coroexistsuitlocale}, on obtient l'existence de la suite locale encadr\'ee par rapport \`a $\nu$ et ind\'ependante de $v$ telle que $w^\gamma$ divise $w^\delta$ dans $R^{(l)}$. En appliquant la Proposition \ref{quidivisequi}, on en d\'eduit que $z,z^{-1}\in R^{(l)}$.
\\ On montre (1) par r\'ecurrence. Plus pr\'ecis\'ement, la Proposition \ref{propmatriceeclat} implique des relations de d\'ependance entre les images des $\beta_i$ dans les $R^{(i')}$ et les images des $\delta$ et $\gamma$, ainsi que $z^{\pm\:1}=\frac{w_j^{(i')}}{w_1^{(i')}}$, o\`u $j\neq 1$ est tel que $\beta_j^{(i')}=\min\limits_{i\in\lbrace 1,...,r,n\rbrace} \left\lbrace\beta_{i}^{(i')}\right\rbrace=\beta_1^{(i')}$. La Proposition \ref{quidivisequi} permet de conclure que $z,z^{-1}\in R^{(i'+1)}$. En particulier si $i'=l-1$ on a (1).
\\En reprenant les notations de la Proposition \ref{propconstruceclatencad}, on remarque qu'on est dans le cas o\`u $h^\times\leqslant 1$ (plus pr\'ecis\'ement, le cas $h^{\times c}=h-1$, $t=h$ et le cas $h^{\times c}=t=h-1$). En utilisant la Proposition \ref{propmatriceeclat} et quitte \`a interchanger, on peut supposer que $\overline z=\frac{x_j^{(l-1)}}{x_1^{(l-1)}}$. Le cas $h^{\times c}=h-1$, $t=h$ se produit si et seulement si $\overline z=\frac{x_j^{(l-1)}}{x_1^{(l-1)}}$ est transcendant sur $k$ (et donc $P=0$). Le cas $h^{\times c}=t=h-1$ se produit quant \`a lui si et seulement si $\overline z$ est alg\'ebrique sur $k$ (c'est-\`a-dire $P\neq 0$). Ainsi (2) et (6) proviennent de l'\'etude directe de ces deux cas particuliers. Dans le cas $P\neq 0$, vu que $z^{\pm\:1}=\frac{w_j^{(l-1)}}{w_1^{(l-1)}}$, on peut prendre $w_n^{(l)}=\sum\limits_{i=0}^d b_iz^i$, ce qui prouve (5). Pour terminer, (3) et (4) sont une application directe de la Proposition \ref{propgenerem} avec $i=0$ et $i'=l$.\\\qed
\begin{prop}\label{proppolycleplus}
Reprenons les notations et les hypoth\`eses de la Proposition \ref{propcaspartpolycle}. Notons $\tilde Q=Q+h$, o\`u $h\in R$ est tel que $\nu_{0,u}(h)>\nu_{0,u}(Q)$. Alors, la Proposition \ref{propcaspartpolycle} est vraie avec $\tilde Q$ \`a la place de $Q$ dans (5).
\end{prop}
\noindent\textit{Preuve}: Par hypoth\`eses, on peut \'ecrire $h$ comme une somme finie $h=\sum\limits_\gamma h_\gamma u^\gamma$ o\`u $h_\gamma\in R$ et $\nu(u^\gamma)>\nu_{0,u}(Q)$. Soit $N=\max\lbrace\vert\gamma\vert\:\vert\: h_\gamma\neq 0\rbrace$. Apr\`es une suite locale encadr\'ee ind\'ependante de $u_{r+1},...,u_{n}$, on peut supposer que:
\[\nu(w_1)<\dfrac{1}{N}\left(\nu_{0,u}(h)-\nu_{0,u}(Q)\right).\]
Pour tout $i\in\lbrace r+1,...,n\rbrace$, effectuons $\left\lfloor\frac{\nu(u_i)}{\nu(w_1)}\right\rfloor$ \'eclatements le long de l'id\'eal $(u_i,w_1)$. On peut alors supposer que, pour $h_\gamma\neq 0$, $u^\gamma$ est divisible par un mon\^ome $\varpi_\gamma$ en $w_1,...,w_r$ tel que $\nu(\varpi_\gamma)\geqslant\nu_{0,u}(Q)$. En appliquant le Corollaire \ref{coroideal} \`a l'id\'eal monomial engendr\'e par $\lbrace y^d\rbrace\cup\lbrace\varpi_\gamma\:\vert\:h_\gamma\neq 0\rbrace$, on construit une suite locale encadr\'ee monomiale ind\'ependante de $u_{r+1},...,u_n$ telle que $y^d$ divise $h$.
\\Ainsi, avec cette hypoth\`ese, on peut consid\'erer la suite locale encadr\'ee construite dans la Proposition \ref{propcaspartpolycle}. Comme $y^d$ divise $Q$ dans $R^{(l)}$ et comme $y^d$ divise $h$, on en d\'eduit que $y^d$ divise $\tilde Q$ dans $R^{(l)}$. Le $w_n^{(l)}$ de la Proposition \ref{propcaspartpolycle} diff\`ere alors de $\frac{\tilde Q}{y^d}$ par des \'el\'ements appartenant \`a l'id\'eal $\left(u_1^{(l)},...,u_{n-1}^{(l)}\right)$.\\\qed
\begin{defi}\label{suiteelementaire}
Reprenons les notations et les hypoth\`eses de la Proposition \ref{proppolycleplus}. La suite locale encadr\'ee par rapport \`a $\nu$ et ind\'ependante de $v$ construite dans la Proposition \ref{proppolycleplus} sera appel\'ee \textbf{la suite \'el\'ementaire uniformisante associ\'ee \`a} $\boldsymbol{\left(R,u,\nu,n,\tilde Q\right)}$, ou plus simplement, \textbf{la} $\boldsymbol n$\textbf{-suite \'el\'ementaire uniformisante}, lorsque il n'y a pas d'ambigu\"it\'e possible dans les choix de $R,u,\nu$ et $\tilde Q$.
\end{defi}
\begin{rem}\label{remjsuiteelemeunif}
\textup{L'entier $n$ de la D\'efinition \ref{suiteelementaire} fait r\'ef\'erence au fait que la suite est d\'ependante uniquement des variables $u_1,...,u_r,u_n$. Pour $j\in\lbrace r+1,...,n\rbrace$, on peut d\'efinir une $j$-suite \'el\'ementaire en rempla\c{c}ant les variables $u_1,...,u_r,u_n$ par $u_1,...,u_r,u_j$.}
\end{rem}

\subsection{Suites formelles encadr\'ees}\label{suiteformencaddefchap3}
~\smallskip ~\\\indent Soit $(R,\mathfrak{m},k)$ un anneau local r\'egulier complet de dimension $n$ avec $\mathfrak{m}=\left(u_1,...,u_{n}\right)$. Soient $\nu$ une valuation de $K=Frac(R)$, centr\'ee en $R$, de groupe des valeurs $\Gamma$ et $\Gamma_1$ le plus petit sous-groupe isol\'e non nul de $\Gamma$. On note: 
\[H=\lbrace f\in R\:\vert\: \nu(f)\notin\Gamma_{1}\rbrace.\]
$H$ est un id\'eal premier de $R$ (voir Preuve du Th\'eor\`eme \ref{thmprelimcar0}). On suppose de plus que:
\[n=e(R,\nu)=emb.dim\left(R/H\right),\]
 c'est-\`a-dire que:
 \[H\subset\mathfrak{m}^2.\]
On note \'egalement $r=r(R,u,\nu)=\dim_\mathbb Q\left(\sum\limits_{i=1}^n\mathbb Q\nu(u_i)\right)$.
\\La valuation $\nu$ est unique si $ht(H)=1$, cas auquel on va se ramener gr\^ace au Corollaire \ref{hauteurcar0}. C'est la compos\'ee de la valuation $\mu:L^{*}\rightarrow\Gamma_{1}$ de rang $1$ centr\'ee en $R/H$, o\`u $L=Frac(R/H)$, avec la valuation $\theta :K^{*}\rightarrow \Gamma / \Gamma_{1}$, centr\'ee en $R_{H}$, telle que $k_{\theta}\simeq \kappa(H)$.
\\Par abus de notation, pour $f\in R$, on notera $\mu(f)$ au lieu de $\mu(f\mod H)$.
\begin{defi}
Soit $(R,u,k)\rightarrow (R',u',k')$ une suite locale encadr\'ee, on note $H'_{0}$ le transform\'e strict de $H$ dans $R'$. On dit que $\mu$ est \textbf{centr\'ee} en $R'$ si $\mu$ est centr\'ee en $R'/H'_{0}$. Dans ce cas, on dit que la suite locale encadr\'ee est une \textbf{suite locale encadr\'ee par rapport \`a} $\boldsymbol{\mu}$.
\end{defi}
\begin{defi}
Soit $(R,u)\rightarrow \left( R^{(1)},u^{(1)}\right) $ un \'eclatement local encadr\'e. Le morphisme induit par compl\'etion formelle est appel\'e un \textbf{\'eclatement formel encadr\'e par rapport \`a} $\boldsymbol{\mu}$.
\\Soient $\widehat{K^{(1)}}=Frac\left( \widehat{R^{(1)}}\right) $, $H^{(0)}$ le transform\'e strict de $H$ dans $R^{(1)}$ et $\overline{H}^{(1)}$ l'id\'eal premier implicite de $\widehat{R^{(1)}}/H^{(0)}\widehat{R^{(1)}}$. 
\\On appelle \textbf{transform\'e} de $H$ dans $\widehat{R^{(1)}}$, not\'e $H^{(1)}$, la pr\'eimage de $\overline{H}^{(1)}$ dans $\widehat{R^{(1)}}$. 
\\Enfin, on appelle \textbf{valuation induite par} $\boldsymbol{\mu}$ en $\widehat{R^{(1)}}$, not\'ee  $\mu^{(1)}$, l'unique extension de $\mu$ de $\kappa\left( H\right) $ \`a $\kappa\left( H^{(1)}\right)$, centr\'ee en $\widehat{R^{(1)}}/H^{(1)}$ et donn\'ee par le Th\'eor\`eme \ref{valetenuniqidealimpl}.
\end{defi}
\begin{defi}
Une suite de morphismes locaux:
\[ \xymatrix{\left( R,u\right) \ar[r]^-{\pi_{0}} & \left( R^{(1)},u^{(1)}\right) \ar[r]^-{\pi_{1}} & \ldots \ar[r]^-{\pi_{l-2}} & \left( R^{(l-1)},u^{(l-1)}\right)  \ar[r]^-{\pi_{l-1}} & \left( R^{(l)},u^{(l)}\right)} \]
est appel\'ee une \textbf{suite formelle encadr\'ee par rapport \`a} $\boldsymbol{\mu}$ si, la suite:
\[ \xymatrix{\left( R,u\right)   \ar[r]^-{\pi_{0}} & \left( R^{(1)},u^{(1)}\right)  \ar[r]^-{\pi_{1}} & \ldots \ar[r]^-{\pi_{l-2}} & \left( R^{(l-1)},u^{(l-1)}\right) } \]
est une suite formelle encadr\'ee par rapport \`a $\mu$ et si $\pi_{l-1}$ est un \'eclatement formel encadr\'e par rapport \`a la valuation $\mu^{(l-1)}$, induite par $\mu$ sur $R^{(l-1)}$.
\end{defi}
Pour tout \'eclatement local encadr\'e de la forme $(R,u)\rightarrow \left( R^{(1)},u^{(1)}\right) $, on d\'efinit une valuation $\nu^{(1)}$ centr\'ee en $\widehat{R^{(1)}}$ comme suit: fixons une valuation $\theta^{(1)}$ de $\widehat{K^{(1)}}$, centr\'ee en $\left( \widehat{R^{(1)}}\right) _{H^{(1)}}$ et telle que $k_{\theta^{(1)}}\simeq\kappa\left( H^{(1)}\right) $. On pose alors $\nu^{(1)}=\theta^{(1)}\circ\mu^{(1)}$.
\\Etant donn\'ee une suite formelle encadr\'ee:
\[ \xymatrix{\left( R,u\right) \ar[r]^-{\pi_{0}} & \left( R^{(1)},u^{(1)}\right) \ar[r]^-{\pi_{1}} & \ldots \ar[r]^-{\pi_{l-2}} & \left( R^{(l-1)},u^{(l-1)}\right)  \ar[r]^-{\pi_{l-1}} & \left( R^{(l)},u^{(l)}\right)}; \]
on peut, par r\'ecurrence sur $1\leqslant i\leqslant l-1$, construire une valuation $\mu^{(i)}$, centr\'ee en $R^{(i)}$ telle que le plus petit sous-groupe non nul du groupe des valeurs de $\nu^{(i)}$ soit $\Gamma_{1}$, et d\'efinir le transform\'e de $H$ dans $R^{(i)}$, not\'e $H^{(i)}$. Par construction, on a:
\[H^{(i)}=\lbrace f\in R^{(i)}\:\vert\: \nu^{(i)}(f)\notin\Gamma_{1}\rbrace.\]
On note alors:
\[e\left( R^{(i)},\nu^{(i)}\right)=emb.dim\left(R^{(i)}/H^{(i)}\right);\]
\[r\left( R^{(i)},u^{(i)},\nu^{(i)}\right)=\dim_\mathbb Q\left(\sum\limits_{j=1}^{n_i}\mathbb Q\nu^{(i)}\left(u_j^{(i)}\right)\right),\]
o\`u $u^{(i)}=(u_1^{(i)},...,u_{n_i}^{(i)})$.
\\\indent On note $\nu_{0,u}$ la valuation monomiale centr\'ee en $R$ et associ\'ee \`a $u=\left( u_{1},...,u_{n}\right)$ et \`a $\nu (u_{1}),...,\nu (u_{n})$. Par la Remarque \ref{valmonpluspetite}, pour tout $f\in R$, on a:
\[\nu_{0,u}(f)\leqslant\nu(f).\]
\begin{rem}\label{Hnulcar0}
\textup{Supposons que $n=r$. Pour une suite formelle encadr\'ee de la forme:
\[ \xymatrix{\left( R,u\right) \ar[r]^-{\pi_{0}} & \left( R^{(1)},u^{(1)}\right) \ar[r]^-{\pi_{1}} & \ldots \ar[r]^-{\pi_{l-2}} & \left( R^{(l-1)},u^{(l-1)}\right)  \ar[r]^-{\pi_{l-1}} & \left( R^{(l)},u^{(l)}\right)}, \]
on note $\nu_{0,u^{(l)}}$ la valuation monomiale centr\'ee en $R^{(l)}$ associ\'ee \`a $u^{(l)}$ et \`a  $\nu^{(l)}\left( u_{1}^{(l)}\right),...,\nu^{(l)}\left( u_{n}^{(l)}\right)$.
\\Si $f\in H\setminus\lbrace 0\rbrace$, alors:
\[\nu_{0,u^{(l)}}(f)<\nu(f),\]
pour toute suite formelle encadr\'ee de la forme pr\'ec\'edente et telle que: 
\[\nu^{(l)}\left( u_{1}^{(l)}\right),...,\nu^{(l)}\left( u_{n}^{(l)}\right)\in\Gamma_{1}.\]
Ainsi, comme $\nu_{0,u^{(l)}}(f)\in\Gamma_{1}$ et $\nu(f)\notin\Gamma_{1}$, on en d\'eduit que $H^{(i)}=(0)$, pour tout $i$.}
\end{rem}

\section{Un th\'eor\`eme de monomialisation en \'egale caract\'eristique}\label{monocar0}
\indent Soit $(R,\mathfrak{m},k)$ un anneau local r\'egulier complet \'equicaract\'eristique de dimension $n$ avec $\mathfrak{m}=\left(u_1,...,u_{n}\right)$. Soient $\nu$ une valuation de $K=Frac(R)$, centr\'ee en $R$, de groupe des valeurs $\Gamma$ et $\Gamma_1$ le plus petit sous-groupe isol\'e non nul de $\Gamma$. On note: 
\[H=\lbrace f\in R\:\vert\: \nu(f)\notin\Gamma_{1}\rbrace.\]
$H$ est un id\'eal premier de $R$ (voir Preuve du Th\'eor\`eme \ref{thmprelimcar0}). On suppose de plus que:
\[n=e(R,\nu)=emb.dim\left(R/H\right),\]
 c'est-\`a-dire que:
 \[H\subset\mathfrak{m}^2.\]
On note \'egalement $r=r(R,u,\nu)=\dim_\mathbb Q\left(\sum\limits_{i=1}^n\mathbb Q\nu(u_i)\right)$.
\\La valuation $\nu$ est unique si $ht(H)=1$, cas auquel on va se ramener gr\^ace au Corollaire \ref{hauteurcar0}. C'est la compos\'ee de la valuation $\mu:L^{*}\rightarrow\Gamma_{1}$ de rang $1$ centr\'ee en $R/H$, o\`u $L=Frac(R/H)$, avec la valuation $\theta :K^{*}\rightarrow \Gamma / \Gamma_{1}$, centr\'ee en $R_{H}$, telle que $k_{\theta}\simeq \kappa(H)$.
\\Par abus de notation, pour $f\in R$, on notera $\mu(f)$ au lieu de $\mu(f\mod H)$.
Par le th\'eor\`eme de Cohen, on peut supposer que $R$ s'\'ecrit sous la forme:
\[R=k\left[ \left[ u_{1},...,u_{n}\right] \right].\]
\subsection{Un premier th\'eor\`eme de monomialisation}
~\smallskip ~\\\indent Pour $j\in \lbrace r+1,...,n\rbrace$, on note $\lbrace Q_{j,i}\rbrace_{i\in\Lambda_{j}}$ l'ensemble des polyn\^omes-cl\'es de l'extension $k\left( \left( u_{1},...,u_{j-1}\right) \right)\hookrightarrow k\left( \left( u_{1},...,u_{j-1}\right) \right)(u_{j})$, $\textbf{Q}_{j,i}=\left\lbrace Q_{j,i'}\vert i'\in\Lambda_{j},i'<i\right\rbrace $, $\Gamma^{(j)}$ le groupe des valeurs de $\nu_{\vert k\left( \left( u_{1},...,u_{j}\right) \right)}$ et $\nu_{j,i}$ la $i$-troncature de $\nu$ pour cette extension.
\begin{thm}\label{thmeclatformcar0}
Reprenons les hypoth\`eses pr\'ec\'edentes et supposons de plus que pour $j\in \lbrace r+1,...,n\rbrace$, il existe un ensemble $\lbrace Q_{j,i}\rbrace_{i\in\Lambda_{j}}$ de polyn\^omes-cl\'es $1$-complet pour la valuation $\nu_{\vert k\left( \left( u_{1},...,u_{j}\right) \right)}$ ne poss\'edant pas de polyn\^ome-cl\'e limite. Deux cas se pr\'esentent:
\begin{enumerate}
\item Ou bien $H\neq(0)$ et il existe une suite formelle encadr\'ee:
\[ \xymatrix{\left( R,u\right) \ar[r]^-{\pi_{0}} & \left( R^{(1)},u^{(1)}\right) \ar[r]^-{\pi_{1}} & \ldots \ar[r]^-{\pi_{l-2}} & \left( R^{(l-1)},u^{(l-1)}\right)  \ar[r]^-{\pi_{l-1}} & \left( R^{(l)},u^{(l)}\right)} \]
telle que:
\[\left( e\left( R^{(l)},\nu^{(l)}\right) , e\left( R^{(l)},\nu^{(l)}\right) - r\left( R^{(l)},u^{(l)},\nu^{(l)}\right)\right) <_{lex}\left(e(R,\nu),n-r\right);\]
\item Ou bien $H=(0)$ et pour tout $f\in R$, il existe une suite formelle encadr\'ee:
\[ \xymatrix{\left( R,u\right) \ar[r]^-{\pi_{0}} & \left( R^{(1)},u^{(1)}\right) \ar[r]^-{\pi_{1}} & \ldots \ar[r]^-{\pi_{l-2}} & \left( R^{(l-1)},u^{(l-1)}\right)  \ar[r]^-{\pi_{l-1}} & \left( R^{(l)},u^{(l)}\right)} \]
telle que $f$ soit un mon\^ome en $u^{(l)}$ multipli\'e par une unit\'e de $R^{(l)}$.
\end{enumerate}
\end{thm}
\noindent\textit{Preuve}: On proc\`ede par r\'ecurrence sur $n-r$. Si $n=r$ alors $\nu(u_{1}),...,\nu (u_{n})$ sont $\mathbb{Q}$-lin\'eairement ind\'ependants et donc, tout $f\in R$ contient un unique mon\^ome de valuation minimale. En particulier,
\[\forall\:f\in R,\:\nu_{0,u}(f)=\nu(f).\]
Par la Remarque \ref{Hnulcar0}, $H=(0)$. Prenons alors un \'el\'ement $f\in R$, par le Th\'eor\`eme \ref{uniflocnondeg}, il existe une suite locale encadr\'ee:
\[ \xymatrix{\left( R,u\right) \ar[r]^-{\pi_{0}} & \left( R^{(1)},u^{(1)}\right) \ar[r]^-{\pi_{1}} & \ldots \ar[r]^-{\pi_{i-2}} & \left( R^{(i-1)},u^{(i-1)}\right)  \ar[r]^-{\pi_{i-1}} & \left( R^{(i)},u^{(i)}\right)} \]
telle que $f$ soit un mon\^ome en $u^{(i)}$ multipli\'e par une unit\'e de $R^{(i)}$. En passant au compl\'et\'e \`a chaque pas, on obtient la suite formelle encadr\'ee satisfaisant (2).
\\ \indent Supposons que $n-r>0$ et que l'on a d\'ej\`a construit une suite formelle encadr\'ee pour toutes les valeurs strictement plus petites et satisfaisant la conclusion du Th\'eor\`eme \ref{thmeclatformcar0}. On veut montrer le r\'esultat pour l'anneau $k\left[\left[u_1,...,u_{n-r}\right]\right]\left[\left[u_{n-r+1}\right]\right]$.
\\ Soit $\left(R^{(i)},\mathfrak m^{(i)},k^{(i)}\right)$ un anneau local apparaissant dans une suite formelle encadr\'ee. Par le Th\'eor\`eme de Cohen, on peut supposer que:
\[R^{(i)}=k^{(i)}\left[\left[u_1^{(i)},...,u_{n_i}^{(i)}\right]\right].\]
Dans un premier temps, montrons que l'on peut toujours se ramener aux hypoth\`eses suivantes:
\[H^{(i)}\cap R^{(i)}=(0)\:,n_i=n,\:r_i=r,\]
o\`u $r_i=r\left(R^{(i)},u^{(i)},\nu^{(i)}\right)$. En effet, si pour un certain $i$, on a:
\[H^{(i)}\cap R^{(i)}\neq(0),\]
en notant:
\[R_{n_i-1}^{(i)}=k^{(i)}\left[\left[u_1^{(i)},...,u_{n_i-1}^{(i)}\right]\right]\textup{ et }\overline u^{(i)}=\left(u_1^{(i)},...,u_{n_i-1}^{(i)}\right),\]
on peut appliquer l'hypoth\`ese de r\'ecurrence sur $n-r$ pour construire une suite formelle encadr\'ee:
\[ \xymatrix{\left(R_{n_i-1}^{(i)},\overline u^{(i)}\right) \ar[r] & \left( R_{n_i-1}^{(i,1)},\overline u^{(i,1)}\right) \ar[r] & \ldots \ar[r] & \left( R_{n_i-1}^{(i,l)},\overline u^{(i,l)}\right)} \]
telle que $e\left( R_{n_i-1}^{(i,l)},\nu^{(i,l)}\right)<e\left(R_{n_i-1}^{(i)},\nu^{(i)}\right)=n_i-1$. Notons alors:
\[\tilde R^{(j)}=R_{n_i-1}^{(i,j)}\left[\left[u_{n_i}^{(i)}\right]\right],\:1\leqslant j\leqslant l,\]
on obtient une suite formelle encadr\'ee:
\[ \xymatrix{\left( R^{(i)},u\right) \ar[r] & \left( \tilde R^{(1)},u^{(1)}\right) \ar[r] & \ldots \ar[r]& \left( \tilde R^{(l-1)},u^{(l-1)}\right)  \ar[r] & \left(\tilde R^{(l)},u^{(l)}\right)} \]
telle que:
\[e\left(\tilde R^{(l)},\nu^{(l)}\right)< e\left(R^{(i)},\nu^{(i)}\right)\leqslant e(R,\nu).\]
De m\^eme, s'il existe un $i$ tel que $n_i<n$ ou $r_i>r$, la suite formelle encadr\'ee recherch\'ee est d\'ej\`a construite et il n'y a rien \`a faire.
\\\indent Ainsi, on peut supposer que, pour tous les anneaux $R^{(i)}$ apparaissant dans n'importe quelle suite formelle encadr\'ee:
\[n_i=n,\:r_i=r,\:H^{(i)}\cap R^{(i)}=(0).\]
En particulier, $\nu_{\vert\:R_{n_i-1}^{(i)}}^{(i)}$ est de rang $1$ et $e\left(R_{n_i-1}^{(i)},\mu^{(i)}\right)=n-1=e\left(R_{n-1},\mu\right)$.
\\\indent\`A partir de maintenant, quitte \`a poser $N=n-r+1$, on peut supposer qu'on est dans le cas o\`u $r=1$.
\\ \\\indent Pour achever la preuve, nous allons montrer qu'\`a une suite formelle pr\`es, les id\'eaux $H^{(i)}$ sont principaux engendr\'es par un polyn\^ome unitaire en $u_n$ et qu'il suffit de monomialiser ce type d'\'el\'ements.
\subsection{L'id\'eal premier implicite est engendr\'e par un polyn\^ome unitaire}
~\smallskip ~\\\indent Pour finir la preuve du Th\'eor\`eme \ref{thmeclatformcar0}, on a doit d\'emontrer la Proposition \ref{princcar0}. On suppose toujours que l'hyoth\`ese de r\'ecurrence est v\'erifi\'ee \`a l'\'etape $n-1$, c'est-\`a-dire que le Th\'eor\`eme \ref{thmeclatformcar0} est vraie pour $R$ de dimension $n-1$.
\begin{prop}\label{princcar0}
Reprenons les hypoth\`eses pr\'ec\'edentes et supposons que le Th\'eor\`eme \ref{thmeclatformcar0} est vrai en dimension $n-1$.
Notons $R_{n-1}=k\left[ \left[ u_{1},...,u_{n-1}\right] \right]$ et supposons que $H\not\subset R_{n-1}$ et $H\cap R_{n-1}=(0)$.
\\Soit $f\in H\setminus\lbrace 0\rbrace$. \`A une suite formelle encadr\'ee pr\`es, $f$ s'\'ecrit sous la forme:
\[f=\alpha f_{n-1}P;\]
o\`u $\alpha\in R^{\times}$, $f_{n-1}\in R_{n-1}$ et $P$ est un polyn\^ome unitaire en $u_{n}$.
\end{prop}
\noindent\textit{Preuve}: La valuation $\mu$ de rang $1$ centr\'ee en $R/H$ induit une valuation de rang $1$ centr\'ee en $R_{n-1}$.
\\Soit $f\in H$, $f\neq 0$, on peut \'ecrire $f=\sum\limits_{j\geqslant 0}b_{j}u_{n}^{j}$, avec $b_{j}\in R_{n-1}$. Par hypoth\`ese, on peut supposer qu'il existe $j\geqslant 0$ tel que $b_{j}\notin H\cap R_{n-1}$. On pose alors:
\[\beta=\min_{j\geqslant 0} \lbrace \mu (b_{j})\:\vert\: b_{j}\notin H\cap R_{n-1}\rbrace.\]
Soit $d$ le plus petit entier naturel tel que $\mu(b_{d})=\beta$ (donc $b_{d}\notin H\cap R_{n-1}$).
\\Soit $N\geqslant d$ un entier naturel non nul tel que, pour tout $j>N$, \[b_{j}\in\left( b_{0},...,b_{N}\right).\]
Soit $j\in\lbrace 0,...,N\rbrace$, comme, par hypoth\`eses, $R_{n-1}$ est un anneau local, r\'egulier et complet de dimension $n-1$, on peut appliquer l'hypoth\`ese de r\'ecurrence (Th\'eor\`eme \ref{thmeclatformcar0} en dimension $n-1$) \`a cet anneau muni de la valuation $\mu$ et \`a l'\'el\'ement $b_{j}$. Il existe donc une suite locale encadr\'ee:
\[\pi:\left( R_{n-1},\left( u_{1},...,u_{n-1}\right) \right) \rightarrow ... \rightarrow \left( R',(u'_{1},...,u'_{n-1})\right)\]
telle que $b_{j}$ soit un mon\^ome en $(u'_{1},...,u'_{n-1})$ (multipli\'e par une unit\'e de $R'$). En passant \`a chaque pas de la suite au compl\'et\'e formel, on obtient une suite formelle encadr\'ee telle que $b_{j}$ est un mon\^ome en $(u'_{1},...,u'_{n-1})$ (multipli\'e par une unit\'e de $R'$). 
\\De plus, par le (1) de la Proposition \ref{propgenerem}, la propri\'et\'e d'\^etre un mon\^ome multipli\'e par une unit\'e est pr\'eserv\'ee pour tous les \'eclatements encadr\'es suivants. Ainsi, on peut choisir $\pi$ de telle sorte que les $b_{0},..., b_{N}$ soient simultan\'ement des mon\^omes en $(u'_{1},...,u'_{n-1})$.
\\Par les choix de $\beta$ et de $d$ et par le Corollaire \ref{coroideal}, apr\`es une suite locale encadr\'ee de plus, on peut se ramener \`a la situation o\`u $b_{d}$ divise $b_{j}$, $0\leqslant j \leqslant N$ et donc, $b_{d}$ divise $b_{j}$ pour tout $j\geqslant 0$.
\\Ainsi, $\frac{f}{b_{d}}\in R'\left[ \left[ u_{n}\right] \right] $ et satisfait les hypoth\`eses du th\'eor\`eme de pr\'eparation de Weierstrass.
\\ \qed
\begin{coro}\label{hauteurcar0}
Sous les m\^emes hypoth\`eses que la Proposition \ref{princcar0}, on a:
\[ht\left( H\right) \leqslant 1.\]
\end{coro}
\noindent\textit{Preuve}: Si $H=(0)$, il n'y a rien \`a montrer. Sinon, prenons $f\in H$ tel que $f\neq 0$. Comme la hauteur de $H$ est croissante lorsque l'on fait des suites locales ou formelles encadr\'ees, (Corollaire \ref{coroht}), par la Proposition \ref{princcar0}, on peut supposer que $f$ est un polyn\^ome unitaire en $u_{n}$ \`a coefficients dans $R_{n-1}$. Ainsi, l'extension d'anneaux:
\[\sigma: R_{n-1}\hookrightarrow  R_{n-1} \left[ \left[ u_{n}\right] \right] / (f) \]
est finie. La pr\'eimage de l'id\'eal $H/ \left( f\right) $ par $\sigma$ est $(0)$. Comme la hauteur est pr\'eserv\'ee par les extensions finies d'anneaux (\cite{matsualg}, Theorem 20), on a:
\[ht\left(  H/(f)\right) = ht((0))=0.\]
Ainsi, $ht\left( H\right)=1$.
\\ \qed
\begin{coro}\label{engcar0}
Sous les hypoth\`eses de la Proposition \ref{princcar0}, \`a une suite formelle encadr\'ee pr\`es, l'id\'eal $H$ est principal engendr\'e par un polyn\^ome unitaire en $u_{n}$.
\end{coro}
\noindent\textit{Preuve}: C'est une cons\'equence directe du Corollaire \ref{hauteurcar0}.
\\ \qed
\begin{rem}
\textup{En appliquant le Corollaire \ref{engcar0} \`a chaque anneau $R^{(i)}$ apparaissant dans la preuve du Th\'eor\`eme \ref{thmeclatformcar0}, on peut supposer, \`a une suite formelle encadr\'ee pr\`es, que pour tout $i$, l'id\'eal $H^{(i)}$ est principal, engendr\'e par un polyn\^ome unitaire en $u_n^{(i)}$.
\\Remarquons qu'il existe alors une unique valuation $\theta^{(i)}$ centr\'ee en $\left(R^{(i)}\right)_{H^{(i)}}$ (c'est la valuation triviale si $ht\left(H^{(i)}\right)=0$ et la valuation discr\`ete centr\'ee en $ \left(R^{(i)}\right)_{H^{(i)}}$ si $ht\left(H^{(i)}\right)=1$). Ainsi, l'extension $\nu^{(i)}$ de $\nu$ \`a $R^{(i)}$ est d\'etermin\'ee de mani\`ere unique.
\\\indent Pour achever la preuve du Th\'eor\`eme \ref{thmeclatformcar0}, il suffit d'obtenir le r\'esultat pour des polyn\^omes en $u_n$.}
\end{rem}
\subsection{Monomialisation des polyn\^omes}
\begin{prop}\label{eclatpolycar0}
Sous les hypoth\`eses pr\'ec\'edentes et en supposant que le Th\'eor\`eme \ref{thmeclatformcar0} est vrai en dimension $n-1$, pour tout polyn\^ome $f\in k\left[\left[ u_{1},...,u_{n-1}\right]\right]\left[u_n\right]$, il existe une suite formelle encadr\'ee $(R,u)\rightarrow (R',u')$ telle que $f$ soit un mon\^ome en $u'$ multipli\'e par une unit\'e de $R'$.
\\Supposons de plus que $f$ soit irr\'eductible dans $k\left[\left[ u_{1},...,u_{n-1}\right]\right]\left[\left[u_n\right]\right]$, la suite formelle encadr\'ee pr\'ec\'edente peut alors \^etre choisie de telle sorte que $u'_n$ divise $f$ et $u_n'^2$ ne divise pas $f$ dans $R'$.
\end{prop}
\noindent\textit{Fin de la preuve du Th\'eor\`eme \ref{thmeclatformcar0} en supposant la Proposition \ref{eclatpolycar0} vraie}: \\Si $H\neq (0)$, prenons $f\in H\cap k\left[\left[u_{1},...,u_{n-1}\right]\right]\left[u_n\right]$, $f\neq 0$; sinon prenons $f\in k\left[\left[u_{1},...,u_{n-1}\right]\right]\left[u_n\right]\setminus\lbrace 0\rbrace$. Par hypoth\`eses, il existe une suite formelle encadr\'ee $(R,u)\rightarrow (R',u')$ telle que $f$ soit un mon\^ome en $u'$ multipli\'e par une unit\'e de $R'$. Notons $H'$ le transform\'e de $H$ dans $R'$.
\\Si $H\neq (0)$, alors, par d\'efinition, $\nu(f)\notin\Gamma_1$ et donc, il existe un $j$ tel que $\nu(u'_j)\notin\Gamma_1$, c'est-\`a-dire, $u'_j\in H'$.
Ainsi, $e(R',\nu)\leqslant n-1<n=e(R,\nu)$ et on est dans la situation (1) du Th\'eor\`eme \ref{thmeclatformcar0}. Si $H=(0)$ et $f\in k\left[\left[u_{1},...,u_{n-1}\right]\right]\left[u_n\right]\setminus\lbrace 0\rbrace$, on se retrouve dans la situation (2) par hypoth\`eses.
\\ Enfin, si $H=(0)$ et $f\in R\setminus\lbrace 0\rbrace$, non n\'ecessairement un polyn\^ome en $u_n$, \'ecrivons $f=f'+f''$ avec $\nu_{0,u}(f'')>\nu(f)$ (et donc $\nu(f)=\nu(f')$). Par le cas polynomial vu avant, il existe une suite formelle encadr\'ee $(R,u)\rightarrow (R',u')$ telle que $f'$ soit un mon\^ome en $u'$ multipli\'e par une unit\'e de $R'$. Or $\nu_{0,u'}(f'')\geqslant\nu_{0,u}(f'')>\nu(f)=\nu(f')$. Par le Corollaire \ref{coroideal}, quitte \`a compl\'eter, il existe une suite formelle encadr\'ee $(R',u')\rightarrow (R'',u'')$ telle que $f$ soit un mon\^ome en $u''$ multipli\'e par une unit\'e de $R''$. \\ \qed
\\ \\\noindent\textit{Preuve de la Proposition \ref{eclatpolycar0}}: On va montrer le r\'esultat par r\'ecurrence sur le degr\'e de $f$. Si $d_{u_n}^{\:\circ}(f)=1$, la Proposition \ref{eclatpolycar0} est alors \'evidente.
\\Soit $f\in k\left[\left[u_{1},...,u_{n-1}\right]\right]\left[u_n\right]$ de degr\'e $d>1$. Par hypoth\`ese de r\'ecurren\-ce, on suppose que la Proposition \ref{eclatpolycar0} est vraie pour tout polyn\^ome de degr\'e strictement inf\'erieur \`a $d$.
\\Par 
hypoth\`ese, vu que l'ensemble de polyn\^omes-cl\'es est $1$-complet et ne poss\`ede pas de polyn\^ome-cl\'e limite,  il existe $i\in\mathbb N^*$ tel que $\nu(f)=\nu_{n,i}(f)$. Ceci veut dire qu'il existe un d\'eveloppement $(n,i)$-standard de $f$ de la forme:
\[f=\sum\limits_{j=0}^Nc_jQ_{n,i}^j,\]
o\`u les $c_j$ sont des d\'eveloppements $(n,i)$-standards n'impliquant pas $Q_{n,i}$ et $\nu(f)=\nu_{n,i}(f)$.
\\On rappelle que pour $i\in\mathbb N^*$, on note $\alpha_{n,i}=d_{Q_{n,i-1}}^{\:\circ}(Q_{n,i})$.
\\Supposons qu'il existe $l\in\mathbb N^*$ tel que $\alpha_{n,l}>1$, prenons alors ce $l$. S'il n'en existe pas, prenons un $l$ suffisamment grand tel que $f=\sum\limits_{j=0}^Nc_jQ_{n,l}^j$ et $\nu(f)=\nu_{n,l}(f)$. Dans tous les cas, par d\'efinition des polyn\^omes-cl\'es et par 
hypoth\`ese, $l<\omega$. 
\\\indent Pour achever la preuve de la Proposition \ref{eclatpolycar0}, il nous suffit donc d'obtenir le r\'esultat voulu sur les polyn\^omes-cl\'es comme nous allons le voir dans la sous-section \ref{sectmonopolyclecar0} et la Proposition \ref{eclatpolyclecar0}.
\subsection{Monomialisation des polyn\^omes-cl\'es}\label{sectmonopolyclecar0}
\begin{prop}\label{eclatpolyclecar0}
Sous les hypoth\`eses pr\'ec\'edentes et en supposant que le Th\'eor\`eme \ref{thmeclatformcar0} est vrai en dimension $n-1$, il existe une suite formelle encadr\'ee: \[(R,u)\rightarrow (R',u'),\] o\`u $u=(u_{1},...,u_{n})$, $u'=(u'_{1},...,u'_{n})$, et v\'erifiant les propri\'et\'es suivantes:
\begin{enumerate}
\item Pour tout $q\in\mathbb{N}^*$ tel que $1\leqslant q \leqslant l$, le polyn\^ome-cl\'e $Q_{n,q}$ est un mon\^ome en $u'$ multipli\'e par une unit\'e de $R'$;
\item Dans $R'$, $u'_{n}$ divise $Q_{n,l}$ mais $u'^{2}_{n}$ ne divise pas $Q_{n,l}$.
\end{enumerate}
\end{prop}
\noindent\textit{Preuve de la Proposition \ref{eclatpolycar0} en supposant la Proposition \ref{eclatpolyclecar0} vraie}:
\\Par hypoth\`ese de r\'ecurrence sur $n-r$, n'importe quelle collection d'\'el\'ements de $k\left(\left( u_{1},...,u_{n-1}\right)\right)$ peut \^etre transform\'ee simultan\'ement en mon\^omes via une suite formelle encadr\'ee. De plus, en appliquant $n-r-1$ fois la Proposition \ref{propcaspartpolycle}, on peut supposer que seuls les $u'_1,...,u'_r$ apparaissent dans ces mon\^omes.
\\Si $\alpha_{n,l}=1$, on applique la Proposition \ref{propcaspartpolycle} \`a chaque polyn\^ome-cl\'e $Q_{n,1},...,Q_{n,l}$ et la Proposition \ref{eclatpolycar0} est d\'emontr\'ee.
\\Supposons que $\alpha_{n,l}>1$. Notons $f=\sum\limits_{j=0}^d a_ju_n^j$, $a_j\in k\left[\left[ u_{1},...,u_{n-1}\right]\right]$. Soit $j_0$ le plus grand $j\in\lbrace 0,...,d\rbrace$ tel que $\nu(a_{j_0})=\min\limits_{0\leqslant j\leqslant d}\lbrace\nu(a_j)\rbrace$. Par le Corollaire \ref{coroideal}, apr\`es une suite locale encadr\'ee ind\'ependante de $u_n$, et quitte \`a compl\'eter, on peut supposer que $a_{j_0}$ divise $a_j$, pour tout $j\in\lbrace 0,...,d\rbrace$. En appliquant le th\'eor\`eme de pr\'eparation de Weierstrass, on peut supposer que $f$ est un polyn\^ome unitaire en $u_n$ de degr\'e $d$.
\\Soit $f=\sum\limits_{j=0}^Nc_jQ_{n,l}^j$, $N=\left\lfloor\frac{d}{\alpha_{n,l}}\right\rfloor$, le d\'eveloppement $(n,l)$-standard de $f$. Par la Proposition \ref{eclatpolyclecar0}, il existe une suite formelle encadr\'ee telle que le d\'eveloppement $(n,l)$-standard de $f$ dans $R'$ soit de la forme $\sum\limits_{j=0}^Nc'_ju'^j_n$, $c'_j\in k'\left[\left[ u'_{1},...,u'_{n-1}\right]\right]$, multipli\'e par une unit\'e de $R'$.
\\Notons $j'_0$ le plus grand $j\in\lbrace 0,...,N\rbrace$ tel que $\nu(c'_{j'_0})=\min\limits_{0\leqslant j\leqslant N}\lbrace\nu(c'_j)\rbrace$. Toujours par le Corollaire \ref{coroideal}, apr\`es une suite locale encadr\'ee ind\'ependante de $u'_n$, et quitte \`a compl\'eter, on peut supposer que $c'_{j'_0}$ divise $c'_j$, pour tout $j\in\lbrace 0,...,N\rbrace$. En appliquant le th\'eor\`eme de pr\'eparation de Weierstrass, on peut supposer que $f$ est un polyn\^ome unitaire en $u'_n$ de degr\'e inf\'erieur ou \'egal \`a $N<d$. Pour conclure, il nous suffit juste d'appliquer l'hypoth\`ese de r\'ecurrence.\\\qed
\\ \\
\noindent\textit{Preuve de la Proposition \ref{eclatpolyclecar0}}: Comme $l\in\mathbb{N}^*$, le d\'eveloppement standard de $Q_{n,l}$ est:
\[Q_{n,l}=Q_{n,l-1}^{\alpha_{n,l}}+\sum\limits_{j=0}^{\alpha_{n,l}-1}\left(\sum\limits_{\overline{\gamma}_{n,l-1}}c_{n,l,j,\overline{\gamma}_{n,l-1}}\textbf{Q}_{n,l-1}^{\overline{\gamma}_{n,l-1}}\right)Q_{n,l-1}^j.\]
Par hypoth\`ese de r\'ecurrence, pour des valeurs strictement inf\'erieures \`a $n-r$, il existe une suite formelle encadr\'ee $(R,u)\rightarrow (R',u')$, ind\'ependante de $u_{n}$ telle que chaque \'el\'ement $c_{n,l,j,\overline{\gamma}_{n,l-1}}$ soit un mon\^ome en $u'_{1},...,u'_{n-1}$ multipli\'e par une unit\'e de $R'$. 
\\Pour chaque $j\in\lbrace r+1,...,n-1\rbrace$, appliquons la $j$-suite \'el\'ementaire uniformisante de la Remarque \ref{remjsuiteelemeunif}, suivie \`a chaque fois d'une compl\'etion formelle. On arrive alors \`a la situation o\`u les $\sum\limits_{\overline{\gamma}_{n,l-1}}c_{n,l,j,\overline{\gamma}_{n,l-1}}\textbf{Q}_{n,l-1}^{\overline{\gamma}_{n,l-1}}$ sont des mon\^omes en $u'_1,...,u'_{r}$ multipli\'es par une unit\'e de $R'$.
\\Appliquons $l-1$ fois la Proposition \ref{propcaspartpolycle}, on peut supposer de plus que:
\[Q_{n,l-1}=\eta u'_n,\]
o\`u $\eta$ est un mon\^ome en $u'_1,...,u'_{n-1}$ multipli\'e par une unit\'e de $R'$.
\\ En appliquant la Proposition \ref{propcaspartpolycle} \`a $u'_1,...,u'_r,u'_n$, quitte \`a passer au compl\'et\'e, on obtient une suite formelle encadr\'ee $(R',u')\rightarrow (R'',u'')$ telle que $Q_{n,l}$ soit un mon\^ome en $u''_1,...,u''_{r},u''_n$. On en d\'eduit imm\'ediatement (1) et (2) par construction, ceci ach\`eve la preuve de la Proposition \ref{eclatpolyclecar0} et donc celle du Th\'eor\`eme \ref{thmeclatformcar0}.\\\qed

\section{Un th\'eor\`eme de monomialisation en caract\'eristique mixte}\label{carmixtemono}
\indent Soient $(R,\mathfrak{m},k)$ un anneau local r\'egulier complet de caract\'eristique mixte de dimension $n$ avec $\mathfrak{m}=(x)=\left(x_1,...,x_{n}\right)$ et $\nu$ une valuation de $K=Frac(R)$ centr\'ee en $R$, de groupe des valeurs $\Gamma$. Soit $\Gamma_1$ le plus petit sous-groupe isol\'e non nul de $\Gamma$. On note: 
\[H=\lbrace f\in R\:\vert\: \nu(f)\notin\Gamma_{1}\rbrace.\]
$H$ est un id\'eal premier de $R$ (voir Preuve du Th\'eor\`eme \ref{thmprelimcar0}). On suppose de plus que:
\[n=e(R,\nu)=emb.dim\left(R/H\right),\]
 c'est-\`a-dire que:
 \[H\subset\mathfrak{m}^2.\]
On note \'egalement $r=r(R,x,\nu)=\dim_\mathbb Q\left(\sum\limits_{i=1}^n\mathbb Q\nu(x_i)\right)$ et $p=car(k)$.
\\La valuation $\nu$ consid\'er\'ee est la compos\'ee de la valuation $\mu:L^{*}\rightarrow\Gamma_{1}$ de rang $1$ centr\'ee en $R/H$, o\`u $L=Frac(R/H)$, avec la valuation $\theta :K^{*}\rightarrow \Gamma / \Gamma_{1}$, centr\'ee en $R_{H}$, telle que $k_{\theta}\simeq \kappa(H)$.
\\Par abus de notation, pour $f\in R$, on notera $\mu(f)$ au lieu de $\mu(f \mod H)$.

\begin{rem}\label{ppasdansH}
\textup{Si $p\in H$, alors $R/H$ est \'equicaract\'eristique et on est sous les hypoth\`eses de la Section \ref{monocar0}. Dans la suite on supposera donc que $p\notin H$.}
\end{rem}
\subsection{Suites formelles encadr\'ees et anneaux de caract\'eristique mixte}
\begin{lem}\label{pasm2}
Il existe $g\in W\left[\left[u_1,...,u_n\right]\right]$ \`a coefficients dans $W^{\times}$ tel que:
\[R\simeq W\left[\left[u_1,...,u_n\right]\right]/(p-g),\]
o\`u $W$ est un anneau local r\'egulier complet de dimension $1$ dont l'id\'eal maximal est engendr\'e par $p$.
\end{lem}
\noindent\textit{Preuve}: On sait qu'il existe un morphisme surjectif: \[\varphi:W\left[\left[u_1,...,u_n\right]\right]\twoheadrightarrow R,\] tel que $\varphi(u_i)=x_i$ et $\varphi_{\vert\:W}=id_W$. Comme $R$ est int\`egre (voir \cite{EGA4-1}, Corollaire 17.1.3), $\ker \varphi$ est un id\'eal premier et, en comparant les dimensions, on en d\'eduit que $ht(\ker \varphi)\leqslant 1$. Or, $ W\left[\left[u_1,...,u_n\right]\right]$ est factoriel donc, $\ker \varphi$ est un id\'eal principal engendr\'e par $f$. 
Comme $p\in\mathfrak m$, il existe $a_1,...,a_n\in R$ tels que:
\[p=a_1x_1+...+a_nx_n.\]
Or, $\varphi$ est surjective, il existe donc $b_1,...,b_n\in W\left[\left[u_1,...,u_n\right]\right]$ tes que $\varphi(b_i)=a_i$, $1\leqslant i\leqslant n$. Notons $g=b_1u_1+...+b_nu_n$, alors, $p-g\in \ker\varphi$.
\\Si un des $b_i$ est divisible par $p$, en notant $b_i=pb'_i$, $b'_i\in W\left[\left[u_1,...,u_n\right]\right]$, on remplace $b_i$ par $b'_ig$. En it\'erant ce processus, apr\`es un nombre au plus d\'enombrable de pas, on peut supposer que tous les $b_i$ sont non divisibles par $p$ et donc $b_i\in W^\times$, $1\leqslant i\leqslant n$.
\\Vu que $p$ est un param\`etre r\'egulier de $W\left[\left[u_1,...,u_n\right]\right]$, l'id\'eal $(p-g)$ est premier, de hauteur $1$ et inclu dans $\ker\varphi$, d'o\`u: $$\ker\varphi=(p-g).$$
avec $g\in (u_{1},...,u_{n})$ \`a coefficients non divisibles par $p$, donc dans $W^{\times}$.
\\\qed
\begin{rem}
\textup{Si $R$ est ramifi\'e ($p\in\mathfrak m^2$) alors $g\in\left(u_1,...,u_n\right)^2$.}
\end{rem}
\`A partir de maintenant on suppose que:
\[R= W\left[\left[u_1,...,u_n\right]\right]/(p-g),\]
avec $g\in\left(u_1,...,u_n\right)$ \`a coefficients dans $W^{\times}$. L'id\'eal maximal de $R$ est $\mathfrak{m}=\left(u_1,...,u_{n}\right)$.

~\\\indent Pour $j\in\lbrace 1,...,n\rbrace$, notons $K_j$ le corps des fractions de $W\left[\left[u_1,...,u_j\right]\right]$. Pour $j\in \lbrace r+1,...,n\rbrace$, on note $\lbrace Q_{j,i}\rbrace_{i\in\Lambda_{j}}$ l'ensemble des polyn\^omes-cl\'es de l'extension $K_{j-1}\hookrightarrow K_{j-1}(u_{j})$. Si $\gamma=(\gamma_1,...,\gamma_l)$, on note:
\[\boldsymbol Q_{j,l+1}^\gamma=\prod\limits_{i=1}^lQ_{j,i}^{\gamma_i}.\]
Pour $j\in\lbrace r,...,n\rbrace$ et $\boldsymbol{\mathrm i}=(i_{r+1},...,i_j)\in\Lambda_{r+1}\times...\times\Lambda_j$, on d\'efinit par r\'ecurrence une valuation $\nu_{\boldsymbol{\mathrm i}}$ de $K_j$ comme suit:
\\Si $j=r$, on pose $\nu_\emptyset=\nu_{\vert\:K_r}$. Supposons que la valuation $\nu_{(i_{r+1},...,i_{j-1})}$ de $K_{j-1}$ soit d\'ej\`a construite. Si $f\in K_{j-1}\left[u_j\right]$, $\nu_{\boldsymbol{\mathrm i}}(f)$ est d\'efini comme l'extension de $\nu_{(i_{r+1},...,i_{j-1})}(f)$ d\'etermin\'ee par $\boldsymbol Q_{j,i_j+1}$. Si $f\in K_r\left[\left[u_{r+1},...,u_j\right]\right]$, posons $N$ suffisamment grand de telle sorte que $f=f_1+f_2$ avec:
\begin{enumerate}
\item[(1)] $f_1\in K_r\left[\left[u_{r+1},...,u_{j-1}\right]\right]\left[u_j\right]$;
\item[(2)] $f_2\in\left(u_j^N\right)K_r\left[\left[u_{r+1},...,u_j\right]\right]$;
\item[(3)] $\nu_{0,u}(f_2)>\nu_{\boldsymbol{\mathrm i}}(f_1)$.
\end{enumerate}
On pose alors $\nu_{\boldsymbol{\mathrm i}}(f)=\nu_{\boldsymbol{\mathrm i}}(f_1)$.

\begin{lem}\label{remregulier} Supposons que, pour $j\in\lbrace r+1,...,n\rbrace$, l'ensemble $\lbrace Q_{j,i}\rbrace_{i\in\Lambda_{j}}$ soit un ensemble $1$-complet de polyn\^omes-cl\'es pour la valuation $\nu_{\vert K_j}$ ne poss\'edant pas de polyn\^ome-cl\'e limite. Si $\nu(p)\notin p\Gamma$, alors, \`a une suite formelle encadr\'ee pr\`es, on peut supposer $R$ de la forme:
\[R=R[r]\left[ \left[ u_{r+1},...,u_{n}\right] \right] ,\]
o\`u $R[r]$ est un anneau local r\'egulier complet (\'eventuellement ramifi\'e) de dimension $r$ et tel que $\nu_{\vert R[r]}$ soit monomiale par rapport au syst\`eme r\'egulier de param\`etres de $R[r]$ et de rang rationnel maximal.
\end{lem}
\noindent\textit{Preuve}: 
Consid\'erons l'\'el\'ement $g\in W\left[\left[u_1,...,u_n\right]\right]$ du Lemme \ref{pasm2}.
Par le Th\'eor\`eme \ref{existpolycle}, pour tout $j\in\lbrace r+1,...,n\rbrace$, la collection $\lbrace Q_{j,i}\rbrace_{i\in\Lambda_{j}}$ forme un ensemble complet de polyn\^omes-cl\'es, il existe donc $\boldsymbol{\mathrm i}=(i_{r+1},...,i_n)\in\Lambda_{r+1}\times...\times\Lambda_n$ tel que:
\[\nu_{\boldsymbol{\mathrm i}}(g)=\nu(g)\textup{ et }\nu\left(Q_{j,i}\right)\leqslant p,\]
pour tout $i\leqslant i_j$, $j\in\lbrace r+1,...,n\rbrace$ (on rappelle que, vu la Remarque \ref{ppasdansH} et comme $p=g$ dans $R$, $\nu(g)\in\Gamma_1$).
\\Notons $\overline g$ l'image de $g$ modulo $p$ dans $k\left[\left[u_1,...,u_n\right]\right]$ et $\overline\nu_{\boldsymbol{\mathrm i}}$ la valuation d\'efinie sur $k\left(\left(u_1,...,u_r\right)\right)\left[\left[u_{r+1},...,u_n\right]\right]$ comme la valuation $\nu_{\boldsymbol{\mathrm i}}$ mais en regardant les \'el\'ements modulo $p$. En appliquant le Th\'eor\`eme \ref{thmeclatformcar0} \`a la valuation $\overline\nu_{\boldsymbol{\mathrm i}}$, 
 il existe une suite formelle encadr\'ee $k\left[\left[u_1,...,u_n\right]\right]\rightarrow k'\left[\left[u'_1,...,u'_n\right]\right] $ telle que $\overline{g}$ soit un mon\^ome en $u'$ multipli\'e par une unit\'e de $ k'\left[\left[u'_1,...,u'_n\right]\right]$. On a alors:
\[\nu(g)=\nu_{\boldsymbol{\mathrm i}}(g)=\overline\nu_{\boldsymbol{\mathrm i}}\left(\overline{g}\right)=\nu_{0,u'}\left(\overline{g}\right).\]
En appliquant \`a chaque \'etape de l'algorithme du Th\'eor\`eme \ref{thmeclatformcar0} 
les m\^emes changements de variables \`a $W\left[\left[u_1,...,u_n\right]\right]$, on obtient une suite $W\left[\left[u_1,...,u_n\right]\right]\rightarrow (R^{(2)},u^{(2)})$ telle que:
\[g=u_{1}^{(2) \alpha_{1}}...u_{r}^{(2) \alpha_{r}}z+ph,\]
o\`u $\alpha_{1},...,\alpha_{r}\in\mathbb{Z}$, $z\in R^{(2) \times}$, $h\in R^{(2)}$. Or l'algorithme du Th\'eor\`eme \ref{thmeclatformcar0} consiste en une r\'ep\'etition de $n$-suites \'el\'ementaires uniformisantes (D\'efinition \ref{suiteelementaire}), ainsi, par la Proposition \ref{propconstruceclatencad} et par choix de $\boldsymbol{\mathrm i}$: \[h\in \left(u_{1}^{(2)},...,u_{n}^{(2)}\right).\]
On en d\'eduit donc que:
\[h\not\in R^{(2) \times}.\]
On peut alors \'ecrire:
\[p-g=p(1-h)-u_{1}^{(2) \alpha_{1}}...u_{r}^{(2) \alpha_{r}}z=w(p-u_{1}^{(2) \alpha_{1}}...u_{r}^{(2) \alpha_{r}}z'),\]
o\`u $w=1-h\in R^{(2) \times}$ et $z'=zw^{-1}\in R^{(2) \times}$.
\\\`A une suite formelle encadr\'ee pr\`es, on peut donc supposer que, dans $R$, on a:
\[p=u_{1}^{\alpha_{1}}...u_{r}^{\alpha_{r}}z\]
avec $\alpha_{1},...,\alpha_{r}\in\mathbb{Z}$, $z\in R^{\times}$.
\\Par hypoth\`eses, comme $\nu(p)\notin p\Gamma$ et $R$ est complet, il existe $\alpha_{i}\notin p\mathbb Z$ et donc:
\[z^{1/\alpha_{i}}\in R.\]
Quitte \`a faire le changement de variable:
\[v_{i}=u_{i}z^{1/\alpha_{i}},\]
on peut supposer que:
\[p=u_{1}^{\alpha_{1}}...u_{r}^{\alpha_{r}}\in \mathfrak{m}^{2}.\]
Ainsi, \`a une suite formelle encadr\'ee pr\`es, $R$ s'\'ecrit sous la forme:
\[R=W\left[\left[u_1,...,u_n\right]\right]/\left(p-u_{1}^{\alpha_{1}}...u_{r}^{\alpha_{r}}\right)\simeq R[r]\left[ \left[ u_{r+1},...,u_{n}\right] \right] ,\]
o\`u $R[r]=W\left[\left[u_1,...,u_r\right]\right]/\left(p-u_{1}^{\alpha_{1}}...u_{r}^{\alpha_{r}}\right)$ est un anneau local r\'egulier complet (\'eventuellement ramifi\'e) de dimension $r$ tel que $\nu_{\vert R[r]}=\nu_{0,(u_1,...,u_r)}$ et $rg.rat\left(\nu_{\vert R[r]}\right)=r$.\\\qed
\subsection{Un premier th\'eor\`eme de monomialisation}
~\smallskip ~\\\indent \`A partir de maintenant et ce jusqu'\`a la fin de la Section \ref{carmixtemono}, on suppose que:
\[\nu(p)\notin p\Gamma,\]
\[R=R[r]\left[ \left[ u_{r+1},...,u_{n}\right] \right] ,\]
o\`u $R[r]$ est un anneau local r\'egulier complet (\'eventuellement ramifi\'e) de dimension $r$ et tel que $\nu_{\vert R[r]}$ soit monomiale par rapport au syst\`eme r\'egulier de param\`etres de $R[r]$ et de rang rationnel maximal.
\\\indent Pour $j\in \lbrace r+1,...,n\rbrace$, on note $\lbrace Q_{j,i}\rbrace_{i\in\Lambda_{j}}$ l'ensemble des polyn\^omes-cl\'es de l'extension $Frac\left(R[r]\left[ \left[ u_{r+1},...,u_{j-1}\right] \right]\right)\hookrightarrow Frac\left(R[r]\left[ \left[ u_{r+1},...,u_{j-1}\right] \right]\right)(u_{j})$.

\begin{thm}\label{thmeclatform}
Reprenons les hypoth\`eses pr\'ec\'edentes et supposons que, pour $j\in \lbrace r+1,...,n\rbrace$, l'ensemble $\lbrace Q_{j,i}\rbrace_{i\in\Lambda_{j}}$ est un ensemble $1$-complet de polyn\^omes-cl\'es ne poss\'edant pas de polyn\^ome-cl\'e limite. Deux cas se pr\'esentent:
\begin{enumerate}
\item[(1)] Ou bien $H\neq (0)$ et il existe une suite formelle encadr\'ee:
\[ \xymatrix{\left( R,u\right) \ar[r]^-{\pi_{0}} & \left( R^{(1)},u^{(1)}\right) \ar[r]^-{\pi_{1}} & \ldots \ar[r]^-{\pi_{l-2}} & \left( R^{(l-1)},u^{(l-1)}\right)  \ar[r]^-{\pi_{l-1}} & \left( R^{(l)},u^{(l)}\right)} \]
telle que:
\[\left( e\left( R^{(l)},\nu^{(l)}\right) , e\left( R^{(l)},\nu^{(l)}\right) - r\left( R^{(l)},u^{(l)},\nu^{(l)}\right)\right) <_{lex}\left(e(R,\nu),n-r\right);\]
\item[(2)] Ou bien $H=(0)$ et pour tout $f\in R$, il existe une suite formelle encadr\'ee:
\[ \xymatrix{\left( R,u\right) \ar[r]^-{\pi_{0}} & \left( R^{(1)},u^{(1)}\right) \ar[r]^-{\pi_{1}} & \ldots \ar[r]^-{\pi_{l-2}} & \left( R^{(l-1)},u^{(l-1)}\right)  \ar[r]^-{\pi_{l-1}} & \left( R^{(l)},u^{(l)}\right)} \]
telle que $f$ soit un mon\^ome en $u^{(l)}$ fois une unit\'e de $R^{(l)}$.
\end{enumerate}
\end{thm}
\noindent\textit{Preuve}: On proc\`ede par r\'ecurrence sur $n-r$. Si $n=r$ alors $R=R[r]$. En particulier,
\[\forall\:f\in R,\:\nu_{0,u}(f)=\nu(f).\]
Par la Remarque \ref{Hnulcar0}, $H=(0)$. Prenons alors un \'el\'ement $f\in R$, par le Th\'eor\`eme \ref{uniflocnondeg}, il existe une suite locale encadr\'ee:
\[ \xymatrix{\left( R,u\right) \ar[r]^-{\pi_{0}} & \left( R^{(1)},u^{(1)}\right) \ar[r]^-{\pi_{1}} & \ldots \ar[r]^-{\pi_{i-2}} & \left( R^{(i-1)},u^{(i-1)}\right)  \ar[r]^-{\pi_{i-1}} & \left( R^{(i)},u^{(i)}\right)} \]
telle que $f$ soit un mon\^ome en $u^{(i)}$ multipli\'e par une unit\'e de $R^{(i)}$. En passant au compl\'et\'e \`a chaque pas, on obtient la suite formelle encadr\'ee satisfaisant (2).
\\ \indent Supposons que $n-r>0$ et que l'on ait d\'ej\`a construit une suite formelle encadr\'ee pour toutes les valeurs strictement plus petites et satisfaisant la conclusion du Th\'eor\`eme \ref{thmeclatform}. On va proc\'eder comme dans la preuve du Th\'eor\`eme \ref{thmeclatformcar0}.

\begin{prop}\label{princ}
Par le Lemme \ref{remregulier}, on peut supposer que $R$ s'\'ecrit sous la forme:
\[R=R[r]\left[ \left[ u_{r+1},...,u_{n}\right] \right] ,\]
o\`u $R[r]$ est un anneau local r\'egulier complet (\'eventuellement ramifi\'e). Notons $R_{n-1}=R[r]\left[ \left[ u_{r+1},...,u_{n-1}\right] \right]$ et supposons que $H\not\subset R_{n-1}$ et $H\cap R_{n-1}=(0)$. Supposons enfin que le Th\'eor\`eme \ref{thmeclatform} est vrai en dimension $n-1$.
\\Soit $f\in H\setminus\lbrace 0\rbrace$. \`A une suite formelle encadr\'ee pr\`es, $f$ s'\'ecrit sous la forme:
\[f=\alpha f_{n-1}P;\]
o\`u $\alpha\in R^{\times}$, $f_{n-1}\in R_{n-1}$ et $P$ est un polyn\^ome unitaire en $u_{n}$.
\end{prop}
\noindent\textit{Preuve}: La preuve est la m\^eme que celle de la Proposition \ref{princcar0}.
\\ \qed
\begin{coro}\label{hauteur}
Sous les m\^emes hypoth\`eses que la Proposition \ref{princ}, on a:
\[ht\left( H\right) \leqslant 1.\]
\end{coro}
\noindent\textit{Preuve}: La preuve est la m\^eme que celle de la Proposition \ref{hauteurcar0}.
\\ \qed
\begin{coro}\label{eng}
Sous les hypoth\`eses de la Proposition \ref{princ}, \`a une suite formelle encadr\'ee pr\`es, l'id\'eal $H$ est principal engendr\'e par un polyn\^ome unitaire en $u_{n}$.
\end{coro}
\noindent\textit{Preuve}: C'est une cons\'equence directe du Corollaire \ref{hauteur}.
\\ \qed
\\ \indent Soit $R^{(i)}$ un anneau local apparaissant dans une suite formelle encadr\'ee. Par le Lemme \ref{remregulier}, on peut \'ecrire $R^{(i)}$ sous la forme $B\left[ \left[ u_{n_{i}}\right] \right] $ o\`u $B$ est un anneau r\'egulier (\'eventuellement ramifi\'e) et si $H^{(i)}\cap R^{(i)}\neq (0)$, alors $H^{(i)}\subset \mathfrak{m}^{(i) 2}$. Par le Corollaire \ref{eng}, $H^{(i)}$ est engendr\'e par un polyn\^ome unitaire en $u_{n_{i}}$.
\begin{prop}\label{eclatpoly}
Sous les hypoth\`eses du Th\'eor\`eme \ref{thmeclatform}, pour tout polyn\^ome $f\in R[r]\left[\left[ u_{r+1},...,u_{n-1}\right]\right]\left[u_n\right]$, il existe une suite formelle encadr\'ee $(R,u)\rightarrow (R',u')$ telle que $f$ soit un mon\^ome en $u'$ multipli\'e par une unit\'e de $R'$.
\end{prop}
\noindent\textit{Preuve du Th\'eor\`eme \ref{thmeclatform} en supposant la Proposition \ref{eclatpoly} vraie}: 
C'est la m\^eme que celle de la preuve du Th\'eor\`eme \ref{thmeclatformcar0}, il suffit de remplacer $k\left[\left[u_{1},...,u_{n-1}\right]\right]\left[u_n\right]$ par $R[r]\left[\left[u_{r+1},...,u_{n-1}\right]\right]\left[u_n\right]$.\\ \qed
\\ \\\noindent\textit{Preuve de la Proposition \ref{eclatpoly}}: 
C'est la m\^eme que celle de la preuve de la Proposition \ref{eclatpolycar0}, il suffit de remplacer $k\left[\left[u_{1},...,u_{n-1}\right]\right]\left[u_n\right]$ par $R[r]\left[\left[u_{r+1},...,u_{n-1}\right]\right]\left[u_n\right]$.
\\Pour conclure, il nous suffit d'avoir le r\'esultat voulu sur les polyn\^omes-cl\'es comme nous allons le voir dans la Proposition \ref{eclatpolycle}.
\begin{prop}\label{eclatpolycle}
Sous les hypoth\`eses du Th\'eor\`eme \ref{thmeclatform}, il existe une suite formelle encadr\'ee: \[(R,u)\rightarrow (R',u'),\] o\`u $u=(u_{1},...,u_{n})$, $u'=(u'_{1},...,u'_{n})$, v\'erifiant les propri\'et\'es suivantes:
\begin{enumerate}
\item Pour tout $q\in\mathbb{N}^*$ tel que $1\leqslant q \leqslant l_0$, les polyn\^omes-cl\'es $Q_{n,q}$ et $Q_{n,l}$ sont des mon\^omes en $u'$ multipli\'es par une unit\'e de $R'$;
\item Dans $R'$, $u'_{n}$ divise $Q_{n,l}$ mais $u'^{2}_{n}$ ne divise pas $Q_{n,l}$.
\end{enumerate}
\end{prop}
\noindent\textit{Preuve de la Proposition \ref{eclatpoly} en supposant la Proposition \ref{eclatpolycle} vraie}:
C'est la m\^eme que celle de la preuve de la Proposition \ref{eclatpolycar0}, il suffit de remplacer $k\left[\left[u_{1},...,u_{n-1}\right]\right]$ par $R[r]\left[\left[u_{r+1},...,u_{n-1}\right]\right]$.\\\qed
\\ \\
\noindent\textit{Preuve de la Proposition \ref{eclatpolycle}}: 
C'est la m\^eme que celle de la preuve de la Proposition \ref{eclatpolyclecar0}.\\\qed

\section{Un th\'eor\`eme de monomialisation sans compl\'etion}\label{soussectpasformel}
\indent Soient $\left( R,\mathfrak m,k\right)$ un anneau local r\'egulier complet de dimension $n$ tel que $\mathfrak m=(u)=(u_1,...,u_n)$. Soit $\nu$ une valuation de $K=Frac(R)$ centr\'ee en $R$ et de groupe des valeurs $\Gamma$. Notons $\Gamma_1$ le plus petit sous-groupe isol\'e non nul de $\Gamma$. On pose: 
\[H=\lbrace f\in R\:\vert\: \nu(f)\notin\Gamma_{1}\rbrace.\]
On suppose de plus que:
\[n=e(R,\nu)=emb.dim\left(R/H\right),\]
 c'est-\`a-dire que:
 \[H\subset\mathfrak{m}^2.\]
La valuation $\nu$ consid\'er\'ee est la compos\'ee de la valuation $\mu:L^{*}\rightarrow\Gamma_{1}$ de rang $1$ centr\'ee en $R/H$, o\`u $L=Frac(R/H)$, avec la valuation $\theta :K^{*}\rightarrow \Gamma / \Gamma_{1}$, centr\'ee en $R_{H}$, telle que $k_{\theta}\simeq \kappa(H)$. 
\\\indent Consid\'erons un sous-anneau local $\left(T,\mathfrak m_T\right)$ de $R$, non n\'ecessairement noeth\-\'erien, contenant $u_1,...,u_n$ et tel que $T/\mathfrak m_T\simeq k$.
Soient $J\subset\lbrace 1,...,n\rbrace$ et $j\in J$ tels que:
\[\nu(u_j)\leqslant \nu(u_i),\:i\in J.\]
Soit $\pi_0:(R,u)\rightarrow \left(R^{(1)},u^{(1)}\right)$ l'\'eclatement encadr\'e le long de $(u_J)$ par rapport \`a $\nu$ (D\'efinition \ref{eclatlocparrapnu}), notons $\mathfrak m^{(1)}$ l'id\'eal maximal de $R^{(1)}$.
\begin{defi}\label{defisuitedefsur}
Le \textbf{transform\'e de } $\boldsymbol T$ \textbf{par} $\boldsymbol{\pi_0}$ est l'anneau:
\[T^{(1)}=T\left[u'_{j\setminus\lbrace 0\rbrace}\right]_{\mathfrak m_1\cap T\left[u'_{j\setminus\lbrace 0\rbrace}\right]}.\]
On dit que l'\'eclatement $\pi_0$ est \textbf{d\'efini sur} $\boldsymbol T$ si $u^{(1)}\subset T^{(1)}$.
\\Pour une suite locale encadr\'ee de la forme:
\[ \xymatrix{\left( R,u\right)=\left( R^{(0)},u^{(0)}\right) \ar[r]^-{\pi_{0}} & \left( R^{(1)},u^{(1)}\right) \ar[r]^-{\pi_{1}} & \ldots   \ar[r]^-{\pi_{l-1}} & \left( R^{(l)},u^{(l)}\right)}, \]
les notions de \textbf{transform\'e} $T^{(i)}$ de $T$ et de \textbf{d\'efinie sur} $\boldsymbol T$ sont d\'efinies par r\'ecurrence sur $1\leqslant i \leqslant l$. Plus pr\'ecis\'ement, la notion de transform\'e $T^{(i)}$ n'est d\'efinie qu'en supposant la suite locale encadr\'ee $(R,u)\rightarrow \left(R^{(i-1)},u^{(i-1)}\right)$ d\'efinie sur T.
\end{defi}
Nous allons montrer qu'ind\'ependamment de la caract\'eristique de $k_\nu$, s'il existe un ensemble $1$-complet de polyn\^omes-cl\'es sans polyn\^ome-cl\'e limite pour un anneau local r\'egulier complet $(R,\mathfrak m,k)$, il va exister une suite locale encadr\'ee (et non plus formelle encadr\'ee) qui fasse d\'ecro\^itre l'invariant $e(R,\nu)$. Si $R$ est de caract\'eristique mixte, il faut supposer de plus que $\nu(p)\not\in p\Gamma$, $p=car(k_\nu)$.
\begin{thm}\label{monomialisationpasformel}
Reprenons les hypoth\`eses pr\'ec\'edentes et supposons de plus qu'il existe un ensemble $1$-complet de polyn\^omes-cl\'es ne poss\'edant pas de polyn\^ome-cl\'e limite. Si $R$ est de caract\'eristique mixte, on suppose \'egalement que $\nu(p)\not\in p\Gamma$, o\`u $p=car(k_\nu)$. Alors:
\begin{enumerate}
\item 
\begin{enumerate}
\item Ou bien $H\neq (0)$ et il existe une suite locale encadr\'ee $(R,u)\rightarrow (R',u')$ telle que: $$e(R',\nu)<e(R,\nu);$$
\item Ou bien $H=(0)$ et pour tout $f\in R$, il existe une suite locale encadr\'ee $(R,u)\rightarrow (R',u')$ telle que $f$ soit un mon\^ome en $u'$ multipli\'e par une unit\'e de $R'$.
\end{enumerate}
\item La suite locale encadr\'ee $(R,u)\rightarrow (R',u')$ de (1) peut \^etre choisie d\'efinie sur $T$.
\end{enumerate}
\end{thm}
\noindent\textit{Preuve}: Par les th\'eor\`emes \ref{thmeclatformcar0} et \ref{thmeclatform}, pour $f\in R$, il existe une suite formelle encadr\'ee:
\[ \xymatrix{\left( R,u\right)=\left( R^{(0)},u^{(0)}\right) \ar[r]^-{\pi_{0}} & \left( R^{(1)},u^{(1)}\right) \ar[r]^-{\pi_{1}} & \ldots   \ar[r]^-{\pi_{l-1}} & \left( R^{(l)},u^{(l)}\right)} \]
telle que, ou bien $e\left(R^{(l)},\nu\right)<e(R,\nu)$ si $H\neq (0)$, ou bien $f$ est un mon\^ome en $u^{(l)}$ multipli\'e par une unit\'e de $R^{(l)}$ si $H=(0)$. \`A partir de cette suite formelle encadr\'ee, nous allons construire, par approximation $\left(u^{(l)}\right)$-adique, la suite locale encadr\'ee $(R,u)\rightarrow (R',u')$ recherch\'ee. 
\\Plus pr\'ecis\'ement, pour $s\in\lbrace 1,...,l\rbrace$, consid\'erons $\pi_{s-1}:\left(R^{(s-1)},u^{(s-1)}\right)\rightarrow \left(R^{(s)},u^{(s)}\right)$ une des transformations de la suite formelle encadr\'ee, elle consiste en une suite \'el\'ementaire uniformisante $\pi_{0,s}$ (D\'efinition \ref{suiteelementaire}), qui r\'esout les singularit\'es d'un certain polyn\^ome-cl\'e, suivie d'une compl\'etion formelle. Ainsi, quitte \`a renum\'eroter les variables, $R^{(s)}$ est obtenu \`a partir de $R^{(s-1)}$ en adjoignant des expressions rationnelles $u_1^{(s)},...,u_r^{(s)},u_n^{(s)}$ en terme d'\'el\'ements de $R^{(s-1)}$ (dont les d\'enominateurs sont des mon\^omes en $u^{(s-1)}$), puis par passage au compl\'et\'e en le centre de la valuation $\nu$.
\\Pour $j\in\lbrace 1,...,n\rbrace$, notons $\mu_{j,s}$ la somme des valuations pour $\nu$ des num\'erateurs et d\'enominateurs de $u_j^{(s)}$, vu en tant que mon\^ome en $u^{(s-1)}$. On note alors:
\[\mu_s=\max_{1\leqslant j \leqslant n}\lbrace\mu_{j,s}\rbrace.\]
Soit $\beta\in\Gamma_1$ tel que $\beta>\sum\limits_{q=1}^l\mu_q$. Notons $I_s$ le $\nu_{0,u^{(s)}}$-id\'eal de $R^{(s)}$ d\'efini par:
\[I_s=\left\lbrace g\in R^{(s)}\:\left\vert\:\nu_{0,u^{(s)}}(g)\geqslant\beta-\sum\limits_{q=1}^s\mu_q\right.\right\rbrace.\]
Nous allons construire, par r\'ecurrence sur $s\in\lbrace 1,...,l\rbrace$, une suite locale encadr\'ee:
\[ \xymatrix{\left( R,u\right)=\left( \tilde R^{(0)},\tilde u^{(0)}\right) \ar[r]^-{\tilde \pi_{0}} & \left( \tilde R^{(1)},\tilde u^{(1)}\right) \ar[r]^-{\tilde \pi_{1}} & \ldots   \ar[r]^-{\tilde \pi_{l-1}} & \left( \tilde R^{(l)},\tilde u^{(l)}\right)}, \]
d\'efinie sur $T$ telle que, pour tout $s$ et $j\in\lbrace 1,...,n\rbrace$, on ait:
\[\textup{(HR)}:\:\:\nu_{0,u^{(s-1)}}\left(\tilde u_j^{(s)}-u_j^{(s)}\right)>\sum\limits_{q=s+1}^l\mu_q.\]
Supposons que la suite locale encadr\'ee soit construite \`a l'\'etape $s-1$. Quitte \`a renum\'eroter les variables si n\'ecessaire, on peut supposer que:
\[u_j^{(s)}=u_j^{(s-1)},\:r+1\leqslant j\leqslant n-1.\]
L'hypoth\`ese de r\'ecurrence: 
\[\nu_{0,u^{(s-2)}}\left(\tilde u_j^{(s-1)}-u_j^{(s-1)}\right)>\sum\limits_{q=s}^l\mu_q,\]
et le fait que les $u_1^{(s)},...,u_r^{(s)}$ s'expriment de mani\`ere rationnelle en fonction de $u^{(s-1)}$ entra\^ine que (HR) est vrai pour $j\in\lbrace 1,...,n-1\rbrace$.
\\Reprenons les notations de la sous-section \ref{sectionsuiteelemunif}. Consid\'erons:
\[\sum\limits_{i=1}^r\alpha_j\nu\left(u_j^{(s-1)}\right)=\overline\alpha\nu\left(u_n^{(s-1)}\right),\]
la plus petite combinaison $\mathbb Z$-lin\'eaire de $\nu\left(u_1^{(s-1)}\right),...,\nu\left(u_r^{(s-1)}\right),\nu\left(u_n^{(s-1)}\right)$ telle que $\overline\alpha\in\mathbb N^*$. Notons:
\[y=\left(u_1^{(s-1)}\right)^{\alpha_1}...\left(u_r^{(s-1)}\right)^{\alpha_r}\]
et
\[Q^{(s)}=\sum\limits_{i=0}^db_iy^{d-i}\left(u_n^{(s-1)}\right)^{i\overline\alpha},\]
le polyn\^ome $Q$ apparaissant dans la Proposition \ref{propcaspartpolycle} correspondant \`a la suite \'el\'ementaire uniformisante $\pi_{0,s}$. Pour chaque $b_i$ apparaissant dans $Q^{(s)}$, choisissons $\tilde b_i\in\left(\tilde R^{(s-1)}\right)^\times \cap T^{(s-1)}$ tel que:
\[\nu_{0,u^{(s-1)}}\left(b_j^{(s)}-\tilde b_j^{(s)}\right)>\sum\limits_{q=s+1}^l\mu_q.\]
Posons:
\[\tilde Q^{(s)}=\sum\limits_{i=0}^d\tilde b_iy^{d-i}\left(u_n^{(s-1)}\right)^{i\overline\alpha},\]
et $\tilde\pi_{s-1}$ la $n$-suite \'el\'ementaire uniformisante d\'etermin\'ee par ces donn\'ees. Ainsi, avec $Q^{(s)}$ et $\tilde Q^{(s)}$, on montre que (HR) est vraie pour $j=n$.
\\La suite locale encadr\'ee que l'on vient de construire par r\'ecurrence:
\[ \xymatrix{\left( R,u\right)=\left( \tilde R^{(0)},\tilde u^{(0)}\right) \ar[r]^-{\tilde \pi_{0}} & \left( \tilde R^{(1)},\tilde u^{(1)}\right) \ar[r]^-{\tilde \pi_{1}} & \ldots   \ar[r]^-{\tilde \pi_{l-1}} & \left( \tilde R^{(l)},\tilde u^{(l)}\right)} \]
d\'efinie sur $T$, telle que, pour tout $s$ et $j\in\lbrace 1,...,n\rbrace$:
\[\nu_{0,u^{(s-1)}}\left(\tilde u_j^{(s)}-u_j^{(s)}\right)>\sum\limits_{q=s+1}^l\mu_q,\]
implique que $in_\nu\left(u_j^{(s)}\right)=in_\nu\left(\tilde u_j^{(s)}\right)$, vus en tant qu'\'el\'ements de $\left(gr_\nu\left(R^{(l)}\right)\right)^*$, alg\`ebre qui contient la sous-alg\`ebre $gr_\nu\left(\tilde R^{(s)}\right)$. De m\^eme, la valuation monomiale $\nu_{0,u^{(l)}}$ de $Frac\left(R^{(l)}\right)$, restreinte \`a $\tilde R^{(l)}$, co\"incide avec la valuation monomiale $\nu_{0,\tilde u^{(l)}}$. On a alors $in_{\nu_{0,u^{(l)}}}\left(u_j^{(s)}\right)=in_{\nu_{0,u^{(l)}}}\left(\tilde u_j^{(s)}\right)$, vus en tant qu'\'el\'ements de $\left(gr_{\nu_{0,u^{(l)}}}\left(R^{(l)}\right)\right)^*$, alg\`ebre qui contient la sous-alg\`ebre $gr_{\nu_{0,u^{(l)}}}\left(\tilde R^{(s)}\right)$.
\\Si $H=(0)$ et $f\in R\setminus\lbrace 0\rbrace$ est monomialis\'e par la suite formelle encadr\'ee, l'\'egalit\'e pr\'ec\'edente implique que:
\[f=\varpi + \tilde f,\]
o\`u $\varpi$ est un mon\^ome en $\tilde u^{(l)}$ et $\nu_{0,\tilde u^{(l)}}(\tilde f)>\nu(\varpi)$.
\\Si $H\neq (0)$ et $f\in H\setminus\lbrace 0\rbrace$ dont le transform\'e strict devient un param\`etre r\'egulier dans $R^{(l)}$, alors:
\[f=\tilde Q^{(l)}+\tilde f,\]
o\`u $\nu_{0,\tilde u^{(l)}}(\tilde f)>\nu(\tilde Q^{(l)})$. En appliquant le Corollaire \ref{coroideal}, apr\`es une suite monomiale $\left(\tilde R^{(l)},\tilde u^{(l)}\right)\rightarrow (R',u')$ (respectivement, en appliquant la Proposition \ref{proppolycleplus}, apr\`es une suite locale encadr\'ee ind\'ependante de $\tilde u_n^{(l)}$, dans le cas $H\neq (0)$), on est ramen\'e \`a la situation o\`u $\varpi$ divise $\tilde f$, c'est-\`a-dire \`a la situation o\`u $f=z\varpi$, $z$ unit\'e de $R'$ (respectivement, $f=zg$, o\`u $g$ est un param\`etre r\'egulier de $R'$ et $z$ une unit\'e de $R'$, dans le cas $H\neq (0)$).\\\qed
\section{Th\'eor\`emes d'uniformisation locale en caract\'eristique nulle}
Soit $S$ un anneau local noeth\'erien. Pour montrer que $S$ est transform\'e en un anneau r\'egulier via une suite locale encadr\'ee, il faut montrer que $\widehat S_{\overline H}$ et $\widehat S/\overline H$ le sont, $\overline H$ \'etant l'id\'eal premier implicite de $\widehat S$. Par le Th\'eor\`eme \ref{RHestreg}, si $S$ est quasi-excellent alors $\widehat S_{\overline H}$ est r\'egulier. Dans un premier temps, nous allons montrer que, sous certaines hypoth\`eses, $\widehat S/\overline H$ est aussi r\'egulier. Enfin, gr\^ace \`a ces deux r\'esultats nous montrerons le th\'eor\`eme d'uniformisation locale pour des valuations de rang $1$ puis pour des valuations de rang quelconque gr\^ace \`a \cite{spivanova}.
\subsection{Un th\'eor\`eme pr\'eliminaire d'uniformisation locale}
\begin{thm}\label{thmprelimcar0}
Soient $(S,\mathfrak m,k)$ un anneau local noeth\'erien int\`egre de corps des fractions $L$ et $\mu$ une valuation de $L$ de rang $1$ et de groupe des valeurs $\Gamma_1$, centr\'ee en $S$, telle que $car\left(k_\mu\right)=0$.
\\Notons $u=(u_1,...,u_n)$ un ensemble minimal de g\'en\'erateurs de $\mathfrak m$ et $\overline H$ l'id\'eal premier implicite de $\widehat S$.
\\Soient $f_1,...,f_s\in\mathfrak m$ tels que $\mu(f_1)=\min_{1\leqslant j\leqslant s}\lbrace\mu(f_j)\rbrace$. Il existe alors une suite locale encadr\'ee:
\[ \xymatrix{\left( S,u,k\right)=\left( S^{(0)},u^{(0)},k^{(0)}\right) \ar[r]^-{\rho_{0}} & \left( S^{(1)},u^{(1)},k^{(1)}\right) \ar[r]^-{\rho_{1}} & \ldots   \ar[r]^-{\rho_{i-1}} & \left( S^{(i)},u^{(i)},k^{(i)}\right)}, \]
ayant les propri\'et\'es suivantes:
\\Notons $\overline H^{(i)}$ l'id\'eal premier implicite de $\widehat{S^{(i)}}$ et $\overline{f_j}$ l'image de $f_j\mod\overline H^{(i)}$, $1\leqslant j\leqslant s$, alors:
\begin{enumerate}
\item $\widehat{S^{(i)}}/\overline H^{(i)}$ est r\'egulier;
\item Pour $1\leqslant j\leqslant s$, $\overline{f_j}$ est un mon\^ome en $u^{(i)}$ multipli\'e par une unit\'e de $\widehat{S^{(i)}}/\overline H^{(i)}$;
\item Pour $1\leqslant j\leqslant s$, $\overline{f_1}$ divise $\overline{f_j}$ dans $\widehat{S^{(i)}}/\overline H^{(i)}$.
\end{enumerate}
\end{thm}
\noindent\textit{Preuve}: Notons $\sigma:S\rightarrow \widehat S$ le morphisme de compl\'etion formelle. Par le Th\'eor\`eme \ref{valetenuniqidealimpl}, $\mu$ s'\'etend de mani\`ere unique en une valuation $\widehat\mu$ centr\'ee en $\widehat S/\overline H$. Notons $u=(y,x)$ tel que $x=(x_1,...,x_l)$, $l=emb.dim\left(\widehat{S}/\overline{H}\right)$, $y=(y_1,...,y_{n-l})$ et les images des $x_1,...,x_l$ dans $\widehat S/\overline H$ induisent un ensemble minimal de g\'en\'erateurs de $(\mathfrak m \widehat S)/\overline H$.
\\Par le Th\'eor\`eme de structure de Cohen, on sait qu'il existe un anneau local r\'egulier complet de caract\'eristique nulle $R$ et un morphisme $\varphi$ surjectif:
\[\varphi:R\twoheadrightarrow \widehat S/\overline H.\]
Notons $H=\ker \varphi$, comme $\overline H$ est un id\'eal premier (Th\'eor\`eme \ref{valetenuniqidealimpl}), $H$ est un id\'eal premier de $R$. On choisit $R$ de telle sorte que $\dim (R)=l$. Notons $K$ le corps des fractions de $R$. Soit $\theta$ une valuation de $K$, centr\'ee en $R_H$, telle que $k_\theta=\kappa(H)$. Si l'on regarde $\widehat\mu$ comme une valuation centr\'ee en $R/H$ via le morphisme $\varphi$, on peut consid\'erer la valuation $\nu=\widehat{\mu}\circ\theta$ centr\'ee en $R$ et de groupe des valeurs $\Gamma$. Alors, $\Gamma_1$ est le plus petit sous-groupe isol\'e non nul de $\Gamma$ et:
\[H=\lbrace f\in R\:\vert\:\nu(f)\not\in\Gamma_1\rbrace.\]
De plus, $car\left(k_\nu\right)=car\left(k_\mu\right)=0$. On s'est donc ramen\'e aux hypoth\`eses du Th\'eor\`eme \ref{thmeclatformcar0}.
\\Soit $T=\varphi^{-1}(\sigma(S))$, c'est un sous-anneau local de $R$ d'id\'eal maximal $\varphi^{-1}(\sigma(\mathfrak m))=\mathfrak m\cap T$. Ainsi, $T$ contient $x_1,...,x_l$ et:
\[T/(\mathfrak m\cap T)\simeq k.\]
Comme le Th\'eor\`eme \ref{thmeclatformcar0} est vrai en caract\'eristique $0$, on peut appliquer le Th\'eor\`eme \ref{monomialisationpasformel}. Ainsi, plusieurs cas se pr\'esentent:
\begin{enumerate}
\item Si $H\neq (0)$, il existe une suite locale encadr\'ee $(R,x)\rightarrow \left(R^{(i)},x^{(i)}\right)$ telle que $e(R,\nu)$ d\'ecroisse strictement. En particulier, ce cas ne peut arriver qu'un nombre fini de fois, on arrive ainsi \`a la situation o\`u $H=(0)$ et donc $R/H$ est r\'egulier.
\item Si $H=(0)$, alors, pour chaque $f_j$, $1 \leqslant j\leqslant s$, il existe une suite locale encadr\'ee $(R,x)\rightarrow \left(R^{(i)},x^{(i)}\right)$ telle que $f_j$ soit un mon\^ome en $x^{(i)}$ multipli\'e par une unit\'e de $R^{(i)}$.
\end{enumerate}
Par la Proposition \ref{propgenerem}, la propri\'et\'e d'\^etre un mon\^ome multipli\'e par une unit\'e est pr\'eserv\'ee par les suites locales encadr\'ees. Ainsi, en it\'erant la proc\'edure de (2), on arrive \`a la situation o\`u tous les $f_1,...,f_s$ sont simultan\'ement des mon\^omes en $x^{(i)}$. Apr\`es une suite locale encadr\'ee de plus $(R,x)\rightarrow \left(R',x'\right)$, on peut supposer que les $f_j$ sont des mon\^omes uniquement en $x'_1,...,x'_r$, $1\leqslant j\leqslant s$, $r=r(R,x,\nu)$. Enfin, en appliquant plusieurs fois le Corollaire \ref{coroexistsuitlocale}, on est ramen\'e \`a la situation o\`u chaque $f_j$ est un mon\^ome en $x'_1,...,x'_r$, $1\leqslant j\leqslant s$ et, pour $j,j'\in\lbrace 1,...,s\rbrace$, $f_j$ divise $f_{j'}$ ou $f_{j'}$ divise $f_j$. De plus, toutes ces suites locales encadr\'ees sont d\'efinies sur $T$.
Consid\'erons le diagramme suivant:
\[ \xymatrix{\left( R,x,k\right) \ar[r]^-{\pi_{0}} \ar[d] & \left( R^{(1)},x^{(1)},k^{(1)}\right) \ar[r]^-{\pi_{1}} \ar[d] & \ldots   \ar[r]^-{\pi_{i-1}} & \left( R^{(i)},x^{(i)},k^{(i)}\right)\ar[d] \\ \left(\widehat{S}/\overline{H},x,k\right)\ar[r]^-{\tilde\pi_{0}}  & \left( \tilde S^{(1)},x^{(1)},k^{(1)}\right) \ar[r]^-{\tilde\pi_{1}} & \ldots   \ar[r]^-{\tilde\pi_{i-1}} & \left( \tilde S^{(i)},x^{(i)},k^{(i)}\right)\\ \left( S,u,k\right) \ar[r]^-{\rho_{0}} \ar[u] & \left( S^{(1)},u^{(1)},k^{(1)}\right) \ar[r]^-{\rho_{1}} \ar[u] & \ldots   \ar[r]^-{\rho_{i-1}} & \left( S^{(i)},u^{(i)},k^{(i)}\right)\ar[u] } \]
Par ce que l'on vient de voir, la premi\`ere colonne et la premi\`ere ligne ont d\'ej\`a \'et\'e construites. En passant au transform\'e strict de $R/H\simeq\widehat S/\overline H$ \`a chaque \'etape de la suite $(\pi_j)_{1\leqslant j \leqslant i-1}$, on construit la suite d'\'eclatements encadr\'es $\left(\tilde\pi_j\right)_{1\leqslant j \leqslant i-1}$ de $\widehat S/\overline H$ d\'efinie sur $S$. Enfin, la suite $\left(\tilde\pi_j\right)_{1\leqslant j \leqslant i-1}$ se rel\`eve en une suite locale encadr\'ee $(\rho_j)_{1\leqslant j \leqslant i-1}$.
\\ Si $R/H$ est singulier, par le Th\'eor\`eme \ref{monomialisationpasformel}, il existe une suite locale encadr\'ee $(\pi_j)_{1\leqslant j \leqslant i-1}$ qui fasse d\'ecro\^itre $e(R,\nu)$. Ainsi, la suite locale encadr\'ee $(\rho_j)_{1\leqslant j \leqslant i-1}$ r\'esultante poss\`ede la propri\'et\'e:
\[emb.dim\left(\widehat{S^{(i)}}/\overline{H}^{(i)}\right)<emb.dim\left(\widehat{S}/\overline{H}\right).\]
Ceci n'arrive qu'un nombre fini de fois. Apr\`es un nombre fini de pas, on arrive \`a la situation o\`u $\widehat{S^{(i)}}/\overline H^{(i)}$ est r\'egulier. Maintenant, si l'on suppose que $\widehat{S^{(i)}}/\overline H^{(i)}$ est r\'egulier, consid\'erons $f_1,...,f_s$ des \'el\'ements non nuls de $S$ tels que $\mu(f_1)=\min_{1\leqslant j\leqslant s}\lbrace\mu(f_j)\rbrace$, alors, par le (2) vu plus haut, on en d\'eduit que, pour $1\leqslant j \leqslant s$, $f_j\mod\overline H^{(i)}$ sont des mon\^omes en $u^{(i)}$ et $f_1\mod\overline H^{(i)}$ divise $f_j\mod\overline H^{(i)}$.\\\qed
\subsection{Uniformisation locale plong\'ee pour des valuations de rang 1}
~\smallskip ~\\ \indent Avant d'\'enoncer et de d\'emontrer le th\'eor\`eme d'uniformisation locale plong\'ee pour des valuations de rang $1$, nous allons donner un lemme un peu plus g\'en\'eral et ind\'ependant de la caract\'eristique. 
\begin{lem}\label{lemmetechniquespiva}(\cite{spiva}, Lemme 16.3) Soient $(A,\mathfrak m,k)$ un anneau local noeth\'erien, $\nu$ une valuation centr\'ee en $A$ et $J$ un $\nu$-id\'eal premier de $A$ non maximal. Notons $h=ht(J)$. Supposons que $A_J$ et $A/J$ soient r\'eguliers. Notons $u=(u_1,...,u_n)$ un ensemble minimal de g\'en\'erateurs de $\mathfrak m$ et supposons que $u=(x,y)$ avec $x=(x_1,...,x_l)$ et $y=(y_1,...,y_{n-l})$ tels que:
\begin{enumerate}
\item $x$ induit un syst\`eme r\'egulier de param\`etres de $A/J$;
\item il existe un ensemble minimal de g\'en\'erateurs $(\widehat{y}_1,...,\widehat{y}_{n-l})$ de $J$ et des mon\^omes $\varpi_1,...,\varpi_{n-l}$ en $x$ tels que $\varpi_1/.../\varpi_{n-l}$ de sorte que $(\widehat{y}_{n-l-h+1},...,\widehat{y}_{n-l})$ induit un syst\`eme r\'egulier de param\`etres de $A_J$ et, pour tout $N\in\mathbb N^*$, il existe $v_j\in A^\times$ tel que:
\[\widehat{y}_j-y_j-\varpi_jv_j\in\varpi_j\mathfrak m^N,\]
$1\leqslant j\leqslant n-l$. Remarquons que, par convention, on peut avoir $y_j=\widehat y_j$, $\varpi_j=0$, $v_j=1$ et $(y)=J$.
\end{enumerate}
Soient $f_1,...,f_s\in A$ tels que:
\[\nu(f_1)\leqslant...\leqslant \nu(f_s).\]
Soit $(T,\mathfrak m_T)$ un sous-anneau local de $A$ non n\'ecessairement noeth\'erien tel que $T/\mathfrak m_T=k$. Enfin, supposons que pour tout $g_1,...,g_t\in A$ tels que:
\[\nu(g_1)\leqslant...\leqslant \nu(g_t),\]
il existe une suite locale encadr\'ee $(A,u)\rightarrow (A',u')$ ind\'ependante de $y$ et d\'efinie sur $T$ telle que, pour tout $1\leqslant j \leqslant t$, $g_j\:mod\: J'$ soit un mon\^ome en $u'$ et $g_q\:mod\: J'$ divise $g_i\:mod\: J'$, $1\leqslant q\leqslant i\leqslant t$, o\`u $J'$ est le transform\'e strict de $J$ dans $A'$.
\\Il existe alors une suite locale encadr\'ee $(A,u)\rightarrow (A'',u'')$ par rapport \`a $\nu$ et d\'efinie sur $T$ telle que $A''$ soit r\'egulier.
\\Supposons de plus que l'une au moins des deux conditions suivantes soit v\'erifi\'ee:
\begin{enumerate}
\item[(3)] $f_i\notin J$, $1\leqslant i \leqslant s$;
\item[(4)] $y_j=\widehat{y}_j$, $1\leqslant j \leqslant n-l$ (donc $J=(y)$), $T=A$ et, pour tout $1\leqslant i \leqslant s$, $f_i$ est un mon\^ome en $(y_{n-l-h+1},...,y_{n-l})$ et $f_i/f_{i+1}$ dans $A_J$.
\end{enumerate}
La suite locale encadr\'ee $(A,u)\rightarrow (A'',u'')$ pr\'ec\'edente peut alors \^etre choisie de telle sorte que les $f_i$ soient des mon\^omes en $u''$ multipli\'es par une unit\'e de $A''$ et telle que $f_i/f_{i+1}$ dans $A''$, $1\leqslant i \leqslant s$.
\end{lem}
\noindent\textit{Preuve}: Nous ne donnerons qu'une id\'ee de preuve, pour plus de d\'etails, on peut consulter \cite{spiva}. Si $J=(0)$, il n'y a rien \`a montrer; supposons donc que $J\neq (0)$. \`A partir de la suite locale encadr\'ee $(A,u)\rightarrow (A',u')$, on veut construire une suite locale encadr\'ee $(A,u)\rightarrow (A'',u'')$ d\'efinie sur $T$ telle que $A''$ soit r\'egulier. Pour cela il suffit d'avoir:
\[A''=Fitt_h(J''/J''^2),\]
o\`u $Fitt_h(J''/J''^2)$ est le $h$-i\`eme id\'eal de Fitting de $J''/J''^2$. Par hypoth\`ese et apr\`es une suite locale encadr\'ee n'impliquant que des variables en $x$, on peut se ramener \`a la situation o\`u $Fitt_h(J'/J'^2)$ est principal et engendr\'e par un mon\^ome en $x$, not\'e $a$. Quitte \`a renum\'eroter les variables de $y$, on peut supposer qu'il existe $n-l-h$ relations de la forme:
\[\psi_j=a\widehat y_j+\sum\limits_{q=n-l-h+1}^{n-l}a_{j,q}\widehat y_q+g_j,\]
o\`u $g_j\in J'^2$ et $a$ divise $a_{j,q}$ pour $1\leqslant j\leqslant n-l-h$ et $n-l-h+1\leqslant q\leqslant n-l$. Si on a (4), alors:
\[\nu_{0,u}(y_j)>\nu(a),\:1\leqslant j\leqslant n-l.\]
Comme $J$ est un $\nu$-id\'eal alors $y_j\in J$ et $a\notin J$. Supposons que l'on n'ait pas (4) et prenons $N\in\mathbb N^*$ tel que:
\[N>\dfrac{1}{\nu_{0,u}(\mathfrak m)}\left(\nu(\varpi_{n-l})+\max\limits_{\underset{f_q\notin J}{1\leqslant q\leqslant s}}\lbrace \nu(f),\nu(a)\rbrace\right).\]
Consid\'erons une variable $x_j$ de $x$ telle que $x_j^\alpha$ divise $\varpi_1$ pour un certain $\alpha\in\mathbb N^*$. On \'eclate l'id\'eal $(y,x_j)$ et on r\'ep\`ete cette proc\'edure $\alpha$ fois. On fait de m\^eme pour toutes les autres variables divisant $\varpi_1$. On arrive alors \`a la situation o\`u:
\[\nu_{0,u'}(y'_1)>\nu(a)+\nu(\varpi_{n-l})-\nu(\varpi_1).\]
On refait pareil pour toutes les autres variables de $y$; on est ainsi ramen\'e \`a la situation o\`u, pour ces nouvelles variables:
\[\nu_{0,u}(y'_j)>\nu(a),\:1\leqslant j\leqslant n-l.\]
Pour chaque variable $x_j$ de $x$ divisant $a$, on \'eclate en l'id\'eal $(y,x_j)$. Ces \'eclatements ont pour effet de multiplier $a$ et les $a_{j,q}$ par $x_j$ ainsi que les $g_1,...,g_t$ par $x_j^\gamma$ o\`u $\gamma\geqslant 2$. Apr\`es un nombre fini de fois, $a$ divise $g_j$ et donc $a$ divise $\psi_j$, $1\leqslant j \leqslant n-l-h$. Ainsi, pour $1\leqslant j \leqslant n-l-h$, les $y_j$ s'expriment comme une fonction des variables restantes modulo $\mathfrak m^2$. Ceci fait donc d\'ecro\^itre $emb.dim(A)$ et $A$ est r\'egulier, \`a une suite formelle encadr\'ee pr\`es.
\\\`A partir de maintenant on peut supposer que $h=n-l$; pour terminer il faut montrer que les $f_i$ sont des mon\^omes en $u''$ multipli\'es par une unit\'e de $A''$. Quitte \`a diviser $f_i$ par un mon\^ome en $y$, on peut supposer que (3) est toujours v\'erifi\'ee. Si (4) est v\'erifi\'ee alors:
\[\nu_{0,u}(y'_j)>\nu(f_i),\:1\leqslant j\leqslant n-l,\:1\leqslant i \leqslant s.\]
Si (4) n'est pas v\'erifi\'ee, l'in\'egalit\'e pr\'ec\'edente reste vraie par le choix de $N$. Ainsi, pour $1\leqslant i\leqslant s$, on a:
\[f_i=\rho_i+\tilde f_i,\]
o\`u $\rho_i$ est un mon\^ome en $x$ et $\nu_{0,u}(\tilde f_i)>\nu_{0,u}(\rho_i)$. On applique le Corollaire \ref{coroideal} \`a chaque $f_i$, $1\leqslant i\leqslant s$, et on obtient le r\'esultat cherch\'e.\\\qed
\\\indent Passons maintenant au th\'eor\`eme d'uniformisation locale plong\'ee pour des valuations de rang $1$ sur un anneau \'equicaract\'eristique dont le corps r\'esiduel est de caract\'eristique nulle.
\begin{thm}\label{uniflocalerang1car0}
Soient $(S,\mathfrak m,k)$ un anneau local int\`egre quasi-excellent de corps des fractions $L$ et $\mu$ une valuation de $L$ de rang $1$ et de groupe des valeurs $\Gamma_1$, centr\'ee en $S$, telle que $car\left(k_\mu\right)=0$.
\\Notons $u=(u_1,...,u_n)$ un ensemble minimal de g\'en\'erateurs de $\mathfrak m$.
\\Soient $f_1,...,f_s\in\mathfrak m$ tels que $\mu(f_1)=\min_{1\leqslant j\leqslant s}\lbrace\mu(f_j)\rbrace$. Il existe alors une suite locale encadr\'ee:
\[ \xymatrix{\left( S,u,k\right)=\left( S^{(0)},u^{(0)},k^{(0)}\right) \ar[r]^-{\rho_{0}} & \left( S^{(1)},u^{(1)},k^{(1)}\right) \ar[r]^-{\rho_{1}} & \ldots   \ar[r]^-{\rho_{i-1}} & \left( S^{(i)},u^{(i)},k^{(i)}\right)}, \]
ayant les propri\'et\'es suivantes:
\begin{enumerate}
\item $S^{(i)}$ est r\'egulier;
\item Pour $1\leqslant j\leqslant s$, $f_j$ est un mon\^ome en $u^{(i)}$ multipli\'e par une unit\'e de $S^{(i)}$;
\item Pour $1\leqslant j\leqslant s$, $f_1$ divise $f_j$ dans $S^{(i)}$.
\end{enumerate}
En d'autres termes, $\mu$ admet une uniformisation locale plong\'ee au sens de la Propri\'et\'e \ref{uniflocplongint}.
\end{thm}
\noindent\textit{Preuve}: Reprenons les notations du Th\'eor\`eme \ref{thmprelimcar0}. On a vu qu'il existe un morphisme surjectif:
\[\psi:\widehat S\twoheadrightarrow \widehat{S}/\overline{H}\simeq R/H.\]
Par le Th\'eor\`eme \ref{thmprelimcar0}, apr\`es une suite locale encadr\'ee auxiliaire, on peut supposer que $\widehat{S}/\overline{H}$ est r\'egulier et donc que $R/H\simeq k\left[\left[x_{1},...,x_l\right]\right]$. Ainsi, il existe un ensemble de g\'en\'erateurs $\widehat y=\left(\widehat y_1,...,\widehat y_{n-l}\right)$ de $\overline H$ et des s\'eries formelles $\phi_j\in k\left[\left[x_{1},...,x_l\right]\right]$ tels que:
\[\widehat y_j=y_j+\phi_j \in \widehat S,\:1\leqslant j\leqslant n-l.\]
Quitte \`a renum\'eroter les $y_j$, on peut supposer que:
\[\mu(y_1)\leqslant\mu(y_2)\leqslant...\leqslant\mu(y_{n-l}).\]
En appliquant le Corollaire \ref{coroideal} aux mon\^omes de $\phi_j$, $1\leqslant j\leqslant n-l$, on peut supposer que:
\[\phi_j=\varpi_j\widehat v_j,\]
o\`u les $\varpi_j$ sont des mon\^omes en $x_1,...,x_l$, $\widehat v_j\in k\left[\left[x_{1},...,x_l\right]\right]^\times$ et tels que:
\[\varpi_1/.../\varpi_{n-l}.\]
Ainsi, on en d\'eduit que:
\[\forall\:j\in\lbrace 1,...,n-l\rbrace,\:\forall\:N\in\mathbb N^*,\:\exists\:v_j\in S^\times,\:\widehat y_j-y_j-\varpi_jv_j\in\varpi_j\mathfrak m^N.\]
Enfin, rappelons que, par le Corollaire \ref{RHestreg}, l'anneau $\widehat S_{\overline H}$ est r\'egulier. On applique alors le Lemme \ref{lemmetechniquespiva} \`a $A=\widehat S$, $J=\overline H$, $T=S$ et $\nu=\mu$. On en d\'eduit ainsi une uniformisation locale plong\'ee (Propri\'et\'e \ref{uniflocplongint}) de $\widehat S$. Or $S$ est quasi-excellent, donc le morphisme $S\rightarrow \widehat{S}$ est r\'egulier et on obtient ainsi une uniformisation locale plong\'ee (Propri\'et\'e \ref{uniflocplongint}) de $S$.\\\qed
\subsection{Th\'eor\`emes d'uniformisation locale plong\'ee}
\begin{coro}\label{uniflocaleplongeecar0}
Soient $(S,\mathfrak m,k)$ un anneau local int\`egre quasi-excellent de corps des fractions $L$ et $\nu$ une valuation de $L$ centr\'ee en $S$ et de groupe des valeurs $\Gamma$ telle que $car\left(k_\nu\right)=0$.
\\Alors, $\nu$ admet une uniformisation locale plong\'ee au sens de la Propri\'et\'e \ref{uniflocplongint}.
\end{coro}
\noindent\textit{Preuve}: On applique le Th\'eor\`eme \ref{uniflocalerang1car0} et le Th\'eor\`eme 1.3 de \cite{spivanova}.\\\qed
\begin{coro}\label{thm1.6car0}
Soient $(S,\mathfrak m,k)$ un anneau local int\`egre quasi-excellent de corps des fractions $L$ et $\nu$ une valuation de $L$ centr\'ee en $S$ et de groupe des valeurs $\Gamma$ telle que $car\left(k_\nu\right)=0$.
\\Pour $I$ un id\'eal de $S$, la paire $(S,I)$ admet une uniformisation locale plong\'ee par rapport \`a $\nu$ au sens de la D\'efinition \ref{defiuniflocpourpaire}.
\end{coro}
\noindent\textit{Preuve}: C'est une application imm\'ediate du Corollaire \ref{uniflocaleplongeecar0}.\\\qed
\begin{coro}Soit $S$ un sch\'ema quasi-excellent tel que pour tout $\xi\in S$, $car(\mathcal O_{S,\xi}/\mathfrak m_{S,\xi})=0$. Soient $X$ une composante irr\'eductible de $S_{red}$ et $\nu$ une valuation de $K(X)$
centr\'ee en un point $\xi\in X$. Il existe alors un \'eclatement $\pi: S'\rightarrow S$ le long d'un sous-sch\'ema de $S$, ne contenant aucune composante irr\'eductible de $S_{red}$ et ayant la propri\'et\'e suivante : \\Soient $X'$ le transform\'e strict de $X$ par $\pi$, $\xi'$ le centre de $\nu$ sur $X'$ et $D$ le diviseur exceptionnel de $\pi$, alors $(\mathcal{O}_{X',\xi'},\mathcal{I}_{D,\xi'})$ admet une uniformisation locale plong\'ee par rapport \`a $\nu$ au sens de la D\'efinition \ref{defiuniflocpourpaire}.
\end{coro}
\noindent\textit{Preuve}: C'est une application directe du Corollaire \ref{thm1.6car0}.\\\qed
\begin{thm}\label{thmfinal}
Soit $(S,\mathfrak m,k)$ un anneau local (non n\'ecessairem\-ent int\`egre) quasi-excellent. Soient $P$ un id\'eal premier minimal de $S$ et $\nu$ une valuation du corps des fractions de $S/P$ centr\'ee en $S/P$ et de groupe des valeurs $\Gamma$ telle que $car\left(k_\nu\right)=0$.
\\Il existe alors un \'eclatement local $\pi:S\rightarrow S'$ par rapport \`a $\nu$ tel que $S'_{red}$ soit r\'egulier et $Spec(S')$ soit normalement plat le long de $Spec(S'_{red})$, c'est-\`a-dire que l'anneau $S$ admet une uniformisation locale par rapport \`a $\nu$ au sens de la Propri\'et\'e \ref{uniflocanneaunonint}.
\end{thm}
\noindent\textit{Preuve}: 
Par le Corollaire \ref{thm1.6car0}, il existe une suite locale encadr\'ee $(S,u)\rightarrow (S',u')$ le long de centres ne contenant aucune composante irr\'eductible du transform\'e strict de $Spec\left(S_{red}\right)$, telle que $Spec\left(S'_{red}\right)$ soit r\'egulier. On peut donc supposer que $S_{red}$ est r\'egulier. Il reste \`a montrer qu'il existe une suite locale encadr\'ee telle que $Spec(S')$ soit normalement plat le long de $Spec(S'_{red})$.
\\Soit $(y_1,...,y_h)=\sqrt{(0)}\subset S$, c'est l'id\'eal qui d\'efinit $ Spec\left(S_{red}\right)$ dans $Spec(S)$. 
\\Rappelons que pour un anneau local noeth\'erien $(R,\mathfrak n)$, le \textit{c\^one tangent} de $Spec(R)$ est d\'efini par:
\[Spec\left(\bigoplus\limits_{n\geqslant 0}\mathfrak n^{n}/\mathfrak n^{n+1}\right).\]
Il suffit de construire une suite locale encadr\'ee telle que le c\^one tangent de $Spec(S')$ soit d\'efini par un id\'eal engendr\'e par des \'el\'ements de $k\left[\overline{y'_1},...,\overline{y'_h}\right]$, o\`u $y'_j$ est le transform\'e strict de $y_j$ dans $S'$ et $\overline{y_j}$ est l'image naturelle de $y_j$ dans l'alg\`ebre gradu\'ee de $S'$, $1\leqslant j\leqslant h$.
\\Notons $A=S_{red}$, $A'=S'_{red}$, on peut alors \'ecrire $S$ sous-la forme:
\[S=A\left[y_1,...,y_h\right]/I.\]
Notons $f_1,...,f_s\in A\left[y_1,...,y_h\right]$ un ensemble de g\'en\'erateurs de $I$ et $\left(x_1,...,x_r\right)$ un ensemble minimal de g\'en\'erateurs de l'id\'eal maximal de $A$.
\\Pour $1\leqslant j\leqslant s$, notons $f_j=\sum\limits_\alpha c_{j,\alpha}y^\alpha\in A[y]$. On va construire une suite locale encadr\'ee et une partition $(u')=(y',x')$ de $(u')$ o\`u $(y')$  est le transform\'e strict de $(y)$.
\\Soit $\nu_{0,x'}$ la valuation monomiale de $A'$ associ\'ee \`a $x'$ et \`a $\left\lbrace\nu\left(x'_j\right)\right\rbrace_j$ (Corollaire \ref{defivaluationmono}). Par le Corollaire \ref{thm1.6car0}, on peut construire une suite locale encadr\'ee $(S,u)\rightarrow (S',u')$ telle que les $c_{j,\alpha}$ soient des mon\^omes en $x'$ multipli\'es par une unit\'e de $A'$.
\\Pour tout $j\in\lbrace 1,...,s\rbrace$, notons $\mu_j=\max\lbrace N\in\mathbb N^*\:\vert\:f_j\in(y)^N\rbrace$ et $f'_j=\sum\limits_\alpha c'_{j,\alpha}y'^\alpha\in A'[y']$ le transform\'e strict de $f_j$ dans $S'=A'[y']$. Pour chaque $x'_t$ apparaissant dans $c'_{j,\alpha}$, pour un certain $j$ et un certain $\alpha$ tel que $\vert\alpha\vert=\mu_j$, \'eclatons en l'id\'eal $(y'_1,...,y'_h,x'_t)$ un nombre suffisant de fois. On arr\^ete le processus lorsque, pour $1\leqslant j\leqslant s$ et $\alpha$ tel que $\vert\alpha\vert>\mu_j$, il existe $\tilde\alpha$ tel que $c'_{j,\tilde\alpha}$ divise $c'_{j,\alpha}$, avec $\vert\tilde\alpha\vert=\mu_j$. Par le Corollaire \ref{coroideal}, on sait que, pour chaque $j$, il existe bien $\tilde\alpha$ tel que $\vert\tilde\alpha\vert=\mu_j$ et pour tout $\alpha$, $c'_{j,\tilde\alpha}$ divise $c'_{j,\alpha}$. Ainsi, le c\^one tangent de $Spec(S')$ est d\'efini par des polyn\^omes qui ne d\'ependent que de $\overline{y'_1},...,\overline{y'_h}$. On en conclut que $Spec(S')$ est normalement plat le long de $Spec(S'_{red})$.\\\qed
\begin{coro}
Soit $S$ un sch\'ema quasi-excellent tel que pour tout $\xi\in S$, $car(\mathcal O_{S,\xi}/\mathfrak m_{S,\xi})=0$. Soient $X$ une composante irr\'eductible de $S_{red}$ et $\nu$ une valuation de $K(X)$
centr\'ee en un point $\xi\in X$. Il existe alors un \'eclatement $\pi: S'\rightarrow S$ le long d'un sous-sch\'ema de $S$, ne contenant aucune composante irr\'eductible de $S_{red}$ et ayant la propri\'et\'e suivante : \\Soient $X'$ le transform\'e strict de $X$ par $\pi$ et $\xi'$ le centre de $\nu$ sur $X'$, alors $\xi'$ est un point r\'egulier de $X'$ et $S'$ est normalement plat le long de $X'$ en $\xi'$.
\end{coro}
\bibliographystyle{plain}
\bibliography{biblio}
\end{document}